\documentclass[preprint, 5p]{elsarticle}
\title{Analysis of a Sinclair-type domain decomposition solver\\for atomistic/continuum coupling}
\author{M.~Hodapp}
\address{Skolkovo Institute of Science and Technology (Skoltech), Center for Energy Science and Technology, Moscow (RU)}
\address{Ecole polytechnique f\'{e}d\'{e}rale de Lausanne (EPFL), STI IGM LAMMM, Lausanne (CH)}

\journal{SIAM Multiscale Modeling and Simulation}
\usepackage{amssymb}
\usepackage{mathrsfs}
\usepackage{amsthm}
\usepackage[nonumberlist,nopostdot,sort=def,nowarn,acronym]{glossaries} 
\usepackage{glossary-mcols} 
\usepackage{color,upgreek}
\usepackage[ruled,lined]{algorithm2e}
\usepackage{setspace}
\usepackage{pifont}
\usepackage{array}
\usepackage{multirow}
\usepackage{multicol}
\usepackage{arydshln}
\usepackage{titlesec}
\usepackage{appendix}
\usepackage{fancybox,framed} 
\usepackage{manfnt} 
\usepackage{scrextend}
\usepackage{float}
\usepackage{enumitem}
\usepackage{mathtools} 

\usepackage{hyperref}

\usepackage{empheq}

\usepackage[dvipsnames]{xcolor}
\newcommand{\sty}[1]{\boldsymbol{#1}}
\newcommand{\styy}[1]{\mathbb{#1}}
\newcommand{\ubar}[1]{\mkern 0.5mu\underline{\mkern-0.5mu#1\mkern-0.5mu}\mkern 0.5mu}
\newcommand{\uubar}[1]{\ubar{\ubar{#1}}}
\let\mepsilon\epsilon
\let\epsilon\varepsilon

\let\theta\vartheta
\let\mrho\rho
\let\rho\varrho

\let\phi\varphi
\let\Gamma\varGamma
\let\Delta\varDelta
\let\Theta\varTheta
\let\Lambda\varLambda
\let\Xi\varXi
\let\Pi\varPi
\let\Sigma\varSigma
\let\Upsilon\varUpsilon
\let\Phi\varPhi
\let\Psi\varPsi

\let\Omega\varOmega

\newcommand{\grad}[3]{\nabla_{#1}^{#2} #3} 
\newcommand{\var}[3]{\delta_{#1}^{#2} #3} 

\newcommand{\bnd}[3]{\partial_{#1}^{#2} #3} 

\renewcommand{\det}[1]{\mathrm{det}\left(#1\right)} 

\newcommand{\bmm}{\sty{m}}
\newcommand{\bmn}{\sty{n}}

\newcommand{\bmt}{\sty{t}}

\newcommand{\bmB}{\sty{B}}

\newcommand{\bmD}{\sty{D}}

\newcommand{\bmM}{\sty{M}}

\newcommand{\bmP}{\sty{P}}

\newcommand{\bmalpha}{\sty{\alpha}}

\newcommand{\bmmeps}{\sty{\mepsilon}}

\newcommand{\bbN}{\styy{N}}

\newcommand{\bbR}{\styy{R}}

\newcommand{\bbZ}{\styy{Z}}

\newcommand{\clB}{\mathcal{B}}

\newcommand{\clD}{\mathcal{D}}
\newcommand{\clE}{\mathcal{E}}
\newcommand{\clF}{\mathcal{F}}
\newcommand{\clG}{\mathcal{G}}

\newcommand{\clI}{\mathcal{I}}

\newcommand{\clL}{\mathcal{L}}

\newcommand{\clP}{\mathcal{P}}
\newcommand{\clQ}{\mathcal{Q}}
\newcommand{\clR}{\mathcal{R}}
\newcommand{\clS}{\mathcal{S}}
\newcommand{\clT}{\mathcal{T}}

\newcommand{\clV}{\mathcal{V}}

\newcommand{\clX}{\mathcal{X}}

\newcommand{\sT}{\mathsf{T}} 

\newcommand{\scH}{\mathscr{H}}

\newcommand{\scS}{\mathscr{S}}

\newcommand{\rma}{\mathrm{a}}

\newcommand{\rmc}{\mathrm{c}}
\newcommand{\rmd}{\mathrm{d}}
\newcommand{\rme}{\mathrm{e}}

\newcommand{\rmh}{\mathrm{h}}
\newcommand{\rmi}{\mathrm{i}}

\newcommand{\rmp}{\mathrm{p}}

\newcommand{\rmF}{\mathrm{F}}

\newcommand{\uv}{\ubar{v}}
\newcommand{\uw}{\ubar{w}}

\newcommand{\uNull}{\ubar{0}}
\newcommand{\uuA}{\uubar{A}}

\newcommand{\uuC}{\uubar{C}}

\newcommand{\uuI}{\uubar{I}}
\newcommand{\uuJ}{\uubar{J}}

\newcommand{\uuL}{\uubar{L}}
\newcommand{\uuM}{\uubar{M}}

\newcommand{\uuT}{\uubar{T}}
\newcommand{\uuU}{\uubar{U}}

\usepackage{accents}
\robustify{\accentset} 

\newglossary{symbols}{sym}{sbl}{List of Symbols}
\newglossary[glignoredl]{ignored}{glignored}{glignoredin}{Ignored Glossary} 

\newcommand{\mynewglossary}[4]{%
 \newglossaryentry{#2}{type=#1, name={#3}, description={#4}, sort={#2}}%
 \expandafter\newcommand\expandafter{\csname #2\endcsname}{\gls{#2}\xspace}%
}

\newcommand{\mynewglossaryB}[7]{%
 \newglossaryentry{#2}{type=#1, name={#5, #6}, user1={#5}, user2={#6}, description={#7}, sort={#2}}
 \expandafter\newcommand\expandafter{\csname #3\endcsname}{\glsuseri{#2}\xspace}%
 \expandafter\newcommand\expandafter{\csname #4\endcsname}{\glsuserii{#2}\xspace}%
}

\newcommand{\mynewglossaryC}[9]{%
 \newglossaryentry{#2}{type=#1, name={#6, #7, #8}, user1={#6}, user2={#7}, user3={#8}, description={#9}, sort={#2}}
 \expandafter\newcommand\expandafter{\csname #3\endcsname}{\glsuseri{#2}\xspace}%
 \expandafter\newcommand\expandafter{\csname #4\endcsname}{\glsuserii{#2}\xspace}%
 \expandafter\newcommand\expandafter{\csname #5\endcsname}{\glsuseriii{#2}\xspace}%
}

\newcommand{\indDomC}  {\mathrm{c}}                                     
\newcommand{\indDomA}  {\mathrm{a}}                                     
\newcommand{\indDomP}  {\mathrm{p}}                                     
\newcommand{\indHarm}  {\mathrm{h}}                                     
\newcommand{\indAHarm} {\mathrm{ah}}                                    
\newcommand{\indDomI}  {\mathrm{i}}                                     
\newcommand{\indDomIm} {\textnormal{\ensuremath{\mathrm{i\texttt{-}}}}} 
\newcommand{\indDomIpl}{\textnormal{\ensuremath{\mathrm{i\texttt{+}}}}} 
\newcommand{\indCGF}   {\mathrm{cgf}}                                   
\newcommand{\indLGF}   {\mathrm{lgf}}                                   

\mynewglossary{symbols}{real}     {\ensuremath{\bbR}}   {Set of real numbers}
\mynewglossary{symbols}{nat}      {\ensuremath{\bbN}}   {Set of natural numbers}
\mynewglossary{symbols}{integ}    {\ensuremath{\bbZ}}   {Set of integers}
\mynewglossary{ignored}{euclTwo}  {\ensuremath{\real^2}}{Two-dimensional Euclidean space}
\mynewglossary{ignored}{euclThree}{\ensuremath{\real^3}}{Three-dimensional Euclidean space}
\mynewglossary{ignored}{euclD}    {\ensuremath{\real^d}}{$d$-dimensional Euclidean space}

\mynewglossary{symbols}{bnrmVec}{\ensuremath{\bmn}}   {Unit normal vector}
\mynewglossary{symbols}{bperm}  {\ensuremath{\bmmeps}}{Third-order permutation tensor (Levi-Civita tensor)}
\mynewglossary{symbols}{dirac}  {\ensuremath{\delta}} {Dirac delta function}

\mynewglossary{symbols}{lat}       {\ensuremath{\Lambda}}              {Bravais lattice}
\mynewglossary{ignored}{latA}      {\ensuremath{\lat^\indDomA}}        {Atomic lattice}
\mynewglossary{ignored}{latC}      {\ensuremath{\lat^\indDomC}}        {Virtual lattice associated with a continuum domain}
\mynewglossary{ignored}{latP}      {\ensuremath{\lat^\indDomP}}        {Virtual lattice associated with a pad domain}
\mynewglossary{symbols}{domRef}    {\ensuremath{\Omega_0}}             {Reference placement of a continuous body}
\mynewglossary{ignored}{domRefA}   {\ensuremath{\domRef^\indDomA}}     {Reference continuum domain associated with an atomistic region}
\mynewglossary{ignored}{domRefC}   {\ensuremath{\domRef^\indDomC}}     {Reference continuum domain}
\mynewglossary{ignored}{domRefP}   {\ensuremath{\domRef^\indDomP}}     {Reference pad domain}
\mynewglossary{symbols}{dom}       {\ensuremath{\Omega}}               {Current placement of a Continuous body}
\mynewglossary{ignored}{domA}      {\ensuremath{\dom^\indDomA}}        {Continuum domain associated with an atomistic region}
\mynewglossary{ignored}{domC}      {\ensuremath{\dom^\indDomC}}        {Continuum domain}
\mynewglossary{ignored}{domP}      {\ensuremath{\dom^\indDomP}}        {Pad domain}
\mynewglossary{ignored}{domRefCore}{\ensuremath{\domRef^\mathrm{core}}}{Dislocation core domain (reference configuration)}
\mynewglossary{symbols}{domCore}   {\ensuremath{\dom^\mathrm{core}}}   {Dislocation core region}
\mynewglossary{symbols}{intf}      {\ensuremath{\Gamma_\rmi}}          {Artificial sharp interface between the atomistic and the continuum domain}

\mynewglossary{ignored}{mpRef}  {\ensuremath{X}}               {Scalar material point in the reference configuration}
\mynewglossary{ignored}{mpRefB} {\ensuremath{Y}}               {Scalar material point in the reference configuration (2)}
\mynewglossary{symbols}{bmpRef} {\ensuremath{\sty{\mpRef}}}    {Material point in the reference configuration}
\mynewglossary{ignored}{bmpRefB}{\ensuremath{\sty{\mpRefB}}}   {Material point in the reference configuration (2)}
\mynewglossary{ignored}{mpCur}  {\ensuremath{x}}               {Scalar material point in the current configuration}
\mynewglossary{ignored}{mpCurB} {\ensuremath{y}}               {Scalar material point in the current configuration (2)}
\mynewglossary{symbols}{bmpCur} {\ensuremath{\sty{\mpCur}}}    {Material point in the current configuration}
\mynewglossary{ignored}{bmpCurB}{\ensuremath{\sty{\mpCurB}}}   {Material point in the current configuration (2)}
\mynewglossary{symbols}{motion} {\ensuremath{\chi}}            {Motion of a material point \bmpRef}
\mynewglossary{ignored}{displ}  {\ensuremath{u}}               {Scalar displacement field}
\mynewglossary{ignored}{bdispl} {\ensuremath{\sty{\displ}}}    {Displacement field}
\mynewglossary{ignored}{bdisplA}{\ensuremath{\bdispl^\indDomA}}{Atomistic displacement associated with \latA}
\mynewglossary{ignored}{bdisplC}{\ensuremath{\bdispl^\indDomC}}{Atomistic displacement associated with \domC}
\mynewglossary{ignored}{bdisplP}{\ensuremath{\bdispl^\indDomP}}{Atomistic displacement associated with \latP}
\mynewglossary{ignored}{vel}    {\ensuremath{v}}               {Scalar velocity}
\mynewglossary{ignored}{bvel}   {\ensuremath{\sty{\vel}}}      {Velocity vector}

\mynewglossaryC{symbols}{atos} {ato} {atoB} {atoC} {\ensuremath{\xi}}      {\ensuremath{\eta}}      {\ensuremath{\zeta}}      {Atomic indices}
\mynewglossaryC{symbols}{batos}{bato}{batoB}{batoC}{\ensuremath{\sty{\xi}}}{\ensuremath{\sty{\eta}}}{\ensuremath{\sty{\zeta}}}{Positions of atoms \ato, \atoB, \atoC}

\mynewglossary{ignored}{defGrad}   {\ensuremath{F}}                   {Scalar deformation gradient}
\mynewglossary{symbols}{bdefGrad}  {\ensuremath{\sty{\defGrad}}}      {Deformation gradient}
\mynewglossary{symbols}{bdefGradE} {\ensuremath{\bdefGrad^\rme}}      {Elastic contribution of \bdefGrad}
\mynewglossary{symbols}{bdefGradP} {\ensuremath{\bdefGrad^\rmp}}      {Plastic contribution of \bdefGrad}
\mynewglossary{ignored}{strnCG}    {\ensuremath{C}}                   {Scalar right Cauchy-Green strain tensor}
\mynewglossary{symbols}{bstrnCG}   {\ensuremath{\sty{\strnCG}}}       {Right Cauchy-Green strain tensor}
\mynewglossary{ignored}{strnLag}   {\ensuremath{E}}                   {Scalar right Lagrange strain tensor}
\mynewglossary{ignored}{bstrnLag}  {\ensuremath{\sty{\strnLag}}}      {Lagrange strain tensor}
\mynewglossary{ignored}{strnSml}   {\ensuremath{\epsilon}}            {Scalar small strain tensor}
\mynewglossary{symbols}{bstrnSml}  {\ensuremath{\sty{\strnSml}}}      {Small strain tensor}
\mynewglossary{symbols}{bstrnSmlEl}{\ensuremath{\bstrnSml^\rme}}      {Elatic part of the small strain tensor}
\mynewglossary{symbols}{bstrnSmlPl}{\ensuremath{\bstrnSml^\rmp}}      {Plastic part of the small strain tensor}
\mynewglossary{symbols}{bstrnSmlNs}{\ensuremath{\bstrnSml^\mathrm{ns, p}}}{Nonsingular plastic strain}
\mynewglossary{ignored}{bdist}     {\ensuremath{\sty{\beta}}}         {Distortion tensor}
\mynewglossary{symbols}{bdistEl}   {\ensuremath{\sty{\beta}^\rme}}    {Elastic distortion}
\mynewglossary{symbols}{bdistPl}   {\ensuremath{\sty{\beta}^\rmp}}    {Plastic distortion}
\mynewglossary{symbols}{bdistNs}   {\ensuremath{\sty{\beta}^\mathrm{p, ns}}}{Nonsingular plastic distortion}

\mynewglossary{symbols}{Efree}   {\ensuremath{\psi}}               {Helmholtz free energy density}
\mynewglossary{symbols}{Etot}    {\ensuremath{\Pi}}                {Total energy of a mechanical system}
\mynewglossary{ignored}{EtotA}   {\ensuremath{\Etot^\indDomA}}     {Total atomistic energy}
\mynewglossary{ignored}{EtotC}   {\ensuremath{\Etot^\indDomC}}     {Total continuum energy}
\mynewglossary{symbols}{EtotPot} {\ensuremath{\Etot^\mathrm{pot}}} {Potential energy of a mechanical system}
\mynewglossary{symbols}{EtotKin} {\ensuremath{\Etot^\mathrm{kin}}} {Kinetic energy of a mechanical system}
\mynewglossary{symbols}{EtotInt} {\ensuremath{\Etot^\mathrm{int}}} {Internal energy of a mechanical system}
\mynewglossary{symbols}{EtotExt} {\ensuremath{\Etot^\mathrm{ext}}} {External energy of a mechanical system}
\mynewglossary{symbols}{EtotEl}  {\ensuremath{\Etot^\rme}}         {Elastic part of \Etot}
\mynewglossary{symbols}{EtotCore}{\ensuremath{\Etot^\mathrm{core}}}{Dislocation core energy}
\mynewglossary{symbols}{Eato}    {\ensuremath{\clE_\xi}}           {Site energy of an atom \ato}
\mynewglossary{symbols}{Edd}     {\ensuremath{W}}                  {Energy per unit length of the dislocation line}
\mynewglossary{symbols}{Eel}     {\ensuremath{\Edd^\mathrm{e}}}    {Elastic energy per unit length of the dislocation line}
\mynewglossary{symbols}{Ecore}   {\ensuremath{\Edd^\mathrm{core}}} {Core energy per unit length of the dislocation line}

\mynewglossary{ignored}{stressSml}  {\ensuremath{\sigma}}          {Scalar Cauchy stress}
\mynewglossary{symbols}{bstressSml} {\ensuremath{\sty{\stressSml}}}{Cauchy stress tensor}
\mynewglossary{symbols}{stressShear}{\ensuremath{\tau}}            {Scalar shear stress}
\mynewglossary{symbols}{bstressPK}  {\ensuremath{\bmP}}            {First Piola-Kirchhoff stress tensor}
\mynewglossary{symbols}{btract}     {\ensuremath{\bmt}}            {Traction vector}
\mynewglossary{symbols}{bstressEsh} {\ensuremath{\bmB}}            {Eshelby stress tensor}

\mynewglossary{ignored}{force} {\ensuremath{f}}                    {Scalar force}
\mynewglossary{symbols}{bforce}{\ensuremath{\sty{\force}}}         {Force vector}
\mynewglossary{ignored}{bfA}   {\ensuremath{\bforce^\indDomA}}     {Atomic force vector}
\mynewglossary{ignored}{bfC}   {\ensuremath{\bforce^\indDomC}}     {Continuum force vector}
\mynewglossary{symbols}{bfbody}{\ensuremath{\bforce^\mathrm{body}}}{Body force}
\mynewglossary{symbols}{bfPK}  {\ensuremath{\bforce^\mathrm{pk}}}  {Peach-Koehler force}
\mynewglossary{symbols}{bfCore}{\ensuremath{\bforce^\mathrm{core}}}{Force on \DDcur due to core energy contribution \EtotCore}

\mynewglossary{ignored}{stiffMat} {\ensuremath{C}}               {Scalar material stiffness}
\mynewglossary{symbols}{bstiffMat}{\ensuremath{\styy{\stiffMat}}}{Material stiffness tensor}
\mynewglossary{symbols}{shearMod} {\ensuremath{\mu}}             {Shear modulus}
\mynewglossary{symbols}{poisRatio}{\ensuremath{\nu}}             {Poisson ratio}
\mynewglossary{ignored}{green}    {\ensuremath{G}}               {Scalar entry of the Green tensor}
\mynewglossary{symbols}{bgreen}   {\ensuremath{\sty{\green}}}    {Second-order Green tensor}

\mynewglossary{symbols}{slipPlane}     {\ensuremath{\scS}}                        {Slip plane}
\mynewglossary{ignored}{DDref}         {\ensuremath{\Gamma}}                      {Discrete dislocation line with respect to a reference configuration}
\mynewglossary{ignored}{DDrefA}        {\ensuremath{\DDref^\indDomA}}             {Discrete dislocation line in the atomistic domain with respect to a reference configuration}
\mynewglossary{ignored}{DDrefC}        {\ensuremath{\DDref^\indDomC}}             {Discrete dislocation line in the continuum domain with respect to a reference configuration}
\mynewglossary{symbols}{DDrefHyb}      {\ensuremath{\DDref^\mathrm{hyb}}}         {Hybrid discrete dislocation line with respect to a reference configuration}
\mynewglossary{symbols}{DDcur}         {\ensuremath{\gamma}}                      {Discrete dislocation line}
\mynewglossary{ignored}{DDa}           {\ensuremath{\DDcur^\indDomA}}             {Discrete dislocation line in the atomistic domain}
\mynewglossary{ignored}{DDc}           {\ensuremath{\DDcur^\indDomC}}             {Discrete dislocation line in the continuum domain}
\mynewglossary{symbols}{DDhyb}         {\ensuremath{\DDcur^\mathrm{hyb}}}         {Hybrid discrete dislocation line}
\mynewglossary{symbols}{DDvir}         {\ensuremath{\DDcur^\mathrm{v}}}           {Virtual subset of a discrete dislocation line outside a bounded domain}
\mynewglossary{ignored}{mpDDref}       {\ensuremath{S}}                           {Scalar material point on \DDref}
\mynewglossary{ignored}{mpDDcur}       {\ensuremath{s}}                           {Scalar material point on \DDcur}
\mynewglossary{symbols}{bmpDDref}      {\ensuremath{\sty{\mpDDref}}}              {Material point on \DDref}
\mynewglossary{symbols}{bmpDDcur}      {\ensuremath{\sty{\mpDDcur}}}              {Material point on \DDcur}
\mynewglossary{symbols}{btanVecDD}     {\ensuremath{\bmt}}                        {Unit tangent vector of the dislocation line}
\mynewglossary{symbols}{bnrmVecDD}     {\ensuremath{\bmm}}                        {Unit normal vector of the dislocation line}
\mynewglossary{symbols}{burgers}       {\ensuremath{b}}                           {Magnitude of the Burgers vector}
\mynewglossary{symbols}{bburgers}      {\ensuremath{\sty{\burgers}}}              {Burgers vector}
\mynewglossary{symbols}{bdislDens}     {\ensuremath{\bmalpha}}                    {Dislocation density tensor}
\mynewglossary{symbols}{cutOffCore}    {\ensuremath{r_\mathrm{cut}^\mathrm{core}}}{Dislocation core cut-off radius}
\mynewglossary{symbols}{charAngle}     {\ensuremath{\theta}}                      {Character angle of the dislocation line}
\mynewglossary{ignored}{EcoreParam}    {\ensuremath{E^\mathrm{core}}}             {Scalar core energy parameter}
\mynewglossary{ignored}{spreadIso}     {\ensuremath{w}}                           {Isotropic spreading function of the Burgers vector}
\mynewglossary{ignored}{spreadIsoW}    {\ensuremath{a}}                           {Isotropic core spreading width}
\mynewglossary{symbols}{bmobility}     {\ensuremath{\bmM}}                        {Dislocation mobility tensor}
\mynewglossary{symbols}{bdrag}         {\ensuremath{\bmD}}                        {Dislocation drag tensor}
\mynewglossary{ignored}{displDD}       {\ensuremath{\tilde{\displ}}}              {Scalar displacement of an isolated dislocation in an unbounded domain}
\mynewglossary{ignored}{displCorr}     {\ensuremath{\hat{\displ}}}                {Scalar corrective displacement}
\mynewglossary{ignored}{bdisplDD}      {\ensuremath{\tilde{\bdispl}}}             {Displacement field of due to a dislocation}
\mynewglossary{ignored}{bdisplCorr}    {\ensuremath{\hat{\bdispl}}}               {Corrective displacement field}
\mynewglossary{ignored}{bstressSmlDD}  {\ensuremath{\tilde{\bstressSml}}}         {Stress field of an isolated dislocation in an unbounded domain}
\mynewglossary{ignored}{bstressSmlCorr}{\ensuremath{\hat{\bstressSml}}}           {Corrective stress field}

\mynewglossary{ignored}{cutOff}  {\ensuremath{r_\mathrm{cut}}}{Cut-off radius}
\mynewglossary{ignored}{intRg}   {\ensuremath{\clR}}          {Interaction range}
\mynewglossary{ignored}{intRgAto}{\ensuremath{\clR_{\ato}}}   {Interaction range of atom \bato}
\mynewglossary{ignored}{latConst}{\ensuremath{a_0}}           {Lattice constant}

\mynewglossary{symbols}{intvar}{\ensuremath{\hat{\alpha}}}{Vector of internal state variables}

\mynewglossary{symbols}{Prob}      {\ensuremath{\clP}}                       {Shorthand notation of a problem definition}
\mynewglossary{ignored}{Pa}        {\ensuremath{\Prob^\indDomA}}             {Atomistic problem}
\mynewglossary{ignored}{Pc}        {\ensuremath{\Prob^\indDomC}}             {Continuum problem}
\mynewglossary{symbols}{Pp}        {\ensuremath{\Prob^\rmp}}                 {Physical problem}
\mynewglossary{symbols}{PcP}       {\ensuremath{\Prob^{\indDomC/\rmp}}}      {Physical subproblem}
\mynewglossary{symbols}{PcDD}      {\ensuremath{\Prob^{\indDomC/\rmd\rmd}}}  {Material subproblem}
\mynewglossary{symbols}{Pcadd}     {\ensuremath{\Prob^\mathrm{cadd}}}        {CADD problem}
\mynewglossary{symbols}{PcaddApprx}{\ensuremath{\tilde{\Prob}^\mathrm{cadd}}}{Approximate CADD problem}

\mynewglossary{symbols}{bdisplDDt}   {\ensuremath{\tilde{\bdispl}^\mathrm{t}}}         {Core template for an infinite straight dislocation}
\mynewglossary{symbols}{bdisplDDcorr}{\ensuremath{\Delta\tilde{\bdispl}^\mathrm{corr}}}{Core template correction}
\mynewglossary{symbols}{transNode}   {\ensuremath{\bmpDDcur^\mathrm{trans}}}           {Transmission node}
\mynewglossary{symbols}{domDetn}     {\ensuremath{\dom^\mathrm{detn}}}                 {Detection domain}

\mynewglossary{ignored}{latI}  {\ensuremath{\lat^\indDomI}}  {Interface region}
\mynewglossary{ignored}{latIm} {\ensuremath{\lat^\indDomIm}} {Interface (-) region}
\mynewglossary{ignored}{latIpl}{\ensuremath{\lat^\indDomIpl}}{Interface (+) region}

\mynewglossary{ignored}{bdisplI}  {\ensuremath{\bdispl^\indDomI}}  {Atomistic displacement associated with \latI}
\mynewglossary{ignored}{bdisplIm} {\ensuremath{\bdispl^\indDomIm}} {Atomistic displacement associated with \latIm}
\mynewglossary{ignored}{bdisplIpl}{\ensuremath{\bdispl^\indDomIpl}}{Atomistic displacement associated with \latIpl}

\mynewglossary{ignored}{diffOp} {\ensuremath{\clL}}                         {Differential operator}
\mynewglossary{ignored}{Lh}    {\ensuremath{\diffOp_\indHarm}}    {Linearized operator}
\mynewglossary{ignored}{Lah}   {\ensuremath{\diffOp_\indAHarm}}   {Anharmonic operator}
\mynewglossary{ignored}{La}   {\ensuremath{\diffOp^\indDomA}}               {Nonlinear atomistic operator}
\mynewglossary{ignored}{LaH}    {\ensuremath{\diffOp_\indHarm^\indDomA}}    {Linearized atomistic operator}
\mynewglossary{ignored}{LaAH}   {\ensuremath{\diffOp_\indAHarm^\indDomA}}   {Anharmonic atomistic operator}
\mynewglossary{ignored}{Lc}     {\ensuremath{\diffOp_\indHarm^\indDomC}}             {Continuum operator}
\mynewglossary{ignored}{Laa}    {\ensuremath{\diffOp_\indHarm^{\indDomA|\indDomA}}}  {Coupling operator}
\mynewglossary{ignored}{Lcc}    {\ensuremath{\diffOp_\indHarm^{\indDomC|\indDomC}}}  {Coupling operator}
\mynewglossary{ignored}{Lac}    {\ensuremath{\diffOp_\indHarm^{\indDomA|\indDomC}}}  {Coupling operator}
\mynewglossary{ignored}{Lca}    {\ensuremath{\diffOp_\indHarm^{\indDomC|\indDomA}}}  {Coupling operator}
\mynewglossary{ignored}{Lic}    {\ensuremath{\diffOp_\indHarm^{\indDomI|\indDomC}}}  {Coupling operator}
\mynewglossary{ignored}{Lii}    {\ensuremath{\diffOp_\indHarm^{\indDomI|\indDomI}}}  {Coupling operator}
\mynewglossary{ignored}{Lci}    {\ensuremath{\diffOp_\indHarm^{\indDomC|\indDomI}}}  {Coupling operator}
\mynewglossary{ignored}{Liipl}  {\ensuremath{\diffOp_\indHarm^{\indDomI|\indDomIpl}}}{Coupling operator}
\mynewglossary{ignored}{Lipla}  {\ensuremath{\diffOp_\indHarm^{\indDomIpl|\indDomA}}}{Coupling operator}
\mynewglossary{ignored}{Lipli}  {\ensuremath{\diffOp_\indHarm^{\indDomIpl|\indDomI}}}{Coupling operator}
\mynewglossary{ignored}{greenOp}{\ensuremath{\clG}}                         {Green operator}
\mynewglossary{ignored}{Ga}     {\ensuremath{\greenOp^\indDomA}}            {Atomistic Green operator}
\mynewglossary{ignored}{Gc}     {\ensuremath{\greenOp^\indDomC}}            {Continuum Green operator}
\mynewglossary{ignored}{Gaa}    {\ensuremath{\greenOp^{\indDomA|\indDomA}}} {Coupling operator}
\mynewglossary{ignored}{Gcc}    {\ensuremath{\greenOp^{\indDomC|\indDomC}}} {Coupling operator}
\mynewglossary{ignored}{Gac}    {\ensuremath{\greenOp^{\indDomA|\indDomC}}} {Coupling operator}
\mynewglossary{ignored}{Gca}    {\ensuremath{\greenOp^{\indDomC|\indDomA}}} {Coupling operator}
\mynewglossary{ignored}{Gci}    {\ensuremath{\greenOp^{\indDomC|\indDomI}}} {Coupling operator}
\mynewglossary{ignored}{Gcipl}    {\ensuremath{\greenOp^{\indDomC|\indDomIpl}}} {Coupling operator}
\mynewglossary{ignored}{Gipli}    {\ensuremath{\greenOp^{\indDomIpl|\indDomI}}} {Coupling operator}
\mynewglossary{ignored}{Giipl}    {\ensuremath{\greenOp^{\indDomI|\indDomIpl}}} {Coupling operator}
\mynewglossary{ignored}{Giplipl}    {\ensuremath{\greenOp^{\indDomIpl|\indDomIpl}}} {Coupling operator}
\mynewglossary{ignored}{Gii}    {\ensuremath{\greenOp^{\indDomI|\indDomI}}} {Coupling operator}
\mynewglossary{ignored}{Gpi}    {\ensuremath{\greenOp^{\indDomP|\indDomI}}} {Coupling operator}
\mynewglossary{ignored}{Gic}    {\ensuremath{\greenOp^{\indDomI|\indDomC}}} {Coupling operator}

\mynewglossary{ignored}{bK}       {\ensuremath{\sty{K}}}        {Nodal stiffness tensor}
\mynewglossary{ignored}{greenCGF} {\ensuremath{\green^\indCGF}} {Scalar continuum Green function}
\mynewglossary{ignored}{bgreenCGF}{\ensuremath{\bgreen^\indCGF}}{Continuum Green tensor}
\mynewglossary{ignored}{greenLGF} {\ensuremath{\green^\indLGF}} {Scalar lattice Green function}
\mynewglossary{ignored}{bgreenLGF}{\ensuremath{\bgreen^\indLGF}}{Lattice Green tensor}

\mynewglossary{ignored}{forceOp}{\ensuremath{\clF}}{Force operator}
\mynewglossary{ignored}{hdiffOp}{\ensuremath{\hat{L}}}{Matrix representation of \diffOp}
\mynewglossary{ignored}{hgreenOp}{\ensuremath{\hat{G}}}{Matrix representation of \greenOp}
\mynewglossary{ignored}{hforceOp}{\ensuremath{\hat{F}}}{Matrix representation of \forceOp}

\newglossary{newsymbols}{newsym}{newsbl}{List of Symbols}
\glsmoveentry{integ}{newsymbols}
\glsmoveentry{real}{newsymbols}
\glsmoveentry{euclD}{newsymbols}
\glsmoveentry{lat}{newsymbols}
\glsmoveentry{latConst}{newsymbols}
\glsmoveentry{atos}{newsymbols}
\glsmoveentry{batos}{newsymbols}
\glsmoveentry{bdispl}{newsymbols}
\glsmoveentry{bdisplDD}{newsymbols}
\glsmoveentry{bforce}{newsymbols}
\glsmoveentry{intRgAto}{newsymbols}
\glsmoveentry{Eato}{newsymbols}
\glsmoveentry{Etot}{newsymbols}
\glsmoveentry{bstiffMat}{newsymbols}
\glsmoveentry{bK}{newsymbols}
\glsmoveentry{bgreenCGF}{newsymbols}
\glsmoveentry{bgreenLGF}{newsymbols}
\glsmoveentry{burgers}{newsymbols}
\glsmoveentry{cutOff}{newsymbols}
\glsmoveentry{diffOp}{newsymbols}
\glsmoveentry{Lh}{newsymbols}
\glsmoveentry{Lah}{newsymbols}
\glsmoveentry{greenOp}{newsymbols}
\glsmoveentry{forceOp}{newsymbols}
\glsmoveentry{hdiffOp}{newsymbols}
\glsmoveentry{hgreenOp}{newsymbols}
\glsmoveentry{hforceOp}{newsymbols}
\glsmoveentry{Prob}{newsymbols}

\makeglossaries

\glsdisablehyper 
\renewcommand{\glossarysection}[2][]{} 
\glssetwidest[0]{\hspace{3.5em}}

\newtheorem{thm}{Theorem}
\newtheorem{prop}{Proposition}
\newtheorem{rem}{Remark}
\newtheorem{defn}{Definition}

\newtheorem{cor}{Corollary}

\newproof{prf}{Proof}
\newtheorem{innercustomlem}{Lemma}
\newenvironment{customlem}[1]
 {\renewcommand\theinnercustomlem{#1}\innercustomlem}
 {\endinnercustomlem}

\titleformat*{\section}{\bfseries}

\newcommand*{\textcal}[1]{%
 \textit{\large\fontfamily{qzc}\selectfont#1}%
}

\newlength{\commentWidth}

\makeatletter
\def\blfootnote{\gdef\@thefnmark{}\@footnotetext}
\makeatother

\newcommand{\st}{\scriptstyle}
\newcommand{\sst}{\scriptscriptstyle}

\newcommand{\mrm}[1]{\mathrm{#1}}
\newcommand{\inv}[1]{\ensuremath{{#1}^{-1}}}

\newcommand{\latb}{\ensuremath{\bar{\lat}}}

\renewcommand{\l}{\ensuremath{{\scriptscriptstyle\lat}}}
\newcommand{\lb}{\ensuremath{{\bar{\l}}}}
\renewcommand{\a}{\indDomA}
\renewcommand{\c}{\indDomC}
\newcommand{\cb}{\ensuremath{{\bar{\c}}}}
\newcommand{\re}{\ensuremath{{\mrm{re}}}}
\newcommand{\p}{\indDomP}
\newcommand{\pr}{\ensuremath{{\p'}}}
\newcommand{\ip}{\indDomIpl}
\renewcommand{\i}{\indDomI}

\newcommand{\I}{\textnormal{\ensuremath{\mathrm{o}}}}
\renewcommand{\Im}{\textnormal{\ensuremath{\mathrm{o\texttt{-}}}}}
\renewcommand{\aa}{\ensuremath{{\a|\a}}}
\newcommand{\cc}{\ensuremath{{\c|\c}}}
\newcommand{\ac}{\ensuremath{{\a|\c}}}
\newcommand{\ca}{\ensuremath{{\c|\a}}}
\newcommand{\pp}{\ensuremath{{\p|\p}}}
\newcommand{\II}{\ensuremath{{\I|\I}}}

\newcommand{\ci}{\ensuremath{{\c|\i}}}

\newcommand{\ap}{\ensuremath{{\a|\p}}}

\newcommand{\prp}{\ensuremath{{\pr|\p}}}

\newcommand{\ipr}{\ensuremath{{\i|\pr}}}
\newcommand{\ipi}{\ensuremath{{\ip|\i}}}
\newcommand{\cip}{\ensuremath{{\c|\ip}}}
\newcommand{\Ic}{\ensuremath{{\I|\c}}}
\newcommand{\cI}{\ensuremath{{\c|\I}}}
\newcommand{\cp}{\ensuremath{{\c|\p}}}
\newcommand{\Iip}{\ensuremath{{\I|\ip}}}
\newcommand{\Ii}{\ensuremath{{\I|\i}}}
\newcommand{\ipd}{\ensuremath{{\i|\p}}}
\renewcommand{\L}{\diffOp\xspace}
\newcommand{\Lhnl}{\ensuremath{\L_\mrm{hnl}}\xspace}
\newcommand{\Lcpl}{\ensuremath{\diffOp_\mathrm{cpl}}\xspace} 
\newcommand{\G}{\greenOp\xspace}

\newcommand{\F}{\forceOp\xspace}
\newcommand{\B}{\ensuremath{\clB}\xspace}
\newcommand{\projop}{\ensuremath{\clT}\xspace}
\renewcommand{\P}{\projop}
\newcommand{\Pinf}{\ensuremath{\tilde{\projop}}\xspace}
\newcommand{\Pbnd}{\ensuremath{\hat{\projop}}\xspace}
\newcommand{\Id}{\ensuremath{\clI}\xspace}
\newcommand{\nullOp}{\ensuremath{\textcal{0}\,}\xspace}

\newcommand{\displInf}{\ensuremath{\displ_\infty}}
\newcommand{\displHom}{\ensuremath{\displ_\mrm{hom}}}
\newcommand{\latInf}{\ensuremath{\lat_\infty}}

\def\higlight{0}

\newenvironment{modified}{%
\if\higlight1%
\color{blue}%
\fi%
}{}

\newenvironment{moved}{%
\if\higlight1%
\color{red}%
\fi%
}{}

\newenvironment{movedAndModified}{%
\if\higlight1%
\color{magenta}%
\fi%
}{}

\newenvironment{new}{%
\if\higlight1%
\color{OliveGreen}%
\fi%
}{}

\RequirePackage[normalem]{ulem} 
\RequirePackage{color}\definecolor{RED}{rgb}{1,0,0}\definecolor{BLUE}{rgb}{0,0,1} 

\begin{document}

\begin{frontmatter}
 \blfootnote{\textit{E-mail address:} \href{mailto:m.hodapp@skoltech.ru}{m.hodapp@skoltech.ru}}
 \begin{abstract}
  The ``flexible boundary condition'' method, introduced by Sinclair and coworkers in the 1970s, remains among the most popular methods for simulating isolated two-dimensional crystalline defects, embedded in an effectively infinite atomistic domain. In essence, the method can be characterized as a domain decomposition method which iterates between a local anharmonic and a global harmonic problem, where the latter is solved by means of the lattice Green function of the ideal crystal. This local/global splitting gives rise to tremendously improved convergence rates over related alternating Schwarz methods. In a previous publication (Hodapp et al., 2019, Comput. Methods in Appl. Mech. Eng. 348), we have shown that this method also applies to large-scale three-dimensional problems, possibly involving hundreds of thousands of atoms, using fast summation techniques exploiting the low-rank nature of the asymptotic lattice Green function. 
  
  Here, we generalize the Sinclair method to bounded domains and develop an implementation using a discrete boundary element method to correct the infinite solution with respect to a prescribed far-field condition, thus preserving the advantage of the original method of not requiring a global spatial discretization. Moreover, we present a detailed convergence analysis and show for a one-dimensional problem that the method is unconditionally stable under physically motivated assumptions. To further improve the convergence behavior, we develop an acceleration technique based on a relaxation of the transmission conditions between the two subproblems. Numerical examples for linear and nonlinear problems are presented to validate the proposed methodology.
 \end{abstract}
\end{frontmatter}



\section{Introduction}
\label{sec:intro}

Computational physics has become a valuable tool for studying the material behavior on the nanoscale due to the vast advances in computing technology over the past 20--30 years. Presently, quantum-mechanical systems of a few hundred atoms can be simulated with density functional theory in order to predict basic material properties, such as elastic constants and phase stability. For larger systems, density functional theory becomes too expensive and problems on the nanoscale are, in lieu thereof, carried out with atomistic models.

However, atomistic models are also limited in size and time, motivating the so-called hierarchical multiscale approach in which some key parameters are calculated atomistically, such as energy barriers for defect motion, and subsequently passed to the next higher-scale level, e.g., dislocation dynamics models (see, e.g., \citep{argon_strengthening_2007}). This approach works well if the relevant mechanisms are confined to a single scale. However, there are situations where the material behavior on the smaller scale is strongly influenced by processes taking place on higher scales which requires a \emph{concurrent multiscale approach}.

The need for concurrent multiscale models has motivated the development of atomistic/continuum (A/C) coupling methods \citep{kohlhoff_new_1989,tadmor_quasicontinuum_1996,shilkrot_multiscale_2004,shimokawa_matching_2004,kochmann_meshless_2014}. Thereby, only the material behavior in the vicinity of crystalline defects, e.g., interstitials, vacancies, dislocations, cracks or grain boundaries, is treated fully-atomistically. This fully atomistic domain is surrounded by a significantly cheaper continuum elasticity region allowing for much larger computational domains to take scale-bridging effects into account.

A/C coupling methods can broadly be grouped into energy- and force-based methods \citep{curtin_atomistic/continuum_2003}. Energy-based methods define a global energy functional which, ideally, closely reassembles the one of the fully atomistic model. However, the still existing challenge is the construction of a consistent coupling of the energy between the atomistic and continuum domains due to the nonlocal-local mismatch between the two models. To date, a myriad of sophisticated approaches has been developed but fully consistent methods only exist for two-dimensional problems (e.g., \citep{ortner_energy-based_2014}). This motivates the use of unconditionally consistent force-based methods. In turn, however, their stability properties have not yet been fully understood \citep{dobson_stability_2010}. Moreover, force-based methods cannot be formulated as energy minimization problems restricting the margin of possible numerical solvers.

An important task is thus to construct stable solvers for force-based A/C coupling methods. \citet{dobson_iterative_2011} have proposed a monolithic Newton-GMRes method for a linearized toy model which was has been shown to approximate the stability region of the fully atomistic model. \citet{hodapp_lattice_2019} have extended this idea to nonlinear problems by introducing an additional criterion in order to correct the search direction whenever the atomistic Hessian becomes indefinite---which allows to overcome energy barriers, e.g., during dislocation motion. However, these methods are difficult to implement \emph{and} parallelize in practice. A different approach is domain decomposition based on the alternating Schwarz method \citep{toselli_domain_2005}, where the coupled problem is solved sequentially by means of a fixed point iteration \citep{parks_connecting_2008,li_efficient_2009}. The advantage of this method is the iteration between two \emph{symmetric} energy minimization problems which can be solved with standard techniques. The disadvantage is, however, its poor convergence behavior which makes the method impractical as a stand-alone solver since it is potentially even slower than the fully atomistic model.

A fast alternative to the alternating Schwarz method is the \emph{``flexible boundary condition method''} developed by Sinclair and coworkers during the 1970s \citep{sinclair_improved_1971,sinclair_influence_1975,sinclair_flexible_1978}. This method can be described as a fixed point iteration between a \emph{local anharmonic} (a.k.a. the fully atomistic) and a \emph{global harmonic} problem. It can be shown that this splitting gives rise to significantly faster convergence rates which even compete with monolithic Krylov subspace solvers \citep{hodapp_lattice_2019}. The Sinclair method has been proposed for effectively infinite problems, where the solution of the global problem reduces to a single matrix-vector multiplication exploiting the Green function of the harmonic operator, and successfully employed to study the behavior of isolated two-dimensional defects, particularly dislocations (e.g., \citep{rao_greens_1998,woodward_flexible_2002,fellinger_geometries_2018}).

\begin{modified}%
In this work, we generalize the Sinclair method to bounded problems and analyze its convergence properties. Moreover, we present an efficient implementation by correcting the solution of the infinite problem with respect to a prescribed far-field by solving the global harmonic problem using a discrete boundary element method \citep{martinsson_fast_2002,li_atomistic-based_2012,hodapp_lattice_2019}. This approach is advantageous over volume-based methods, e.g., finite element methods, since boundary element methods do not require an explicit discretization of the interior of the domain which is not essential since practitioners are generally interested in the material behavior solely inside the fully atomistic domain. Moreover, the original Sinclair method is recovered if the outer boundary vanishes.

The work is organized as follows. In Section \ref{sec:atom} and \ref{sec:con}, we introduce the atomistic and continuum models which are subsequently used in Section \ref{sec:acc} to formulate the coupled problem. In Section \ref{sec:sinc}, the Sinclair method for bounded problems is described and a detailed description of the proposed solution algorithm is covered. Section \ref{sec:conv} is then devoted to the convergence analysis of the proposed algorithm. We further show for a one-dimensional problem that the algorithm is uniformly stable under physically motivated assumptions. We also introduce a technique to accelerate the fixed point iteration using a relaxation of the transmission conditions. Numerical examples are presented in Section \ref{sec:examples} to validate the convergence analysis and to demonstrate the efficiency of the proposed method.
\end{modified}%

\section{Notation}
\label{sec:notation}

We begin by setting up some generic notation for defining computational problems on discrete domains. Other notation will be introduced as required throughout the manuscript.

\begin{new}%
\subsection{Computational problems defined on discrete domains}
\label{sec:notation.discrete_problem}

In the present manuscript, we will work with functions defined on discrete domains $\lat \subset \real^d$, where $d=1,...,3$ is the dimension. In this section, we exemplify the notation of a computational problem defined on $\lat$ as will be used in the remainder of this work. We further assert that all tensorial quantities are defined with respect to the usual orthonormal basis system $\{ e_i \in \real^d \,\vert\, e_ie_j = \delta_{ij} \}_{i=1,...,d}$.

Let $u \in \clV = \clV(\lat) := \{\, v : \lat \rightarrow \real^d \,\}$.
We will denote operators acting on $\clV$ using calligraphic symbols ($\clL, \clG,\clT$, etc.).
For example, let $\clL : \clV \rightarrow \clV^\ast$, where $\clV^\ast \simeq \clV$ is the dual space of $\clV$. In general we have $f \in \clV^\ast$ given and are looking for $u$ such that
\begin{equation}\label{eq:notation.Lu=f}
 \clL[\displ](\xi) = f(\xi) \qquad \text{in} \; \lat.
\end{equation}
If not explicitly required we will drop the argument $\xi$ and write $\clL[u] = f \; \text{in} \; \lat$ for brevity. Instead of \eqref{eq:notation.Lu=f}, we additionally make use of the following two \emph{equivalent} notations
\begin{align}
 \text{(i)} \quad \clL^{\l|\l}[u] = f^{\l}, && \text{(ii)} \quad \clL^{\l|}[u] = f^{\l}.
\end{align}
The superscripted indices ``$\l|\l$'' in notation (i) specify that the operator $\L$ is acting on functions defined on $\lat$ and produces functions $f^\l$ defined on $\lat$ . Here, $\L$ acts on the entire domain on which it is defined and we may therefore omit the second index (since this information is already contained in the definition of $\L$) and just use ``$\l|$'' according to notation (ii).

This notation will prove useful when working with subspaces of $\clV$. For example, let $\latA,\latC \subset \lat$ such that $\latA \cup \latC = \lat$ and the corresponding decomposition of $\clV$ as $\clV = \clV(\latA) \oplus \clV(\latC)$. Let now $\clL^{\a|} : \clV \rightarrow \clV^\ast(\latA)$ and $\clL^{\c|} : \clV \rightarrow \clV^\ast(\latC)$. Then, we can write \eqref{eq:notation.Lu=f} equivalently as
$
\clL[u] \equiv
\begin{psmallmatrix} \clL^{\a|}[\displ] \\ \clL^{\c|}[u] \end{psmallmatrix}
=
\begin{psmallmatrix} f^\a \\ f^\c \end{psmallmatrix}
$.
More precisely, the superscripted index ``$\a|$'' implies here that the operator is acting on functions defined on the entire domain $\lat$ and produces functions defined on $\latA$; and vice versa for the index ``$\c|$''.

If \eqref{eq:notation.Lu=f} is a linear equation, the operator $\clL$ does not depend on $u$ and can be defined as follows
\begin{equation}
 \begin{aligned}
  \clL : \clV &\rightarrow \clV^\ast \\
         v    &\mapsto     \clL[v] \qquad \text{such that} \; \forall\,\ato\in\lat \qquad
         \clL[v](\ato) = \sum_{\atoB \in \lat} L(\ato,\atoB)v(\atoB),
 \end{aligned}
\end{equation}
where $L(\ato,\atoB)$ denotes the kernel of $\clL$. Linearity allows us to split also the domain space of $\clL$ and define $\clL^{|\a} : \clV(\latA) \rightarrow \clV^\ast$ and $\clL^{|\c} : \clV(\latC) \rightarrow \clV^\ast$. Then, we can equivalently write
$\clL[u] = \clL^{|\a}[u] + \clL^{|\c}[u]$ or, using the expanded notation,
\begin{equation}
 \clL[u](\xi)
 =
 \clL^{|\a}[u](\xi) + \clL^{|\c}[u](\xi)
 =
 \sum_{\eta \in \latA} L(\xi,\eta) u(\eta) + \sum_{\eta \in \latC} L(\xi,\eta) u(\eta),
\end{equation}
where the superscripted index ``$|\a$'' here implies that the operator is acting on functions defined on $\latA$ and produces functions defined on the entire domain $\lat$; and vice versa for the index ``$|\c$''.

This splitting of domain and co-domain spaces naturally leads to the matrix notation
\begin{equation}
 \clL[u] \equiv
 \begin{pmatrix}
  \clL^\aa & \clL^\ac \\ \clL^\ca & \clL^\cc
 \end{pmatrix}
 \begin{bmatrix}
  u^\a \\ u^\c
 \end{bmatrix}
 =
 \begin{pmatrix}
  \clL^\aa[u] + \clL^\ac[u] \\ \clL^\ca[u] + \clL^\cc[u]
 \end{pmatrix}
 =
 \begin{pmatrix}
  \clL^{\a|}[u] \\ \clL^{\c|}[u]
 \end{pmatrix}
 =
 \begin{pmatrix}
  f^\a \\ f^\c
 \end{pmatrix},
\end{equation}
with
$\clL^{\a|\a} : \clV(\latA) \rightarrow \clV^\ast(\latA)$,
$\clL^{\a|\c} : \clV(\latA) \rightarrow \clV^\ast(\latC)$,
$\clL^{\c|\a} : \clV(\latC) \rightarrow \clV^\ast(\latA)$,
and
$\clL^{\c|\c} : \clV(\latC) \rightarrow \clV^\ast(\latC)$.

For an overview of all domains and their associate indices, used for functions and operators throughout this work, the reader is referred to Appendix \ref{apdx:dom_idx}.
\end{new}%

\subsection{Norms}

For all $v \in \clV$ we require the usual $l^2$- and $l^\infty$-norms 
\begin{align}
 \| v \|            &= \left( \sum_{\ato \in \lat} \vert v(\ato) \vert^2 \right)^{1/2}, \\
 \| v \|_{l^\infty} &= \arg{\left\{ \underset{\ato \in \lat}{\max}{\, v(\ato)} \right\}},
\end{align}
as well as the inner product
\begin{equation}
 \langle v, w \rangle = \sum_{\ato \in \lat} v(\ato) w(\ato), \qquad \text{with} \; w \in \clV.
\end{equation}
For norms of functions which are elements of a subspace of $\clV$, e.g., $v^\a \in \clV(\latA)$, we write $\| v^\a \|$ etc.

In some cases we explicitly consider an operator $\clL$ as a mapping $\clL : l^2 \rightarrow l^2$ (cf. Section \ref{sec:conv}). The corresponding operator norm, induced by $l^2$, is defined as
\begin{equation}
 \| \clL \| = \underset{v \ne 0}{\sup}{\, \frac{\| \clL[v] \|}{\| v \|}} = \underset{\| v \| = 1}{\sup}{\, \| \clL[v] \|} = \mrho(\clL),
\end{equation}
where $\mrho(\clL)$ is the largest singular value of $\clL$.


\section{Reference atomistic problem}
\label{sec:atom}

\begin{modified}%
Let $\latInf$ be a Bravais lattice
\begin{equation}\label{eq:atom.bravais}
 \latInf := \left\{\, \sum_{i=1}^d a_iv_i \,\bigg\vert\, a_i \in \integ \,\right\},
\end{equation}
where the set of basis vectors $\{ v_i \in \real^d \}_{i=1,...,d}$ defines the lattice type, e.g., hexagonal, face-centered cubic or body-centered cubic. A deformation of $\latInf$ will be described via displacements $\displ : \latInf \rightarrow \real^d$ and we assume that we can identify $\displ$ with some sufficiently smooth interpolant which allows us to define the gradient $\grad{}{}{\displ}$.

Our \emph{computational atomistic domain} is the subset $\lat \subseteq \latInf$ which may contain crystalline defects as illustrated for a hexagonal lattice in Figure \ref{fig:domains}. In this work, we consider Dirichlet-type problems, that is, the displacements of atoms outside of $\lat$ are fixed to some $\bar{\displ}$. Atoms in $\lat$ interact with the set of boundary atoms in $\partial\lat \subset (\latInf \setminus \lat)$ whose size depends on the range of the interatomic potential (see below). The union of $\lat$ and $\partial\lat$ is $\bar{\lat} = \lat \cup \partial\lat$. Hence, we define the space of \emph{admissible finite energy displacements} as  follows(cf. \citep{luskin_atomistic--continuum_2013})
\begin{equation}
 \clV(\clX) := \{\, v : \clX \rightarrow \real^d \,\vert\, \grad{}{}{v} \in L^2, \; v = \bar{u} \; \text{in} \; \latInf \setminus \lat \,\},
\end{equation}
where $\clX \subseteq \latInf$ is some generic set. In the following we let $\clV = \clV(\latInf)$.

\begin{figure}[t]
 \centering
 \includegraphics[width=0.95\textwidth]{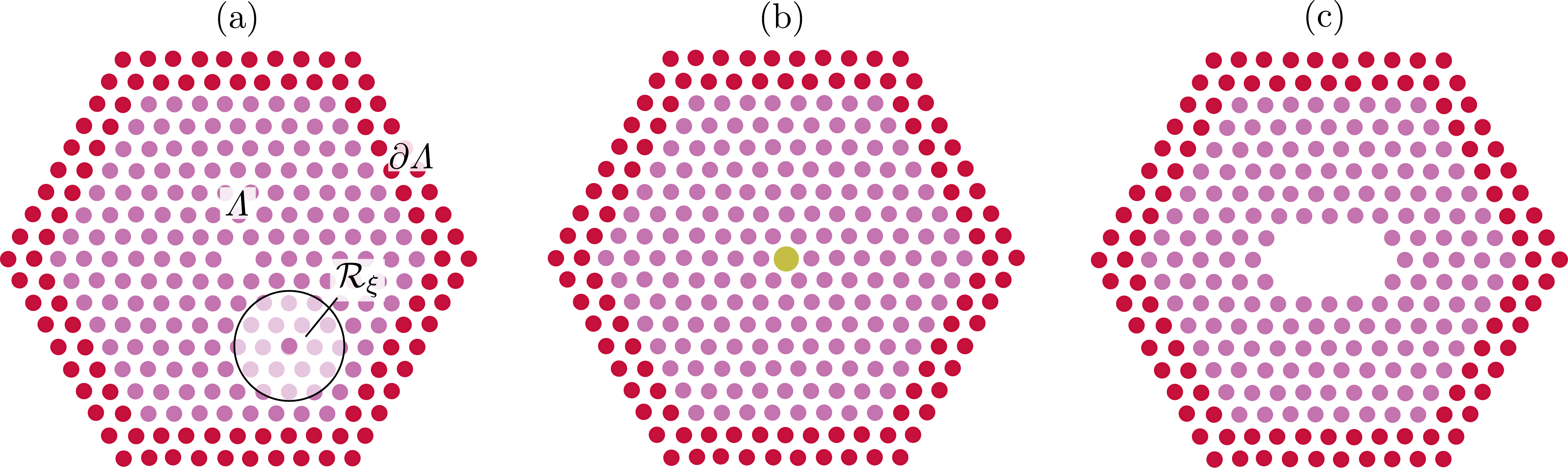}
 \caption{Examples for possible atomistic domains containing crystalline defects, such as in (a) a vacancy, in (b) a substitutional atom, or in (c) a microcrack}
 \label{fig:domains}
\end{figure}

\textit{%
We emphasize that the construction of $\clV$ includes infinite problems, i.e., when $\lat \rightarrow \latInf$. In this case the optimization problem \eqref{eq:atom.ato_problem} is not computable. However, the coupled problem we propose in Section \ref{sec:acc}, which is split into a finite atomistic and a possibly infinite continuum problem, will be---since solutions to the continuum problem can be obtained using Green function techniques (cf. Section \ref{sec:sinc.H_k+1}).
}%

\begin{rem}
 When considering defects with infinite energy, such as dislocations, and the computational domain $\lat$ is also infinite, the final solution is not in $\clV$. One solution procedure is to impose an initial guess $\displ_0 : \lat \rightarrow \real^d$ and subsequently solve the atomistic problem for a corrective (finite-energy) displacement $\displ \in \clV$ (cf. \citep{ehrlacher_analysis_2016}) such that the final solution is given by $\displ_0 + \displ$. This solution procedure is also possible within the present framework but not explicitly accounted for to avoid further notation which is not strictly necessary to describe the details of the proposed method. 
\end{rem}
\end{modified}%

Every atom $\ato \in \lat$ has a site energy \Eato. We assume that \Eato depends on the displacement of atom \ato relative to all other atoms within its interaction range \intRgAto which may extend over a few lattice spacings (usually second- or third-nearest neighbor interactions for metals). This renders the atomistic model nonlocal, but short-range. We write this dependency as $\{ \displ(\atoB) - \displ(\ato) \}_{\atoB \in \intRgAto \setminus \ato} \equiv \{ \displ(\atoB) - \displ(\ato) \}$ such that $\Eato = \Eato(\{ \displ(\atoB) - \displ(\ato) \})$.

The total energy of the system reads
\begin{equation}\label{eq:atom.Etot}
 \Etot(\displ) = \Etot_0 + \Etot_\mrm{int}(\displ) + \Etot_\mrm{ext}(\displ),
\end{equation}
\begin{modified}%
where $\Etot_0$ is the equilibrium cohesive energy.
\end{modified}%
For convenience we assert that $\Etot_0 = 0$ in the following. The internal and external contributions are defined as
\begin{align}
    \Etot_\mrm{int}(\displ) =   \sum_{\ato \in \lat} \Eato(\{ \displ(\atoB) - \displ(\ato) \}),
 && \Etot_\mrm{ext}(\displ) = - \sum_{\ato \in \lat} \force_\mrm{ext}(\ato) \cdot \displ(\ato),
\end{align}
where $\force_\mrm{ext} \in \clV^\ast$ is an external force with compact support. We point out that the formulation of the method does not exclude any specific type of interatomic interaction per se. That is, the definition of the internal energy $\Etot_\mrm{int}$ accounts for the entire class of interatomic many-body potentials.

In this work, attention is drawn to (quasi-)static problems.
\begin{modified}%
That is, starting from some initial guess, we seek for \emph{local solutions} $\displ \in \clV$ of the optimization problem
\end{modified}%
\begin{equation}\label{eq:atom.ato_problem}
 \displ := \arg{\left\{ \, \underset{v \in \clV}{\min} \, \Etot(v) \, \right\}}.
\end{equation}
Solutions to \eqref{eq:atom.ato_problem} solve an Euler-Lagrange equation, subject to the prescribed boundary conditions on $\bnd{}{}{\lat}$, i.e,
\begin{equation}\label{eq:atom.euler-lagrange}\boxed{
 \left\{
 \begin{aligned}
  \; \L[\displ] &= \force_\mrm{ext} &\qquad& \text{in} \; \lat, \\
  \; \displ     &= \bar{\displ}     &      & \text{on} \; \bnd{}{}{\lat},
 \end{aligned}
 \right.
 }
\end{equation}
where the nonlinear operator \L is defined as
\begin{equation}
 \begin{aligned}
  \L : \clV &\rightarrow \clV^\ast \\
       v    &\mapsto     \L[v] \qquad \text{such that} \; \forall\,\ato\in\lat \qquad \L[v](\ato) = \var{}{}{\Etot_\mrm{int}}(\ato), 
 \end{aligned}
\end{equation}
where $\var{}{}{\Etot_\mrm{int}}(\ato)$ is the functional derivative of $\Etot_\mrm{int}(\displ)$ with respect to $\displ$ at \ato. In the ground state, i.e., in the absence of external forces, we have $\var{}{}{\Etot_\mrm{int}}(0) = 0$.
In addition, we require the usual strong stability conditions on the minimizers \displ such that
\begin{equation}\label{eq:atom.stab}
 \forall\,v \in \clV \setminus 0 \qquad \langle \var{}{2}{\Etot_\mathrm{int}}(\displ)[v], v \rangle > 0.
\end{equation}
If \eqref{eq:atom.stab} holds, it is easy to see that solutions to \eqref{eq:atom.euler-lagrange} also solve \eqref{eq:atom.ato_problem}.


\section{Continuum model}
\label{sec:con}

\begin{modified}%
We will now introduce the continuum model as a linearization of the atomistic site energy.%
\footnote{%
\begin{modified}%
The convention of using term ``continuum model'' is not fully precise. In fact, we never use a true continuum model but rather a \emph{discretized} continuum model
\end{modified}%
}
The motivation behind such an approximation stems from the fact that nonlinearities usually only occur in some localized parts of the domain, typically in the vicinity of defects, such that we can replace the atomistic model in the remainder of the domain with a significantly cheaper continuum model. In this work we employ a \emph{local continuum model} which has been successfully used in many variants of atomistic/continuum coupling methods, e.g., the coupled/atomistic discrete dislocation method \citep{shilkrot_multiscale_2004,anciaux_coupled_2018} or the coupling of length scales method \citep{rudd_concurrent_2000}, and is a legitimate choice if the solution in the continuum domain is smooth and/or the interatomic potential is not strongly nonlocal.

In the following we assume a \emph{defect-free continuum model}. Hence, for the sake of deriving the continuum model, we assume in the remainder of this section that every atom has the same (perfect) environment and the site energy $\Eato$ is therefore \emph{independent} of $\ato$.
\end{modified}%

Our continuum model is based on a linearization around a uniformly deformed state $\latInf + \displ_\mrm{F}$, where
\begin{equation}\label{eq:con.displ_hom}
 \forall\,\ato \in \latInf \qquad \displ_\mrm{F}(\ato) = F\ato + C, \qquad \text{with} \; F,C \ge 0.
\end{equation}
Consider now $\displ = \displ_\mrm{F} + \displ'$, where $\displ'$ is a perturbation of $\latInf + \displ_\mrm{F}$.
\begin{modified}%
A Taylor expansion of $\Eato$ to second order around $\displ_\mrm{F}$ then yields the nonlocal harmonic site energy
\begin{equation}\label{eq:con.Etot_hnl}
   \clE_{\mrm{hnl},\ato}(\{ \displ(\atoB) - \displ(\ato) \})
 = \frac{1}{2} \sum_{\atoB \in \intRg_{\ato}} K_\mrm{hnl}(\ato - \atoB) \cdot \Big( \big(\displ_\mrm{F}(\ato) + \displ'(\ato)\big) \otimes \big(\displ_\mrm{F}(\atoB) + \displ'(\atoB)\big) \Big),
\end{equation}
where
$K_\mrm{hnl}(\ato - \atoB) = \sum_{\atoC \in \intRgAto}
\frac{\partial^2 \Eato(\{ \displ(\atoC) - \displ(\ato) \})}
     {\partial\displ(\ato)\partial\displ(\atoB)}$
is the interatomic force constant tensor. It should be noted that this approximation neglects the spurious drift term that arises from the first order term if the deformation gradient is not uniform, e.g., in the case of rotations, and is therefore only valid for \emph{small (nonuniform) deformations}. 

In addition, if the perturbation remains close to homogeneous, we can also consider a linearization of $\displ'$. This requires a more concrete definition of the gradient of $\displ$. To that end, we partition the ideal lattice \eqref{eq:atom.bravais} into a periodic set of simplexes and define an interpolant $\phi_{\ato} \in W^{1,2}(\real^d)$ with compact support on a set $\intRgAto^\mrm{h} \subset \intRgAto$ which spans to the nearest neighbors in the adjacent simplexes. The displacement $\displ$ and its gradient can then be defined $\forall\,x\in\real^d$ as
\begin{align}\label{eq:con.def_u+ugrad}
 \displ(x) = \sum_{\xi \in \intRgAto^\mrm{h}} \phi_{\ato}(x)\displ(\ato),
 &&
 \grad{}{}{\displ}(x) = \sum_{\xi \in \intRgAto^\mrm{h}} \displ(\ato) \otimes \grad{}{}{\phi_{\ato}}(x).
\end{align}
\end{modified}%

Having $\displ(x)$ and $\grad{}{}{\displ}(x)$ well-defined, we then apply a Taylor expansion to $\displ'$ and neglect higher gradients such that
\begin{equation}\label{eq:con.cbh}
 \displ'(\atoB) \approx \displ'(\ato) + \grad{}{}{\displ'}(\ato)(\atoB - \ato),
\end{equation}
which is usually referred to as the Cauchy-Born hypothesis.
\begin{modified}%
Using the approximation \eqref{eq:con.cbh} in \eqref{eq:con.Etot_hnl} yields the classical definition of the local harmonic energy (cf. \citep{hodapp_lattice_2019})
\begin{equation}\label{eq:con.Etot_h}
   \clE_{\mrm{h},\ato}(\{ \displ(\atoB) - \displ(\ato) \})
 = \frac{1}{2} \, \bstiffMat \cdot (\grad{}{}{\displ(\ato)} \otimes \grad{}{}{\displ(\ato)})
 \overset{\eqref{eq:con.def_u+ugrad}}{=}
 \frac{1}{2} \sum_{\atoB \in \intRgAto^\mrm{h}} K_\mrm{h}(\ato - \atoB) \cdot \Big( \big(\displ_\mrm{F}(\ato) + \displ'(\ato)\big) \otimes \big(\displ_\mrm{F}(\atoB) + \displ'(\atoB)\big) \Big),
\end{equation}
where $\bstiffMat$ is the usual fourth-order material stiffness tensor. In the term on the right hand side $K_\mrm{h}(\ato - \atoB)$ can be considered as a local version of $K_\mrm{hnl}(\ato - \atoB)$ which depends on the precise choice of $\phi_{\ato}$.

Using $K_\mrm{h}(\ato - \atoB)$, we then define the linearized version of $\L$ as follows
\begin{equation}\label{eq:con.Lh=fext}
 \forall\,\ato\in\latInf \qquad \Lh[\displ](\ato) =
 \sum_{\atoB \in \intRgAto^\mrm{h}} K_\mrm{h}(\ato - \atoB) \displ(\atoB).
\end{equation}
\end{modified}%


\section{Force-based coupled atomistic/continuum problem}
\label{sec:acc}

\begin{modified}%
When the computational domain becomes sufficiently large, e.g., in order to avoid spurious elastic interactions of the defect with the far-field boundary, solving the fully atomistic problem \eqref{eq:atom.ato_problem} can be too costly. This motivates the introduction of an atomistic/continuum (A/C) coupling scheme which approximates the full problem by restricting atomic resolution to some small part of the computational domain, containing the defect, surrounded by the significantly cheaper defect-free continuum elasticity model.

In this work we use a force-based A/C coupling scheme (e.g. \citep{kohlhoff_new_1989}) which is described in the following. To that end, we first decompose the computational domain $\lat$ into an atomistic domain $\latA \subset \lat$ and a continuum domain $\latC \subset \lat$ such that $\lat := \latA \cup \latC$, as shown in Figure \ref{fig:domain_decomp}. In practical applications we usually have that
\begin{equation}
 \mathrm{diameter}(\latA) \ll \mathrm{diameter}(\latC),
\end{equation}
meaning that the atomistic domain is well-separated from the far-field boundary.

\begin{figure}[hbt]
 \centering
 \includegraphics[width=0.5\textwidth]{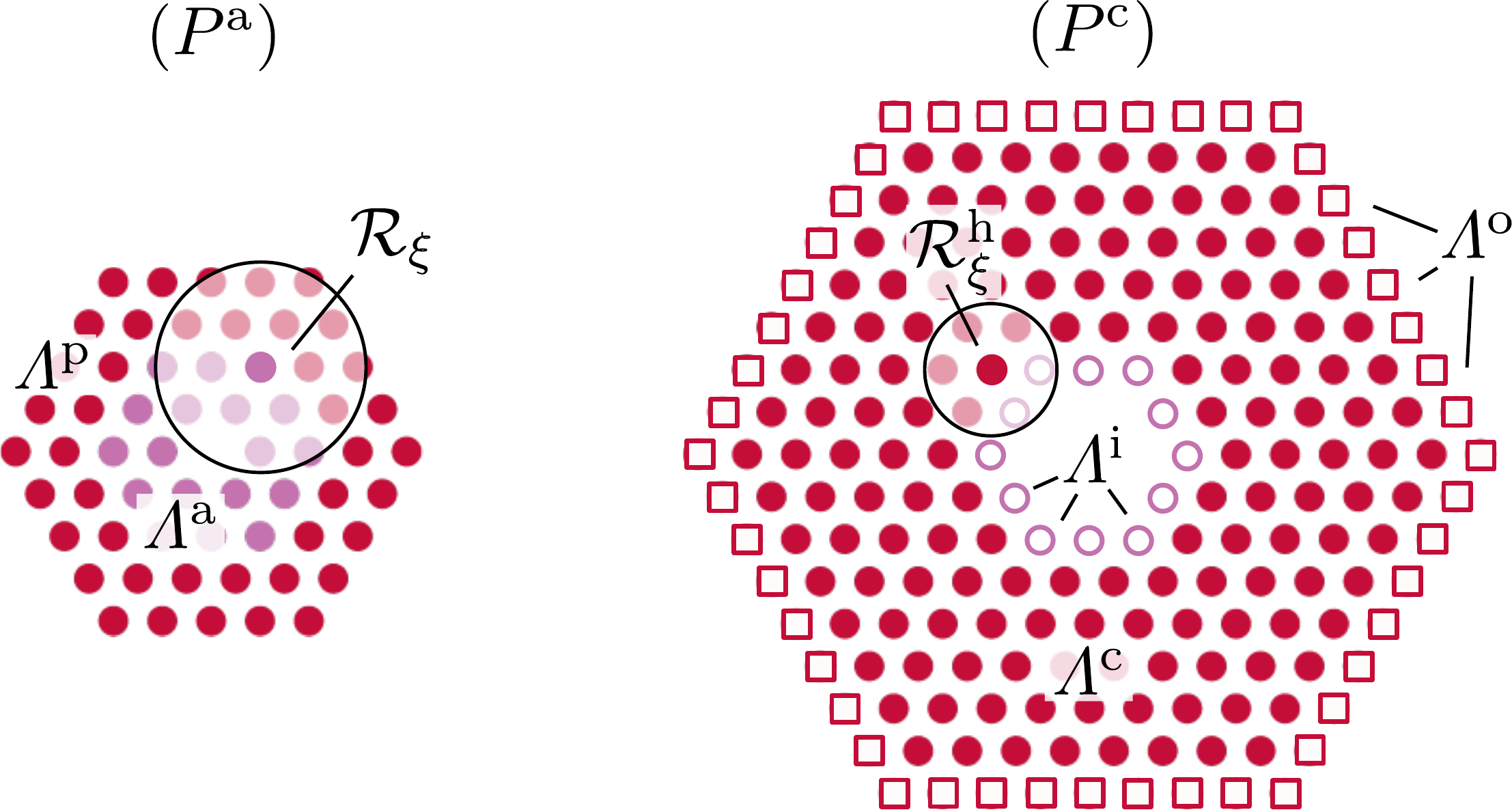}
 \caption{Decomposition of the computational problem into an atomistic problem $(P^\rma)$ and a continuum problem $(P^\rmc)$}
 \label{fig:domain_decomp}
\end{figure}

Atoms in $\latA$ are nonlocal and interact with continuum nodes in $\latP \subset \latC$ according to their interaction range $\intRgAto$. Vice versa, the continuum problem is assumed to be linear elastic, as described in the previous section, such that its interaction range $\intRgAto^\mrm{h}$ spans up to first nearest neighbors (cf. Figure \ref{fig:domain_decomp}). Therefore, continuum nodes in $\latC$ ``interact'' with atoms from the interface domain $\latI$. In addition, for finite problems, the outer boundary of the continuum problem prescribing the far-field boundary condition is denoted by $\lat^\I$.

Mathematically speaking, by using a linear elastic continuum model we assume that the solution in $\latC$ is sufficiently regular in the sense that higher gradients do not dominate the coupling error. We will make this more precise in Section \ref{sec:conv.rate} where we quantify the modeling error of the coupling scheme.

\emph{%
Moreover, it is emphasized that Figure \ref{fig:domain_decomp} is only a schematic and we do not impose any particular restrictions on the shape of the domains beyond the previously stated assumptions.%
}
\end{modified}%

With the above notation, we now state the coupled problem as follows: find $\displ \in \clV$ such that
\begin{empheq}[box=\fbox]{align}\label{eq:acc.coupled_problem}
 (P^\a) \; \left\{
 \begin{aligned}
  \; \L[\{\displ^\a,\displ^\p\}] &= \force_\mrm{ext} &\qquad& \text{in} \; \lat^\a, \\
  \; \displ                      &= \displ^\c        &      & \text{in} \; \lat^\p,
 \end{aligned}
 \right.
 &&
 (P^\c) \; \left\{
 \begin{aligned}
  \; \Lh[\{\displ^\c,\displ^\i,\displ^\I\}] &= \force_\mrm{ext}  &\qquad& \text{in} \; \lat^\c, \\
  \; \displ                                 &= \displ^{\indDomA} &      & \text{on} \; \lat^\i, \\
  \; \displ                                 &= \bar{\displ}      &      & \text{on} \; \lat^\I.
 \end{aligned}
 \right.
\end{empheq}

The distinct advantage of force-based methods is the consistency of the coupling, i.e., no spurious forces arise in the vicinity of the artificial interface due to the nonlocal-local mismatch between the two models. Their disadvantage is, however, the lack of a well-defined energy functional. This restricts the choice of possible monolithic solvers to multidimensional root-finding methods, e.g., the generalized minimal residual method \citep{dobson_iterative_2011}.

Therefore, another popular choice are domain decomposition solvers based on the alternating Schwarz method which iteratively solve two energy minimization problems in $\lat^\a$ and $\lat^\c$, bypassing the concurrent coupling (e.g., \citep{parks_connecting_2008,li_atomistic-based_2012,pavia_parallel_2015}). However, they usually converge very slowly (potentially even slower than a model with full atomistic resolution), even with sophisticated acceleration techniques, such as overlapping subdomains, which makes A/C coupling counterproductive. The aim of this work is thus to develop a new domain decomposition solver with improved convergence behavior over these existing methods.

\section{Sinclair method for bounded problems}
\label{sec:sinc}

\begin{modified}%
In the 1970s, Sinclair and coworkers introduced a fast alternative to the alternating Schwarz method to solve atomistic problems which are embedded in an effectively infinite domain. In \citep{hodapp_lattice_2019}, we have shown that the excellent convergence properties of the Sinclair method are due to the particular splitting of the coupled operator into a \emph{local finite anharmonic} and a \emph{global \textbf{infinite} harmonic} problem.

In Section \ref{sec:sinc.description} we will generalize this splitting procedure to bounded problems. In comparison with the original Sinclair method, this now requires a \emph{global \textbf{finite} harmonic} problem to be computed. In Section \ref {sec:sinc.H_k+1}, we introduce a method to compute this finite harmonic problem using a superposition of an infinite inhomogeneous and a finite homogeneous contribution. Both contributions will be solved using a \emph{discrete boundary element method}.

This strategy offers two advantages. First, we retain the efficiency in terms of the total number of degrees of freedom of the original Sinclair method by not requiring a spatial discretization. Second, our analysis in Section \ref{sec:conv} can be readily used for bounded and unbounded problems, i.e., by simply switching off the finite homogeneous contribution.

Implementation aspects of integrating the discrete boundary element method into the computation of the coupled problem are discussed in Section \ref{sec:sinc.impl}.
\end{modified}%

\subsection{Description of the method}
\label{sec:sinc.description}

\subsubsection{Iteration equation}
\label{sec:sinc.description.iter_eq}

\begin{modified}%
We now derive the iteration equation for the Sinclair method. To do so, we first require a splitting of the coupled problem \eqref{eq:acc.coupled_problem} into an anharmonic and a harmonic part.

\textit{%
Furthermore, we require that the coupled problem \eqref{eq:acc.coupled_problem} is defined on a \textbf{defect-free} domain $\lat$ since the global harmonic problem is considered to be defect-free. In what follows, we therefore tacitly assume for analysis purposes that any \textbf{incorporation of defects} by adding/removing atoms from $\lat^\a$ will be instead achieved by modifying the atomistic energy \eqref{eq:atom.Etot} (that is, e.g., by switching off interactions with a ``vacancy atom'').
}%
\end{modified}%

\begin{defn}[\textbf{Anharmonic/harmonic operator split}]\label{defn:ah/h_split}
 Let $\Lcpl$ be the differential operator associated with the coupled problem \eqref{eq:acc.coupled_problem}. We then denote the anharmonic/harmonic operator split as the additive decompositions
 \begin{align}\label{eq:sinc.Lcpl_split}
  \Lcpl = \Lah + \Lh, && \displ &= \displ_\indAHarm + \displ_\indHarm
 \end{align}
 into anharmonic parts $\Lah,\displ_\indAHarm$, and harmonic parts $\Lh,\displ_\indHarm$, with
 \begin{align}\label{eq:sinc.Lah_&_uah}
  \Lah =
  \left\{\,
  \begin{aligned}
   \L - \Lh & \qquad \text{in} \; \latA, \\
   \nullOp  & \qquad \text{else},
  \end{aligned}
  \right.
  &&
  \displ_\indAHarm = 0 \qquad \text{in} \; \latInf \setminus \latA.
 \end{align}
\end{defn}

\begin{modified}%
We remark here that this definition is not well-defined in the sense that $\displ_\indAHarm$ and $\displ_\indHarm$ are \emph{not unique} in $\lat^\a$. The splitting $\displ = \displ_\indAHarm + \displ_\indHarm$ is used here as a tool which is exploited below to derive the iteration equation. In fact, as we shall see later on in Section \ref{sec:sinc.description}, we never need to compute $\displ_\indAHarm$ and $\displ_\indHarm$ in $\lat^\a$ in practice. 
\end{modified}%

Using \eqref{eq:sinc.Lcpl_split} from Definition \ref{defn:ah/h_split} and exploiting the linearity of the harmonic operator, we can write the coupled problem \eqref{eq:acc.coupled_problem} as follows
\begin{equation}
 \Lcpl[\displ] = \Lah[\displ] + \Lh[\displ_\indAHarm] + \Lh[\displ_\indHarm] = \force_\mathrm{ext} \qquad \text{in} \; \lat,
\end{equation}
omitting the natural boundary conditions.
In matrix notation we can write this
\begin{equation}\label{eq:sinc.Lcpl_matrix}
 \begin{pmatrix} \L^{\a|}[\displ] - \Lh^{\a|}[\displ] \\ 0 \end{pmatrix}
 +
 \begin{pmatrix} \Laa & \nullOp^{\a|\cb} \\ \Lca & \nullOp^{\c|\cb} \end{pmatrix}
 \begin{bmatrix} \displ_\indAHarm^\indDomA \\ 0 \end{bmatrix}
 +
 \begin{pmatrix} \Laa & \Lh^{\a|\cb} \\ \Lca & \Lh^{\c|\cb} \end{pmatrix}
 \begin{bmatrix} \displ_\indHarm^\indDomA \\ \displ^\cb \end{bmatrix}
 =
 \begin{pmatrix} \force_\mathrm{ext}^\indDomA \\ \force_\mrm{ext}^\c \end{pmatrix},
\end{equation}
where the index ``\cb'' refers to functions defined on the domain $\latC \cup \lat^\I$ (where $\displ_\indAHarm = 0$ according to \eqref{eq:sinc.Lah_&_uah}).

After re-arranging some of the terms and moving $\force_\mathrm{ext}$ to the left hand side of the equation, \eqref{eq:sinc.Lcpl_matrix} reads
\begin{equation}\label{eq:sinc.op_split}
 \underbrace{
  \begin{pmatrix} \L^{a|}[\displ] \\ 0 \end{pmatrix}
  - \begin{pmatrix} \force_\mathrm{ext}^\indDomA \\ 0 \end{pmatrix}
 }_{\textbf{\normalsize{(AH)}}}
 +
 \underbrace{
  \begin{pmatrix} \Laa & \Lh^{\a|\cb} \\ \Lca & \Lh^{\c|\cb} \end{pmatrix}
  \begin{bmatrix} \displ_\indHarm^\indDomA \\ \displ^\cb \end{bmatrix}
  - \begin{pmatrix}
     \Laa[\displ_\indHarm] + \Lh^{\a|\cb}[\displ] \\ 
     -\Lca[\displ_\indAHarm] + \force_\mathrm{ext}^\indDomC
    \end{pmatrix}
 }_{\textbf{\normalsize{(H)}}}
 = \begin{pmatrix} 0 \\ 0 \end{pmatrix}.
\end{equation}
Here, the parts labeled as \textbf{(AH)} and \textbf{(H)} are referred to as the anharmonic and the harmonic problem, respectively. They are coupled through the solution in \latP and \latI and will be computed separately.

To solve \eqref{eq:sinc.op_split}, a staggered procedure is proposed which iterates between \textbf{(AH)} and \textbf{(H)}. Therefore, let us denote $k \in \nat$ as the global iteration index and fix an initial guess $\displ_0$. The $k$+1-th iteration then reads
\begin{empheq}[box=\fbox]{align}\label{eq:sinc.sinc}
 \textbf{(AH)}_{k+1} \;
 \left\{
 \begin{aligned}
  \; \L[\{ \displ^\a_{k+1}, \displ^\p_{k+1} \}] &= \force_\mrm{ext} &\;& \text{in} \; \lat^\a, \\
  \; \displ_{k+1}                               &= \displ^\c_k      &  & \text{in} \; \lat^\p,
 \end{aligned}
 \right.
 &&
 \textbf{(H)}_{k+1} \;
 \left\{
 \begin{aligned}
  \; \Lh[\displ_{\mrm{h},k+1}] &= \L^{|\a}[\displ_{\indHarm,0}] + \Lh^{|\cb}[\displ_0]      &\;& \text{in} \; \latA, \\
  \; \Lh[\displ_{\mrm{h},k+1}] &= -\Lh^{|\a}[\displ_{\indAHarm,k+1}] + \force_\mrm{ext} &\;& \text{in} \; \lat^\c, \\
  \; \displ                    &= \bar{\displ}                                            &  & \text{on} \; \lat^\I.
 \end{aligned}
 \right.
\end{empheq}
Note that \textbf{(AH)} is a finite problem, defined in $\lat^\a$, and \textbf{(H)} is a global problem, defined in the entire domain $\lat := \lat^\a \cup \lat^\c$.

At first glance, it seems natural to choose the global harmonic solution as an initial guess. However, there are situations where an initial guess is not available or linear elasticity is a bad approximation of the coupled problem. Fortunately, in such cases we can initialize the algorithm with $\displ_0 = 0$ without actually influencing the asymptotic convergence rate. This will be demonstrated, i.a., in Section \ref{sec:conv}.

\subsubsection{Treatment of the source term \texorpdfstring{$\Lca[\displ_{\indAHarm,k+1}]$}{}}
\label{sec:sinc.source_term}

In practice, one would like to avoid to work with the source term $\Lca[\displ_{\indAHarm,k+1}]$ of $\textbf{(H)}_{k+1}$, taking into account that $\displ_{\indHarm,k+1}^\a$ needs to be evaluated in the \emph{entire} atomistic domain to compute ${\displ^\a_{\indAHarm,k+1} = \displ^\a_{k+1} - \displ^\a_{\indHarm,k+1}}$ explicitly. Fortunately, it is possible to solve $\textbf{(H)}_{k+1}$ using solely the total solution $\displ$. To see this, consider the solution to the harmonic problem in step $k$+1
\begin{equation}\label{eq:sinc.displ_upd}
   \begin{pmatrix} \displ^\a_{\mrm{h},k+1} \\ \displ^\c_{k+1} \end{pmatrix}
 = \inv{\begin{pmatrix} \Lh^\aa & \Lh^\ac \\ \Lh^\ca & \Lh^\cc \end{pmatrix}}
   \begin{bmatrix} \force_\mrm{ext}^\a \\ - \Lh^\ca[\displ_{\mrm{ah},k+1}] - \Lh^\cI[\displ] + \force_\mrm{ext}^\c \end{bmatrix}.
\end{equation}
By subtracting $\displ_{\mrm{h},k}$ from $\displ_{\mrm{h},k+1}$ we obtain
\begin{equation}\label{eq:sinc.displ_upd2}
   \begin{pmatrix} \displ^\a_{\mrm{h},k+1} \\ \displ^\c_{k+1} \end{pmatrix}
 = \begin{pmatrix} \displ^\a_{\mrm{h},k} \\ \displ^\c_k \end{pmatrix}
 - \inv{\begin{pmatrix} \Lh^\aa & \Lh^\ac \\ \Lh^\ca & \Lh^\cc \end{pmatrix}}
   \begin{bmatrix} 0 \\ \Lh^\ca[\Delta\displ_{\mrm{ah},k}] \end{bmatrix},
\end{equation}
where $\Delta\displ_{\mrm{ah},k} = \displ_{\mrm{ah},k+1} - \displ_{\mrm{ah},k}$. To obtain an expression of $\Lh^\ca[\Delta\displ_{\mrm{ah},k}]$ in terms of the full solution $\displ$ only, consider
\begin{equation}\label{eq:sinc.Lca_u_ah_k}
 \Lh^\ca[\displ_{\mrm{ah},k}] = - \Lh^\ca[\displ_{\mrm{h},k}] - \Lh^\cc[\displ_{k}] - \Lh^\cI[\displ] + \force_\mrm{ext}^\c,
\end{equation}
which can be directly obtained from \eqref{eq:sinc.displ_upd}. Now subtract \eqref{eq:sinc.Lca_u_ah_k} from $\Lh^\ca[\displ_{\mrm{ah},k+1}]$ leading to 
\begin{equation}\label{eq:sinc.finh}
 \begin{aligned}
     \Lh^\ca[\Delta\displ_{\mrm{ah},k}]
   = \force^\c_{\mathrm{inh},k+1}
  &= \Lh^\ca[\displ_{\mrm{ah},k+1}] + \Lh^\ca[\displ_{\mrm{h},k}] + \Lh^\cc[\displ_{k}] + \Lh^\cI[\displ] - \force_\mrm{ext}^\c \\
  &= \Lh^\ca[\displ_{k+1}] + \Lh^\cc[\displ_{k}] + \Lh^\cI[\displ] - \force_\mrm{ext}^\c \\
  &= \Lh^\ca[\displ_k + \Delta\displ_k] + \Lh^\cc[\displ_{k}] + \Lh^\cI[\displ] - \force_\mrm{ext}^\c \\
  &= \Lh^\ca[\Delta\displ_k],
 \end{aligned}
\end{equation}
which we denote as the \emph{inhomogeneous force} $f^\c_{\mathrm{inh},k+1}$ in the $k$+1-th iteration (cf. \citep{hodapp_lattice_2019}). In the latter expression only the full solution appears and, therefore, $\Lh^\ca[\Delta\displ_{\mrm{ah},k}]$ can be evaluated conveniently after the anharmonic problem has been solved.
\textit{%
Since the continuum model is local, $\force^\c_{\mathrm{inh},k+1}$ will be nonzero only in the domain $\latIpl$, i.e., the layer of nodes which are nearest neighbors of the interface atoms $\latI$ (see {\normalfont\textbf{(H)}} in Figure \ref{fig:domain_harmonic}). So, we just need to compute $\force^\ip_{\mathrm{inh},k+1}$ in practice.
}%

\begin{rem}
 The name ``inhomogeneous force'' was coined in the original works by Sinclair (e.g., \citep{sinclair_flexible_1978}). The origin of this name stems from the fact that updating the atomistic solution generates a mismatch between both models---which vanishes upon convergence. The idea of Sinclair was thus to update the displacements in $\lat^\c$ corresponding to a force which counteracts $\Lh^\ca[\Delta\displ_{\mrm{ah},k}]$; this is why the minus sign appears on the right hand side of \eqref{eq:sinc.displ_upd2}.
\end{rem}

In principle, the harmonic problem can now be solved with any conventional finite element method. However, volume-based methods require a very fine discretization and can thus become significantly more expensive than the atomistic problem itself. For the class of A/C coupling problems, boundary element techniques seem preferable since the full solution in $\lat^\c$ is usually not explicitly required.

\begin{moved}%
\subsection{Computation of \texorpdfstring{$\normalfont\textbf{(H)}$}{} using a discrete boundary element method}
\label{sec:sinc.H_k+1}

Boundary element methods for lattice problems have been proposed by \citet{martinsson_fast_2002} for the discrete Laplace equation and in the context of A/C coupling by \citet{li_atomistic-based_2012} and \citet{hodapp_lattice_2019}. By using a \emph{discrete boundary element method}, we do not have to resort back to a ``true'' continuum model and could consider arbitrary interaction stencils, i.e., all subsequent developments generically apply to local and nonlocal elasticity models. Moreover, it is worthwhile noting that fast summation techniques, such as the fast multipole method \citep{greengard_fast_1987} or hierarchical matrices \citep{tyrtyshnikov_mosaic-skeleton_1996,hackbusch_sparse_1999}, readily apply to lattice problems as demonstrated in \citep{martinsson_boundary_2009,hodapp_lattice_2019}.

\begin{movedAndModified}%
\subsubsection{Solution procedure}
\label{sec:sinc.H_k+1.sol_proc}

In the following we will drop the subscripted iteration index $k$+$1$ and refer to \textbf{(H)} as some harmonic problem with generic right hand side $r$ and denote its solution by $\displ_\indHarm$. Since $\textbf{(H)}$ is a linear problem we can decompose it into two contributions as illustrated in Figure \ref{fig:domain_harmonic}: an infinite inhomogeneous problem $\textbf{(H$_\infty$)}$ and a finite homogeneous problem $\textbf{(H$_{\rm hom}$)}$ which corrects \textbf{(H$_\infty$)} for the prescribed boundary condition on $\lat^\I$. The solution to $\textbf{(H)}$ is then obtained by superimposing the solutions of $\textbf{(H$_\infty$)}$ and $\textbf{(H$_{\rm hom}$)}$. Again, the idea behind this decomposition is the possibility to switch conveniently between unbounded and bounded problems since both classes of problems are important in terms of applications.

\begin{figure}[hbt]
 \centering
 \includegraphics[width=0.95\textwidth]{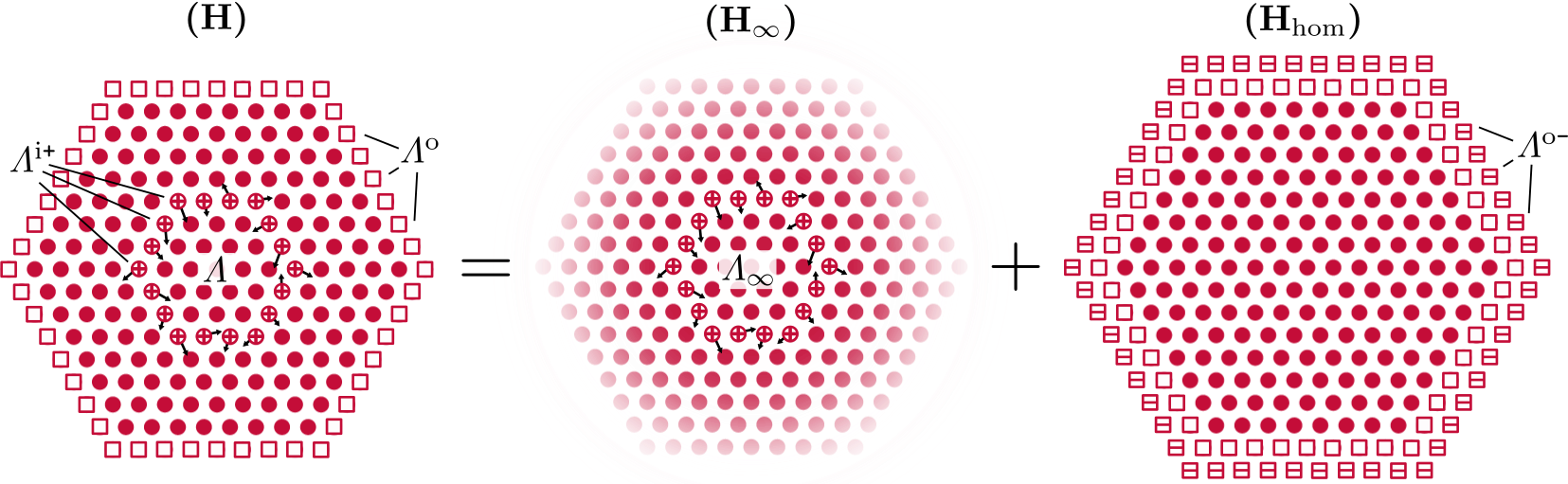}
 \caption{Decomposition of the harmonic problem \textbf{(H)} into an infinite inhomogeneous contribution \textbf{(H$_\infty$)} and a finite homogeneous contribution \textbf{(H$_{\rm hom}$)} (only the forces on the $\latIpl$ nodes are visualized)}
 \label{fig:domain_harmonic}
\end{figure}

To that end, given some $r$, we are looking for solutions $\displInf$ and $\displHom$ such that
\begin{empheq}[box=\fbox]{align}\label{eq:con.split}
 \normalfont\textbf{(H$_\infty$)} \quad \left.
 \begin{aligned}
  \; \Lh[\displInf] &= r &\qquad& \text{in} \; \latInf,
 \end{aligned}
 \right.
 &&
 \normalfont\textbf{(H$_{\rm hom}$)} \; \left\{
 \begin{aligned}
  \; \Lh[\displHom] &= 0                        &\qquad& \text{in} \; \lat, \\
  \; \displHom      &= \bar{\displ} - \displInf &      & \text{on} \; \lat^\I. \\
 \end{aligned}
 \right.
\end{empheq}
The full solution is then given $\forall\,\ato \in \lat$ as $\displ_\indHarm(\ato) = \displInf(\ato) + \displHom(\ato)$.

\subsubsection{Computation of \texorpdfstring{$\normalfont\textbf{(H$_\infty$)}$}{} using the lattice Green function}
\label{sec:sinc.H_k+1.inh}

In order to solve the infinite problem $\normalfont\textbf{(H$_\infty$)}$ we first briefly define the \emph{lattice Green operator} $\greenOp$, the inverse of $\Lh$, of which we will make heavily use of in the following sections.

The lattice Green operator can be formally defined for all $v : \latInf \rightarrow \real^d$ such that
\begin{equation}\label{eq:con.GLh=I}
 \begin{aligned}
  \forall\,\ato\in\latInf \qquad
  (\greenOp\Lh)[v](\ato)
  &= \sum_{\atoB \in \latInf} \left( \sum_{\{\, \atoC \in \latInf \,\vert\, \atoC-\atoB \in \intRg_{\atoC}^\rmh \,\}} \green(\ato - \atoC) K_\rmh(\atoC - \atoB) \right) v(\atoB) \\
  &= \clI[v](\ato) \; \left(\equiv v(\ato)\right),
 \end{aligned}
\end{equation}
where the kernel of $\greenOp$ is denoted as the Green function $\green(\ato - \atoC)$ and $\clI$ is the identity operator. One possibility to obtain $\green$ is by means of semi-discrete Fourier transforms (see, e.g., \citep{thomson_lattice_1992}) but we are not concerned with this aspect here and assume in the following that $\green$ is given.

Using $\greenOp$, the solution of $\normalfont\textbf{(H$_\infty$)}$ can be written as
\begin{equation}\label{eq:sinc.Hinf}
 \displInf = \greenOp[r] \qquad \text{in} \; \latInf.
\end{equation}

\subsubsection{Computation of \texorpdfstring{$\normalfont\textbf{(H$_{\rm hom}$)}$}{} using the boundary summation equation}
\label{sec:sinc.H_k+1.hom}

In order to obtain the solution to the homogeneous problem $\normalfont\textbf{(H$_\mrm{hom}$)}$ we will write $\displHom$ solely in terms of degrees of freedom on the outer boundary $\lat^\I$. This procedure can be considered as the analog to the ``integration by parts'' in continuum mechanics leading to a \emph{boundary summation equation}, a discrete variant of the boundary integral equation, in terms of the boundary displacement and forces (as opposed to displacements and stresses in the case of the boundary integral equation). But instead of multiplying the Euler-Lagrange equation with a test function and integrating by parts twice, we apply the Green operator and then use the relation \eqref{eq:con.GLh=I}.

To derive such an expression for $\displHom$ we closely follow the derivation in \citep{hodapp_lattice_2019} (cf. section 2.5 therein). In the first step, we define the harmonic problem on $\bar{\lat} = \lat \cup \lat^\I$. That is, we assume some $v \in \clV$ which solves
\begin{equation}\label{eq:con.Lh_u_extnd}
 \Lh[v] = 0 \qquad \text{in} \; \bar{\lat}.
\end{equation}
Using the ``superscript notation'' (cf. Section \ref{sec:notation.discrete_problem}), we split $\Lh[v]$ into an inner and an outer contribution as follows
\begin{equation}\label{eq:con.Lh_u_extnd_mat}
 \Lh^{\lb|}[v] = \Lh^{\lb|\lb}[v] + \Lh^{\lb|\Im}[v],
\end{equation}
where the co-domain superscripted by ``$\Im$'' refers to functions defined on the outer layer of nodes $\lat^\Im \subset \latInf \setminus \bar{\lat}$ which interact with the nodes on the outer boundary $\lat^\I$ as shown in Figure \ref{fig:domain_harmonic}. We now apply the part of the lattice Green operator $\greenOp^{\lb|\lb}$, which maps a force in $\bar{\lat}$ to a displacement in $\bar{\lat}$, to \eqref{eq:con.Lh_u_extnd_mat} leading to
\begin{equation}\label{eq:con.GLh_u_extnd_mat}
 (\G^{\lb|\lb}\Lh^{\lb|})[v] = (\G^{\lb|\lb}\Lh^{\lb|\lb})[v] + (\G^{\lb|\lb}\Lh^{\lb|\Im})[v].
\end{equation}

This formulation is still not particular useful since $\G^{\lb|\lb}\Lh^{\lb|\lb}$ \emph{appears to be} a dense operator, meaning that its associate matrix is fully populated due to the long-range nature of the lattice Green function. Fortunately this is not the case and we will rewrite it to reveal its sparsity pattern. To that end, we consider the following decomposition of $\G\Lh$ (eq. \eqref{eq:con.GLh=I}) in matrix notation
\begin{equation}\label{eq:con.GLh=I_2}
 \begin{aligned}
     \G\Lh
  &\equiv \begin{pmatrix} \G^{\lb|\lb}  & \G^{\lb|\re}  \\ \G^{\re|\lb}  & \G^{\re|\re} \end{pmatrix}
     \begin{pmatrix} \Lh^{\lb|\lb} & \Lh^{\lb|\re} \\ \Lh^{\re|\lb} & \Lh^{\re|\re} \end{pmatrix} \\
  &= \begin{pmatrix} \G^{\lb|\lb}\Lh^{\lb|\lb} + \G^{\lb|\re}\Lh^{\re|\lb}  & \G^{\re|\re}\Lh^{\re|\lb}  + \G^{\re|\lb}\Lh^{\lb|\lb} \\ 
                     \G^{\lb|\re}\Lh^{\re|\re}    + \G^{\lb|\lb}\Lh^{\lb|\re} & \G^{\lb|\re}\Lh^{\re|\lb} + \G^{\lb|\lb}\Lh^{\lb|\lb} \end{pmatrix}
   = \begin{pmatrix} \Id^{\lb|\lb} & \nullOp^{\lb|\re} \\ \nullOp^{\re|\lb} & \Id^{\re|\re} \end{pmatrix},
 \end{aligned}
\end{equation}
where the superscripted domain index ``$\rm re$'' refers to functions defined on the outer (infinite) domain, the \emph{remainder} $\lat^\re = \latInf \setminus \bar{\lat}$. Therefore, we can write
\begin{equation}\label{eq:sinc.GlblbLlblb}
 (\G^{\lb|\lb}\Lh^{\lb|\lb})[v] = \Id^{\lb|\lb}[v] - (\G^{\lb|\re}\Lh^{\re|\lb})[v].
\end{equation}
In addition, recall that the operator $\Lh$ satisfies the following relations $\forall\,v$
\begin{align}\label{eq:con.Lrcb_props}
 \textbf{(1)} \qquad
 \Lh^{|\lb}[v] = 0 \qquad \text{in} \; \lat^\re \setminus \lat^\Im,
 &&
 \textbf{(2)} \qquad
 \Lh^{|\l}[v] = 0 \qquad \text{in} \; \lat^\re.
\end{align}
These relations are due to short-range interactions and can be deduced from the construction of the local version of the force constant tensor $K_\mrm{h}$ (cf. Section \ref{sec:con}): \textbf{(1)} Displacements of nodes in $\bar{\lat}$ only exert forces on the remainder nodes $\lat^\Im$; \textbf{(2)} Displacements in $\lat$ do not exert any forces on nodes in the remainder $\lat^\re$ (only on $\lat^\I$). With \eqref{eq:sinc.GlblbLlblb} and \eqref{eq:con.Lrcb_props}, $(\G^{\lb|\lb}\Lh^{\lb|\lb})[v]$ now reads
\begin{equation}\label{eq:con.red1}
 (\G^{\lb|\lb}\Lh^{\lb|\lb})[v]
 = \Id^{\lb|\lb}[v] - (\G^{\lb|\Im}\Lh^{\Im|\lb})[v]
 = \begin{pmatrix}
    \Id^{\I|\I} - \G^{\I|\Im}\Lh^{\Im|\I} & \nullOp^{\I|\l} \\
    \G^{\l|\Im}\Lh^{\Im|\I} & \Id^{\l|\l}
   \end{pmatrix}
   \begin{bmatrix}
    v^\I \\ v^\l
   \end{bmatrix}.
\end{equation}
Hence, $\G^{\lb|\lb}\Lh^{\lb|\lb}$ is obviously at least partially sparse---only the $\lb|\I$ part is dense.

Turning back to \eqref{eq:con.GLh_u_extnd_mat}, we can also rewrite the second term on the right hand side $(\G^{\lb|\lb}\Lh^{\lb|\Im})[v]$. Using the fact that $\Lh^{|\Im}[v] = 0 \; \text{in} \; \lat$, i.e., that displacements of the $\lat^\Im$-nodes do not create forces in the computational domain $\lat$ (which, again, follows from the properties of $K_\mrm{h}$), it follows
\begin{equation}\label{eq:con.red2}
 (\G^{\lb|\lb}\Lh^{\lb|\Im})[v] = (\G^{\lb|\I}\Lh^{\I|\Im})[v].
\end{equation}

Finally, using \eqref{eq:con.red1} and \eqref{eq:con.red2} in \eqref{eq:con.GLh_u_extnd_mat}, we obtain the boundary summation equation which yields the following expression for the displacements in $\bar{\lat}$ after re-arranging some of the terms (cf. \citep{hodapp_lattice_2019})
\begin{equation}\label{eq:con.bse}
 v = \F^{|\I}[v] - \G^{|\I}[\force] \qquad \text{in} \; \latb,
\end{equation}
with $\F^{|\I} = \G^{|\Im}\Lh^{\Im|\I}$ and the boundary forces $\force^\I = \Lh^{\I|\Im}[v]$.

Now, using as boundary condition $v = \displHom = \bar{\displ} - \displInf$ of $\normalfont\textbf{(H$_{\rm hom}$)}$ in \eqref{eq:con.bse}, we can solve it for $\force^\I$, that is, solving the linear system
\begin{equation}\label{eq:sinc.Hom_ls}
 \force^\I = - (\inv{\G^{\I|\I}}(\Id^{\I|\I} - \F^{\I|\I}))[\bar{\displ} - \displInf].
\end{equation}
Having computed $\force^\I$, the solution to $\normalfont\textbf{(H$_{\rm hom}$)}$ in the entire computational domain $\lat$ is obtained via
\begin{equation}\label{eq:sinc.Hom_sol}
 \displHom^\l = \F^{\l|\I}[\bar{\displ} - \displInf] - \G^{\l|\I}[\force].
\end{equation}

\begin{rem}[Treatment of the outer boundary]\label{rem.outer_boundary}
 For large computational domains, especially in three dimensions, the outer boundary contains substantially more degrees of freedom than the inner boundary. Therefore, even with fast summation techniques, it might be beneficial to seek for solutions in a subspace of $\clV(\lat^\I)$. A first step in that direction has been taken by \citet{li_atomistic-based_2012} who proposed a $\mathbb{P}1$ interpolation of $\displ$ over a reduced set of nodes on $\lat^\I$.
\end{rem}
\end{movedAndModified}%
\end{moved}%

\subsection{Implementation}
\label{sec:sinc.impl}

We now formulate the implementation of the Sinclair iteration equation \eqref{eq:sinc.sinc} using the discrete boundary element method for solving the global harmonic problem. The implementation of the discrete boundary element method, as presented in Section \ref{sec:sinc.H_k+1}, is shown by Algorithm \ref{algo:dbem}.

Following Algorithm \ref{algo:sinc}, we first compute the initial guess of the problem, either by using the elastic solution---or zero. In every iteration $k$+1 we solve the anharmonic atomistic problem \textbf{(AH)}$_{k+1}$ and subsequently compute the inhomogeneous force $\force^\ip_{\mathrm{inh},k+1}$ (eq. \ref{eq:sinc.finh}) which arises at the A/C interface due to the mismatch between both models. Then, we recompute the  harmonic problem \textbf{(H)}$_{k+1}$, noting that we only need the solution in the pad domain $\latP$ in order to provide the boundary condition for \textbf{(AH)}$_{k+2}$ (if there are external forces in $\latC$ we also need to account for them---but we do this only once during the first iteration). The previous steps are then repeated until convergence is attained, that is, e.g., when the maximum force on each atom or node is below some tolerance.

\begin{algorithm}[hbt]
 \SetAlgoSkip{bigskip}
 \LinesNumbered
 \SetKwInput{Input}{Input}
 \SetKwInput{Output}{Output}
 \SetKwBlock{Repeat}{repeat}{end}
 \setstretch{1.2}
 \setlength{\commentWidth}{0.6\textwidth} 
 \newcommand{\atcp}[1]{\tcp*[r]{\makebox[\commentWidth]{#1\hfill}}}
 \caption{Discrete boundary element method (\texttt{DBEM})}
 \label{algo:dbem}
 \Input{boundary condition $\bar{\displ}^\I$, right hand side $r$}
 $\displInf^\lb \,\leftarrow\, \G^{\lb|\l}[r]$
  \atcp{solve \textbf{(H)}$_\infty$ (eq.\,\eqref{eq:sinc.Hinf})}
 $w^\I \,\leftarrow\, (\Id^\II - \F^\II)[\bar{\displ} - \displInf]$
  \atcp{compute right hand side of the linear system}
 $\force^\I \,\leftarrow\, - \inv{\G^\II}[w]$
  \atcp{solve linear system \eqref{eq:sinc.Hom_ls} for the boundary forces $\force^\I$}
 $\displHom^\l \,\leftarrow\, \F^{\l|\I}[\bar{\displ} - \displInf] - \G^{\l|\I}[\force]$
  \atcp{compute solution of \textbf{(H)}$_{\rm hom}$ (eq.\,\eqref{eq:sinc.Hom_sol})}
 $\displ_\indHarm^\l \,\leftarrow\, \displInf^\l + \displHom^\l$
  \atcp{assemble the full harmonic solution}
 \Output{solution $\displ_\indHarm^\l$}
\end{algorithm}

\begin{algorithm}[hbt]
 \SetAlgoSkip{bigskip}
 \LinesNumbered
 \SetKwInput{Input}{Input}
 \SetKwInput{Output}{Output}
 \SetKwBlock{Repeat}{repeat}{end}
 \setstretch{1.2}
 \setlength{\commentWidth}{0.5\textwidth}
 \newcommand{\atcp}[1]{\tcp*[r]{\makebox[\commentWidth]{#1\hfill}}}
 \caption{Sinclair method for bounded problems (\texttt{Sinc})}
 \label{algo:sinc}
 \Input{natural boundary condition $\bar{\displ}^\I$, external force $\force_\mrm{ext}$}
 \uIf{\textbf{initial guess}}{
  $\displ^\lb_0 \,\leftarrow\, \texttt{DBEM}(\bar{\displ}^\I, \force_\mrm{ext})$;
 }\Else{
  $\displ^\lb_0 \,\leftarrow\, 0$;
 }
 $k \,\leftarrow\, 0$;\\
 \While{ $\| \var{}{}{\Etot^\a}(\displ_k) \|_{l^\infty} < TOL \; \wedge \; \| \var{}{}{\Etot^\c}(\displ_k) \|_{l^\infty} < TOL$ }{
  $\displ^\a_{k+1} \,\leftarrow\, \arg \, \left\{\, \underset{v^\a}{\min} \, \Etot^\a(\{v^\a,\displ^\p_k\}) \,\right\}$
   \atcp{solve \textnormal{\textbf{(AH)}}}
  $\force^\ip_{\mrm{inh},k+1} \,\leftarrow\, \Lh^{\ip|\a}[\displ_{k+1}] + \Lh^{\ip|\c}[\displ_k]$
   \atcp{compute inhomogeneous force}
  \If{\textnormal{\textbf{not}} \textbf{initial guess} \textnormal{\textbf{and}} $k=0$}{
   $\displ^\c_k \,\leftarrow\, \G^{\c|\c}[\force_\mrm{ext}]$;
  }
  $\displ^\c_{k+1} \,\leftarrow\, \displ^\c_k + \texttt{DBEM}(\bar{\displ}^\I - \displ^\I_k, -\force^\ip_{\mrm{inh},k+1})$
   \atcp{solve \textnormal{\textbf{(H)}}}
  $k \,\leftarrow\, k+1$;\\
 }
 \Output{global solution $\displ_k$}
\end{algorithm}


\section{Convergence analysis}
\label{sec:conv}

\newcommand{\displRef}{\ensuremath{\displ_\mrm{ref}}} 
\newcommand{\errmod}{\ensuremath{e}} 

\begin{modified}%
We now analyze the convergence behavior of the Sinclair method. In this work we focus on analyzing the convergence behavior to solutions around a linearized state. That is, we linearize the coupled problem around some \emph{general} displacement $\displ_\mrm{F}$ and seek for a solution $\displ \in \clV$ such that
\begin{equation}\label{eq:conv.lin_problem}
   \Lcpl(\displ_\mrm{F})[\displ - \displ_\mrm{F}]
 = \left\{
 \begin{aligned}
  \Lhnl[\displ - \displ_\mrm{F}] &= r_\mrm{F} &\quad& \text{in} \; \lat^\a, \\
  \Lh  [\displ - \displ_\mrm{F}] &= r_\mrm{F} &     & \text{in} \; \lat^\c,
 \end{aligned}
 \right.
 \qquad \text{with} \; r_\mrm{F} = \force_\mrm{ext} - \Lcpl(\displ_\mrm{F})[\displ_\mrm{F}],
\end{equation}
where $\Lhnl$ is the nonlocal harmonic operator which is obtained after linearizing the atomistic energy. The objective of this section is to study the convergence to $\displ$ using the iteration equation \eqref{eq:sinc.sinc}, starting from some initial guess $\displ_0$.

\textit{%
We remark that a linearization, though restrictive, allows us to study the approximate structure of the method in a more accessible setting in order to benchmark the convergence behavior, in particular, depending on how strongly $\Lh$ differs from $\Lhnl$.%
}\;%
Nevertheless, we point out that the techniques which will be developed for the linearized problem can be generalized in the sense that a nonlinear problem can be considered as a sequence of linear problems. Moreover, in Section \ref{sec:dyn_relax} we shall see that, using the tools from the linear analysis, we are able to improve the convergence behavior of general nonlinear problems by optimizing the transmission conditions between both problems around intermediate linearized states.

Our analysis is closely related to domain decomposition methods for partial differential equations (e.g. \citep{toselli_domain_2005}) and structured as follows. In Section \ref{sec:conv.proj} we derive the iteration operator which relates the error with respect to the fully atomistic reference solution between two consecutive iterations. In Section \ref{sec:conv.rate} and \ref{sec:conv.stab} we then investigate the convergence properties of the Sinclair method by analyzing this iteration operator. In Section \ref{sec:conv.relaxation} we further generalize the Sinclair method by relaxing the transmission conditions between the two subproblems in order to improve the convergence behavior.

\subsection{Iteration operator}
\label{sec:conv.proj}

In this section we will quantify the error $\displRef - \displ_{k+1}$ in the $k$+1-th iteration with respect to the reference solution $\displRef$, which solves $\Lhnl[\displRef] = \force_\mrm{ext} \; \text{in} \; \lat$, by deriving an equation of the form
\begin{equation}\label{eq:conv.iter_eq}
 \displRef - \displ_{k+1} = \projop[\displRef - \displ_k] + \errmod \qquad \text{in} \; \lat,
\end{equation}
where $\projop$ is the iteration operator and $\errmod$ is the modeling error. In other words, it will be shown that the error can be decomposed into a convergent part $\projop[\displRef - \displ_k]$ (provided that our method is stable), which describes the convergence to $\displ$, and a non-convergent part $\errmod$, independent of $\displ_k$, which occurs due to replacing $\Lhnl$ with $\Lh$ in $\latC$.

Before presenting the result, we will abbreviate the block operator used to update the harmonic problem \textbf{(H)}, i.e., the ``$\c|\c$'' block of $\inv{\Lh^{\l|\l}}$ in \eqref{eq:sinc.displ_upd}, as $\inv{\clS^\cc}$, where
\begin{equation}\label{eq:conv.schur}
 \clS^\cc = \Lh^\cc - \Lh^\ca\inv{(\Lh^\aa)}\Lh^\ac.
\end{equation}
Equation \eqref{eq:conv.schur} is nothing but the Schur complement of $\Lh^\aa$ in $\Lh^{\l|\l}$ (sometimes in the literature also denoted by $\Lh^{\l|\l}/\Lh^\aa$).

We will split the derivation of the iteration operator $\projop$ into two parts. First, in Lemma \ref{lem:itop1}, we will derive the general structure of $\projop$ using the algorithm described in Section \ref{sec:sinc.description}, which is independent of the implementation. Second, in Lemma \ref{lem:itop2} we will then recast $\projop$ from Lemma \ref{lem:itop1} by replacing $\clS^\cc$ with the discrete boundary element method from Section \ref{sec:sinc.H_k+1}.

Using the previous definitions, we now state the first main result of Section \ref{sec:conv}:

\begin{customlem}{1a}[\textbf{Iteration operator}]\label{lem:itop1}
 Let $\displRef$ be a unique solution to \eqref{eq:atom.ato_problem} with $\L = \Lhnl$. Then, using the iteration equation \eqref{eq:sinc.sinc}, the error in the $k$+1-th iteration can be written as
 \begin{equation}\label{eq:conv.proj}
  \displRef - \displ_{k+1} = \projop[\displRef - \displ_k] + \P_\rme[\displRef] \qquad \text{in} \; \lat,
 \end{equation}
 with the iteration operator $\P : \clV(\lat) \rightarrow \clV(\lat)$ and the operator $\P_\rme : \clV(\lat) \rightarrow \clV(\lat)$ given by
 \begin{align}\label{eq:conv.projop}
  \P
  =
  \begin{pmatrix} \nullOp^\aa & \inv{\Lhnl^\aa}\Lhnl^\ac \\ \nullOp^\ca & \Id^\cc - \inv{\clS^\cc} (\Lh^\cc - \Lh^\ca\inv{\Lhnl^\aa}\Lhnl^\ac) \end{pmatrix},
  &&
  \P_\rme
  =
  \begin{pmatrix} \nullOp^\aa & \nullOp^\ac \\ \inv{\clS^\cc}(\Lhnl^\ca - \Lh^\ca)  & \inv{\clS^\cc}(\Lhnl^\cc - \Lh^\cc) \end{pmatrix}.
 \end{align}
 
 Moreover, for the ``$|\c$'' block of the operator $\P$ we have
 \begin{equation}\label{eq:conv.T^|c}
  \P^{|\c}
  =
  \begin{pmatrix} \P^\ac \\ \P^\cc \end{pmatrix}
  =
  \begin{pmatrix} \P^\ap & \P^{\a|\c\setminus\p} \\ \P^\cp & \P^{\c|\c\setminus\p} \end{pmatrix}
  =
  \begin{pmatrix} \P^\ap & \nullOp^{\a|\c\setminus\p} \\ \P^\cp & \nullOp^{\c|\c\setminus\p} \end{pmatrix},
 \end{equation}
 where the index ``$\c\setminus\p$'' refers to functions defined on $\latC \setminus \latP$.
\end{customlem}

\begin{prf}
 To avoid further technicalities, we assume that $k>1$ such that the contributions due to the boundary conditions and the external force are already contained in $\displ_1$.
 
 We begin by writing the global solution $\displ_{k+1/2}$ in $\lat$, that is, after \textbf{(AH)}$_{k+1}$ has been solved, as
 \begin{align}
  \begin{split}
   \begin{pmatrix} \displ^\a_{k+1/2} \\ \displ^\c_{k+1/2} \end{pmatrix}
   =
   \begin{pmatrix} \displ^\a_{k+1} \\ \displ^\c_k \end{pmatrix}
   &=
   \begin{pmatrix} \displ^\a_k \\ \displ^\c_k \end{pmatrix}
   +
   \begin{pmatrix} \inv{\Lhnl^\aa} & \nullOp^\ac \\ \nullOp^\ca & \nullOp^\cc \end{pmatrix}
   \begin{pmatrix} \force_\mrm{ext}^\a - \Lhnl^{\a|}[\displ_k] \\ \force_\mrm{ext}^\c - \Lh^{\c|}[\displ_k] \end{pmatrix} \\
   &=
   \begin{pmatrix} \displ^\a_k \\ \displ^\c_k \end{pmatrix}
   +
   \begin{pmatrix} \inv{\Lhnl^\aa} & \nullOp^\ac \\ \nullOp^\ca & \nullOp^\cc \end{pmatrix}
   \begin{pmatrix} \Lhnl^{\a|}[\displRef] - \Lhnl^{\a|}[\displ_k] \\ \Lhnl^{\c|}[\displRef] - \Lh^{\c|}[\displ_k] \end{pmatrix} \\
   &=
   \begin{pmatrix} \displ^\a_k \\ \displ^\c_k \end{pmatrix}
   +
   \begin{pmatrix} \Id^\aa & \inv{\Lhnl^\aa}\Lhnl^\ac \\ \nullOp^\ca & \nullOp^\cc \end{pmatrix}
   \begin{bmatrix} \displRef^\a - \displ^\a_k \\ \displRef^\c - \displ^\c_k \end{bmatrix}.
  \end{split}
\end{align}
 Analogously, we can write the global solution $\displ_{k+1}$ after solving \textbf{(H)}$_{k+1}$ as
 \begin{align}
  \begin{split}
   \begin{pmatrix} \displ^\a_{k+1} \\ \displ^\c_{k+1} \end{pmatrix}
   =
   \begin{pmatrix} \displ^\a_{k+1/2} \\ \displ^\c_{k+1/2} \end{pmatrix}
   &+
   \begin{pmatrix} \nullOp^\aa & \nullOp^\ac \\ \nullOp^\ca & \inv{\clS^\cc} \end{pmatrix}
   \begin{pmatrix} \force_\mrm{ext}^\a - \Lhnl^{\a|}[\displ_{k+1/2}] \\ \force_\mrm{ext}^\c - \Lh^{\c|}[\displ_{k+1/2}] \end{pmatrix} \\
   =
   \begin{pmatrix} \displ^\a_{k+1/2} \\ \displ^\c_{k+1/2} \end{pmatrix}
   &+
   \begin{pmatrix} \nullOp^\aa & \nullOp^\ac \\ \nullOp^\ca & \inv{\clS^\cc} \end{pmatrix}
   \begin{pmatrix} \Lhnl^{\a|}[\displRef] - \Lhnl^{\a|}[\displ_{k+1/2}] \\ \Lhnl^{\c|}[\displRef] - \Lh^{\c|}[\displ_{k+1/2}] \end{pmatrix} \\
   =
   \begin{pmatrix} \displ^\a_{k+1/2} \\ \displ^\c_{k+1/2} \end{pmatrix}
   &+
   \begin{pmatrix} \nullOp^\aa & \nullOp^\ac \\ \inv{\clS^\cc}\Lh^\ca & \inv{\clS^\cc}\Lh^\cc \end{pmatrix}
   \begin{bmatrix} \displRef^\a - \displ^\a_{k+1/2} \\ \displRef^\c - \displ^\c_{k+1/2} \end{bmatrix}
   \\ &+
   \begin{pmatrix} \nullOp^\aa & \nullOp^\ac \\ \inv{\clS^\cc}(\Lhnl^\ca - \Lh^\ca)  & \inv{\clS^\cc}(\Lhnl^\cc - \Lh^\cc) \end{pmatrix}
   \begin{bmatrix} \displRef^\a \\ \displRef^\c \end{bmatrix}.
  \end{split}
 \end{align}
 Thus, we can write the multiplicative iterates as
 \begin{alignat}{4}
  \label{eq:conv.iter1}
  \displ_{k+1/2} &= \displ_k       &&+ \P_\a[\displRef -\displ_k]       &\qquad& \text{in} \; \lat, \\
  \label{eq:conv.iter2}
  \displ_{k+1}   &= \displ_{k+1/2} &&+ \P_\c[\displRef -\displ_{k+1/2}] + \P_\rme[\displRef] &\qquad& \text{in} \; \lat,
 \end{alignat}
 with the operators
 \begin{align}
  \P_\a &=
  \begin{pmatrix} \Id^\aa & \inv{\Lhnl^\aa}\Lhnl^\ac \\ \nullOp^\ca & \nullOp^\cc \end{pmatrix},
  &
  \P_\c =
  \begin{pmatrix} \nullOp^\aa & \nullOp^\ac \\ \inv{\clS^\cc}\Lh^\ca & \inv{\clS^\cc}\Lh^\cc \end{pmatrix},
 \end{align}
 and $\P_\rme$ as given in \eqref{eq:conv.projop}. Using \eqref{eq:conv.iter1} in \eqref{eq:conv.iter2} now permits to obtain
 \begin{equation}
  \displRef - \displ_{k+1} = \P[\displRef - \displ_k] + \P_\rme[\displRef] \qquad \text{in} \; \lat
 \end{equation}
 with
 \begin{equation}\label{eq:conv.projop_proof}
  \begin{split}
   \P = \Id - \P_\a - \P_\c + \P_\c\P_\a
      = \begin{pmatrix} \nullOp^\aa & \inv{\Lhnl^\aa}\Lhnl^\ac \\ \nullOp^\ca & \Id^\cc - \inv{\clS^\cc} (\Lh^\cc - \Lh^\ca\inv{\Lhnl^\aa}\Lhnl^\ac) \end{pmatrix}.
  \end{split}
 \end{equation}
 
 To show that $\P^{|\c\setminus\p} = \nullOp^{|\c\setminus\p}$ we use that fact that $\forall\,v\in\clV$ we have $\Lhnl^{\a|\c\setminus\p}[v] = 0$ since atoms do interact with nodes outside $\latP$ (cf. Figure \ref{fig:domain_decomp}). Therefore, it follows trivially that $\P^{\a|\c\setminus\p} = \nullOp^{\a|\c\setminus\p}$. To show that $\P^{\c|\c\setminus\p} = \nullOp^{\c|\c\setminus\p}$ we rewrite $\P^\cc$ as
 \begin{align}
  \P^\cc
  =
  \Id^\cc - \inv{\clS^\cc} (\Lh^\cc - \Lh^\ca\inv{\Lhnl^\aa}\Lhnl^\ac)
  &\overset{\phantom{\eqref{eq:conv.schur}}}{=}
  \inv{\clS^\cc}(\clS^\cc - \Lh^\cc + \Lh^\ca\inv{\Lhnl^\aa}\Lhnl^\ac) \\
  &\overset{\eqref{eq:conv.schur}}{=}
  \inv{\clS^\cc} \Lh^{\c|\a} (\inv{\Lhnl^\aa}\Lhnl^\ac - \inv{\Lh^\aa}\Lh^\ac)
 \end{align}
 and use the fact that $\forall\,v\in\clV$ we also have $\Lh^{\a|\c\setminus\p}[v] = 0$.
 \qed
\end{prf}

In the form of Lemma \ref{lem:itop1} $\P$ is not convenient to implement. In particular, for the $\P^\cp$ block, it would be necessary to invert the Schur complement $\clS^\cc$, that is, inverting a $dN^\c\times dN^\c$ block matrix, where $N^\c = \#\latC$, which can be very costly if $N^\c$ is large. Following Algorithm \ref{algo:sinc}, we will now derive a more efficient representation of $\P^\cp$ using the discrete boundary element method which only involves boundary operators. We point out that this has no influence on the convergence behavior of the method and, therefore, the reader may already jump to the following section without missing essential information; should we require any result from Lemma \ref{lem:itop2} we will explicitly refer to it.

\begin{customlem}{1b}[\textbf{Implementation of the block operator $\P^\cp$}]\label{lem:itop2}
 The block operator $\P^\cp$ from \eqref{eq:conv.projop} can be written as $\P^\cp = \P^\cp_1 + \P^\cp_2$, with
 \begin{align}\label{eq:conv.Tcp_splitting}
  \P^\cp_1 = \big(\inv{\clS^\cc}\big)^{\c|\i}\Lh^\ipd, &&
  \P^\cp_2 = \big(\inv{\clS^\cc}\big)^{\c|\ip}\Lh^\ipi\big(\inv{\Lhnl^\aa}\big)^\ipr\Lhnl^\prp,
 \end{align}
 where the index ``\pr'' refers to functions defined on the set of atoms which interact with the continuum nodes in $\latP$ (cf. Figure \ref{fig:domain_decomp}), and
 \begin{equation}\label{eq:conv.Schur_inv}
  \inv{\clS^\cc} = \G^\cc - \B^\cI\G^\Ic, \qquad \text{with} \quad \B^\cI = \F^\cI + \G^\cI\inv{\G^\II}(\Id^\II - \F^\II).
 \end{equation}
\end{customlem}

We remark here that the expression for $\inv{\clS^\cc}$ can be obtained from Algorithm \ref{algo:dbem}. Therein, we denote $\B^{\c|\I}$ as the boundary operator which maps some boundary displacement $v^\I$ to a homogeneous solution in $\latC$ (this operation corresponds to lines 2--4 of Algorithm \ref{algo:dbem}; or, alternatively, to $\texttt{DBEM}(v^\I,0)$).

\begin{prf}
 The strategy of the proof is to write $\P^{\c|\p}$ from \eqref{eq:conv.projop} as $\P^{\c|\p} = \P_1^{\c|\p} + \P_2^{\c|\p}$, with (using \eqref{eq:conv.Schur_inv})
 \begin{align}
  \P^\cp_1 = \clI^\cp - (\G^\cc - \B^\cI\G^\Ic)\Lh^\cp, &&
  \P^\cp_2 = (\G^\cc - \B^\cI\G^\Ic)\Lh^\ca\inv{\Lhnl^\aa}\Lhnl^\ap,
 \end{align}
 and show that this splitting equals \eqref{eq:conv.Tcp_splitting}.
 
 \subsubsection*{Block operator \texorpdfstring{$\P^\cp_1$}{}:}
 
 We first analyze the block $\G^\cc\Lh^\cp$. To that end, we recall from \eqref{eq:con.GLh=I_2} that we can likewise write it using
 \begin{equation}
  \begin{array}{@{}r@{{}\mathrel{}}c@{\mathrel{}{}}l@{}}
       (\G\Lh)^\cp = \Id^\cp
   &=& \G^{\c|\lb}\Lh^{\lb|\p} + \underbrace{\G^{\c|\re}\Lh^{\re|\p}}_{\overset{\eqref{eq:con.Lrcb_props}}{=} \nullOp^\cp} \\
   &=& \G^\ca\Lh^\ap + \G^\cc\Lh^\cp + \G^\cI\Lh^{\I|\p} \\
   &=& \G^\ci\Lh^{\i|\p} + \G^\cc\Lh^\cp + \G^\cI\Lh^{\I|\p},
  \end{array}
 \end{equation}
 where $\G^\ca\Lh^\ap = \G^\ci\Lh^{\i|\p}$ is due to local interactions, i.e., for all $v\in\clV$, $\Lh^{\a|\p}[v] = 0 \; \text{in} \; \latA \setminus \latI$. Thus, using $\G^\cc\Lh^\cp = \Id^\cp - \G^\ci\Lh^{\i|\p} - \G^\cI\Lh^{\I|\p}$, we can write $\P^\cp_1$ as
 \begin{equation}
  \P_1^\cp = \G^\ci\Lh^{\i|\p} + \G^\cI\Lh^{\I|\p} + \B^\cI\G^\Ic\Lh^\cp.
 \end{equation}
 
 Next, we analyze the latter term $\B^\cI\G^\Ic\Lh^\cp$. To obtain an alternative expression for the block $\G^\Ic\Lh^\cp$ we again use \eqref{eq:con.GLh=I_2} and write
 \begin{equation}
  (\G\Lh)^{\I|\p} = \G^\Ii\Lh^{\i|\p} + \G^\Ic\Lh^\cp + \G^\II\Lh^{\I|\p} = \nullOp^{\I|\p}.
 \end{equation}
 With $\G^\Ic\Lh^\cp = -\G^\Ii\Lh^{\i|\p} - \G^\II\Lh^{\I|\p}$, we have
 \begin{equation}\label{eq:conv.term_23}
  \begin{aligned}
   \P_1^{\c|\p}
   &=
   \G^\ci\Lh^{\i|\p} + \G^\cI\Lh^{\I|\p} - \B^\cI\G^\II\Lh^{\I|\p} - \B^\cI\G^\Ii\Lh^{\i|\p} \\
   &=
   \G^\ci\Lh^{\i|\p} -\B^\cI\G^\Ii\Lh^{\i|\p},
  \end{aligned}
 \end{equation}
 where we have used the relation $(\G^\cI - \B^\cI\G^\II)\Lh^{\I|\p} = \nullOp^\cp$ which follows from the fact that, for any input $v^\p$, both operators $\G^\cI$ and $\B^\cI\G^\II$ produce the same solution in $\lat^\c$ (since the boundary condition \emph{is} the Green function).
 
 \subsubsection*{Block operator \texorpdfstring{$\P^\cp_2$}{}:}
 
 The expression for $\P^\cp_2$ from \eqref{eq:conv.Tcp_splitting} can be obtained by recalling the properties of the operators $\Lhnl$ and $\Lh$: for all $v\in\clV$ we have $\Lhnl^\ap[v] = 0 \; \text{in} \; \latA \setminus \lat^\pr$, $\Lh^{\c|\a\setminus\i}[v] = 0 \; \text{in} \; \latC$ (here, the index ``$\a\setminus\i$'' refers to functions defined on $\latA \setminus \latI$), and $\Lh^\ca[v] = 0 \; \text{in} \; \latC \setminus \latIpl$.
 \qed
\end{prf}
\end{modified}%

Now, using Lemma \ref{lem:itop2}, we can split the iteration operator into a contribution $\Pinf$ due to the infinite inhomogeneous problem \textbf{(H)}$_\infty$ and a contribution $\Pbnd$ due to the finite homogeneous problem \textbf{(H)}$_\mrm{hom}$ such that
\begin{equation}\label{eq:conv.Tinf+Tbnd}
 \P = \Pinf + \Pbnd =
 \begin{pmatrix}
  \nullOp^\aa & \P^\ap    & \nullOp^{\a|\c\setminus\p} \\
  \nullOp^\ca & \Pinf^\cp & \nullOp^{\c|\c\setminus\p}
 \end{pmatrix}
 +
 \begin{pmatrix}
  \nullOp^\aa & \nullOp^\ap & \nullOp^{\a|\c\setminus\p} \\
  \nullOp^\ca & \Pbnd^\cp   & \nullOp^{\c|\c\setminus\p}
 \end{pmatrix},
\end{equation}
where
\begin{align}
 \Pinf^\cp
 &=
 \G^\ci\Lh^\ipd + \G^\cip\Lh^\ipi\big(\inv{\Lhnl^\aa}\big)^\ipr\Lhnl^\prp,
 \label{eq:conv.Tinf} \\
 \Pbnd^\cp
 &=
 - \B^\cI\G^\Ii\Lh^\ipd - \B^\cI\G^\Iip\Lh^\ipi\big(\inv{\Lhnl^\aa}\big)^\ipr\Lhnl^\prp.
 \label{eq:conv.Tbnd}
\end{align}
Note that the latter vanishes for infinite problems.

\subsection{Convergence rate and error analysis}
\label{sec:conv.rate}

With Lemma \ref{lem:itop1} we can now prove the convergence rate of the Sinclair method:

\begin{modified}%
\begin{thm}[\textbf{Convergence rate}]\label{thm:conv_rate}
 Under the assumption that the iteration operator $\P$ (eq. \ref{eq:conv.projop}) admits the eigendecomposition $\P = \clQ \clD \clQ^{-1}$, the norm of the error \eqref{eq:conv.proj} in the $k{+}1$-th iteration can be bounded from above as
 \begin{equation}\label{eq:conv.error_bound}
  \| \displRef - \displ_{k+1} \| \le \sigma^{k+1} \| \clQ \| \| \inv{\clQ} \|
  \left( \| \displRef - \displ_0 \|  - \frac{1}{1 - \sigma} \| \P_\rme[\displRef] \| \right)
  +
  \frac{1}{1 - \sigma} \| \clQ \| \| \inv{\clQ} \| \| \P_\rme[\displRef] \|,
 \end{equation}
 where $\sigma = \sigma(\P)$ is the spectral radius of $\P$.
 
 Moreover, the spectral radius of \P is equivalent to the spectral radius of the block $\P^\pp$, that is, $\sigma = \sigma(\P^\pp)$.
\end{thm}

\begin{prf}
 Using \eqref{eq:conv.proj}, we can write the iterates until the $k$+1-th iteration as
 \begin{align}
  \displRef - \displ_1 &= \P[\displRef - \displ_0] + \P_\rme[\displRef] \qquad \text{in} \; \lat, \\
  \vdots \nonumber \\
  \label{eq:conv.proj2}
  \displRef - \displ_{k+1} &= \P[\displRef - \displ_k] + \P_\rme[\displRef] \qquad \text{in} \; \lat.
 \end{align}
 Recursively using the error from the previous iteration(s) in \eqref{eq:conv.proj2}, we obtain
 \begin{equation}\label{eq:conv.u_k+1-u=P[u_0-u]}
  \displRef - \displ_{k+1} = \P^{k+1}[\displRef - \displ_0] + \sum _{i=0}^k\P^i \P_\rme[\displRef], \qquad \text{where} \quad \P^n = \prod_{i=1}^{n\ge 1} \P = \clQ \clD \clQ^{-1} \clQ \clD \clQ^{-1} \, \cdots = \clQ \clD^n \clQ^{-1}.
 \end{equation}
 Taking norms on both sides and applying the Cauchy-Schwarz inequality we get
 \begin{equation}
  \begin{aligned}
   \| \displRef - \displ_{k+1} \|
   &\le
   \sigma^{k+1} \| \clQ \| \| \inv{\clQ} \| \| \displRef - \displ_0 \|
   +
   \sum _{i=0}^k\sigma^i \| \clQ \| \| \inv{\clQ} \| \| \P_\rme[\displRef] \| \\
   &=
   \sigma^{k+1} \| \clQ \| \| \inv{\clQ} \| \| \displRef - \displ_0 \|
   +
   \frac{1 - \sigma^{k+1}}{1 - \sigma} \| \clQ \| \| \inv{\clQ} \| \| \P_\rme[\displRef] \|,
  \end{aligned}
 \end{equation}
 which is the expected bound.
 
 The second statement is obtained by noting that only the block $\P^{\c|\p}$ in $\P$ is nonzero (cf. Lemma \ref{lem:itop1}). Shifting indices such that the ``$\p$'' column appears on the right 
 gives an upper triangular block matrix and the spectral radius of such a matrix is equivalent to the spectral radius of the nonzero diagonal block, which is nothing but $\P^\pp$.
 \qed
\end{prf}

From \eqref{eq:conv.error_bound}, it can be seen that the error is bounded by a term depending on $k$+1 powers of $\sigma$ and a constant term, independent of $k$.
\end{modified}%

Hence, the first term defines the convergence rate of the method and also its stability, that is, provided that $\sigma < 1$, the method will converge (we turn to this question in the following section). For the special case when the atomistic and continuum models coincide, the exact number of required iterations can therefore be immediately deduced from Theorem \ref{thm:conv_rate}:

\begin{cor}\label{cor:conv_rate}
 Let $\Lhnl^\a = \Lh^a$ such that $\Lcpl = \Lh$. Then, Algorithm \ref{algo:sinc} converges in two steps.
\end{cor}

\begin{prf}
 It suffices to show that the iteration operator vanishes. Indeed, since $\inv{\clS^\cc}$ is now nothing but the inverse of the Schur complement $\Lcpl/\Lh^\aa = \Lh^\cc - \Lh^\ca\inv{(\Lh^\aa)}\Lh^\ac$ (cf. \citep{zhang_schur_2005}), it follows immediately that
 \begin{equation}
  \P^\cc = \Id^\cc - \inv{\clS^\cc} (\Lcpl/\Lh^\aa) = 0
 \end{equation}
 and, therefore, $\displ_2 - \displ = 0$.\qed
\end{prf}

This result is general and holds for arbitrary interaction stencils---provided that we can compute the corresponding lattice Green function.

\begin{new}%
The second term in \eqref{eq:conv.error_bound} is nonconvergent and bounds the modeling error of the coupled problem. From the structure of $\P_\rme$ (eq. \eqref{eq:conv.projop}) we see that $\| \P_\rme[\displRef] \|$ is nothing but the mismatch between both models in $\lat^\c$ after applying both $\Lhnl$ and $\Lh$ to the reference solution.

If we linearize \eqref{eq:conv.lin_problem} around a \emph{homogeneous displacement} $\displ_\rmF$, we can bound $\| \P_\rme[\displRef] \|$ by the Cauchy-Born modeling error. That is, assuming that $\displ_\rme$ solves $\Lh[\displ_\rme] = 0 \; \text{in} \; \latC$, given $\displRef$ as a boundary condition, we may write
\begin{equation}
 \begin{aligned}
  \| \P_\rme[\displRef] \|
  &=
  \| (\inv{\clS^\cc}) [ (\Lhnl^{\c|} - \Lh^{\c|})[\displRef] + \Lcc[\displ_\rme] - \Lcc[\displ_\rme] ) \| \\
  &=
  \| (\inv{\clS^\cc}\Lcc)[\displRef - \displ_\rme] \| \\
  &\le
  \| \inv{\clS^\cc}\Lcc \| \| \displRef^\c - \displ^\c_\rme \| \lesssim f(\grad{}{2}{\displRef^\c}, ...),
 \end{aligned}
\end{equation}
where $f$ is some function depending on higher gradients of the reference solution (cf., e.g., \citep{e_cauchyborn_2007} for some in-depth analysis of the Cauchy-Born approximation). This is what we expect when replacing the nonlocal atomistic model with the (Cauchy-Born) elasticity model---if higher gradients remain small in $\latC$, the coupled problem is supposed to be sufficiently accurate.
\end{new}%

\subsection{Stability}
\label{sec:conv.stab}

\subsubsection{General remarks on the stability}
\label{sec:conv.stab.general}

We now turn to the question whether the Sinclair method is stable, that is, under which conditions $\sigma < 1$ holds for the general case when $\Lh \neq \Lhnl$.

\begin{modified}%
A common strategy to prove such a stability result is to show that $\sigma < \| \P \| < 1$ (see, e.g., \citep{benzi_algebraic_2001}). This requires to show that $\P > 0$, meaning that $\P$ is positive (i.e., all elements of its associate matrix are $>0$), and that $\forall\,v>0$, $\P[v] < v$. The author is aware of one such related stability result by \citet{parks_connecting_2008} for the alternating Schwarz method. However, the proof in \citep{parks_connecting_2008} assumes that the atomistic Hessian is an \textit{M}-matrix.%
\footnote{A nonsingular matrix is an \textit{M}-matrix if its off-diagonal components are nonpositive and its inverse is nonnegative (cf. \citep{quarteroni_numerical_2007})}
This is generally not the case since its off-diagonal components are positive \emph{and} negative for physically admissible interatomic potentials. Unfortunately, even if the atomistic Hessian is an \textit{M}-matrix, the operator $\P$ is not strictly positive.
To see this, consider $\P$ from \eqref{eq:conv.projop_proof} factorized to $\P = (\Id - \P_\c)(\Id - \P_\a)$. Assume now that the matrix associated with $\Lh$ is an \textit{M}-matrix. Then we have that $\inv{\clS^\cc} > 0$ and, therefore, $(\Id - \P_\c) \ngtr 0$. Since $(\Id - \P_\a) > 0$, as shown in \citep{parks_connecting_2008}, it follows that $\P \ngtr 0$ in general. We have also observed this in numerical experiments.

Fortunately, in the following section, we will show that in one dimension the situation simplifies because the matrix associated with $\P$ has sufficiently low dimension which allows to place bounds on each of its eigenvalues directly. The 1d proof requires the stability of $\Lhnl$ and certain upper bounds on its inverse matrix elements, and it is likely that similar results would have to be established for higher dimensions.
But this is out-of-scope in the context of the current work.
\end{modified}%

\begin{new}%
Yet to give an intuitive understanding of the behavior of the spectral radius in order to provide some practical guidance for designing a convergent algorithm, we assess $\sigma = \sigma(\Delta) = \sigma(\P(\Delta))$ in terms of the perturbation $\Delta = \Lhnl - \Lh$. From Corollary \ref{cor:conv_rate} it follows
\begin{equation}\label{eq:conv.cont}
 \vert \sigma(\nullOp) - \sigma(\Delta) \vert = \sigma(\Delta) 
 =
 \sigma(\P(\nullOp) - \P(\Delta)) \le \| \P(\nullOp) - \P(\Delta) \|.
\end{equation}
Equation \eqref{eq:conv.cont} can be understood as a continuity result. Roughly speaking, for sufficiently small perturbations $\Delta$, there will be no sudden jump of $\sigma$ beyond 1. 

Moreover, we may write $\sigma(\P(\nullOp) - \P(\Delta))$ as
\begin{equation}
 \begin{aligned}
  \sigma(\P(\nullOp) - \P(\Delta))
  \overset{\text{Thm.\;\ref{thm:conv_rate}}}{=}
  \sigma(\P^{\p|\p}(\nullOp) - \P^{\p|\p}(\Delta))
  &\overset{\phantom{\text{Lem.\;\ref{lem:itop2}}}}{=}
  \sigma\Big(
  \hphantom{\;-}\Id^\pp - \big(\inv{\clS^\cc}\big)^{\p|\c} (\Lh^{\c|\p} - \Lh^{\p|\a}\inv{\Lhnl^\aa}\Lhnl^\ap) \\
  &\overset{\phantom{\text{Lem.\:\ref{lem:itop2}}}}{\phantom{=}}
  \phantom{\sigma\Big(}
  \;-\Id^\pp - \big(\inv{\clS^\cc}\big)^{\p|\c} (\Lh^{\c|\p} - \Lh^{\p|\a}\inv{\Lh^\aa}\Lh^{\a|\p})
  \Big) \\
  &\overset{\text{Lem.\;\ref{lem:itop2}}}{=}
  \sigma\Big(
  \big(\inv{\clS^\cc}\big)^{\p|\ip} \Lh^{\ip|\i} \Big(\big(\inv{\Lhnl^\aa}\big)^\ipr\Lhnl^\prp - \big(\inv{\Lh^\aa}\big)^\ipr\Lh^\prp\Big)
  \Big).
 \end{aligned}
\end{equation}
In particular, the latter term implies that the spectral radius changes with the difference between the solutions to a fully atomistic and a fully continuum problem---computed in $\latA$ and subject to identical boundary conditions in $\lat^\p$---\emph{on the interface} $\lat^\i$ (since $\Lh^{\ip|\i}$ only acts on this part of the solution). 
This implies that the elastic solution in $\lat^\a$ may be very inaccurate in the defect core, but if the defect is sufficiently far from the boundary linear elasticity is usually a good approximation such that the difference of the solutions in $\latI$ remains small. For such cases we expect that $\sigma < 1$.
\end{new}%

\begin{modified}%
\subsubsection{Stability result in one dimension}
\label{sec:conv.stab.1d}

In this section, we will show that the Sinclair method is unconditionally stable in one dimension requiring only the assumption that the atomistic problem is stable.

We will analyze the problem defined on the domain $\lat \subset \integ$, shown in Figure \ref{fig:domain_decomp_1d}, where
\begin{align}
 \latA := \left\{\, -M, -M+1, ..., 0, ..., M-1, M \,\right\},
 &&
 \latC := ( \left\{\, -N, ..., -M-1 \,\right\} \cup \left\{\, M+1, ..., N \,\right\} ).
\end{align}

\begin{figure}[hbt]
 \centering
 \includegraphics[width=0.67\textwidth]{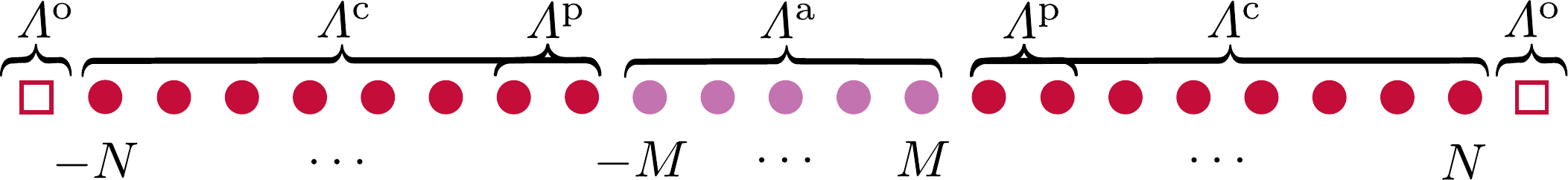}
 \caption{Domain decomposition for the one-dimensional problem considered in Section \ref{sec:conv.stab.1d}}
 \label{fig:domain_decomp_1d}
\end{figure}

We assume second-nearest neighbor interactions for the atomistic model and that the linearization takes place around a homogeneously deformed state. The interaction stencils for the atomistic and continuum problems then read (using linear lattice interpolants $\phi_{\ato}$, cf. \citep{curtin_atomistic/continuum_2003})
\begin{align}\label{eq:conv.stencils_1d}
 K_\mrm{hnl}(\ato - \atoB)
 =
 \left\{
 \begin{aligned}
   \;-k_2 &\qquad&& \text{if} \; \lvert \ato - \atoB \rvert = 2, \\
     -k_1 &\qquad&& \text{if} \; \lvert \ato - \atoB \rvert = 1, \\
     2k   &      && \text{if} \; \lvert \ato - \atoB \rvert = 0, \\
     0    &      && \text{else},
 \end{aligned}
 \right.
 &&
 K_\mrm{h}(\ato - \atoB)
 =
 \left\{
 \begin{aligned}
   \;-\bar{k} &\qquad&& \text{if} \; \lvert \ato - \atoB \rvert = 1, \\
     2\bar{k} &      && \text{if} \; \lvert \ato - \atoB \rvert = 0, \\
     0        &      && \text{else}.
 \end{aligned}
 \right.
\end{align}
where $k = k_1 + k_2$ and $\bar{k} = k_1+4k_2$. The corresponding Euler-Lagrange equation of the coupled problem is given by
\begin{equation}\label{eq:conv.euler_lagrange}
 \Lcpl[\displ](\ato) =
 \left\{
 \begin{aligned}\;
      -k_2\big(\displ(\ato - 2) - 2\displ(\ato) + \displ(\ato + 2)\big)
    -  k_1\big(\displ(\ato - 1) - 2\displ(\ato) + \displ(\ato + 1)\big)
  &= \force_\mrm{ext}(\ato) \qquad \forall\,\ato \in \lat^\a, \\
    - \bar{k}\big(\displ(\ato - 1) - 2\displ(\ato) + \displ(\ato + 1)\big)
  &= \force_\mrm{ext}(\ato) \qquad \forall\,\ato \in \lat^\c.
 \end{aligned}
 \right.
\end{equation}

For this particular system, \citet{dobson_stability_2010} have shown that $k_1>0$ and $k_2<0$ for Lennard-Jones-type interactions and that $k_1+4k_2>0$ is a necessary and sufficient condition to render $\Lhnl^\aa$ positive definite (i.e., stable) and, therefore, \eqref{eq:conv.euler_lagrange} a well-posed problem. Under this assumption, the Sinclair iteration equation \eqref{eq:sinc.sinc} converges: 

\begin{thm}[\textbf{Stability in 1d}]\label{thm:rho_sinc}
 Let $\Lcpl$ be given by \eqref{eq:conv.euler_lagrange}, with $k_1>0$, $k_2<0$ and $k_1+4k_2>0$. Then, $\sigma < 1$.
\end{thm}

\begin{prf}
 The proof is given in Appendix \ref{sec:apdx.proof_thm_rho_sinc}.
\end{prf}
\end{modified}%

\begin{figure}[hbt]
 \centering
 \begin{minipage}{0.5\textwidth}
  \centering
  (a)\\[0.5em]
  \includegraphics[width=0.9\textwidth]{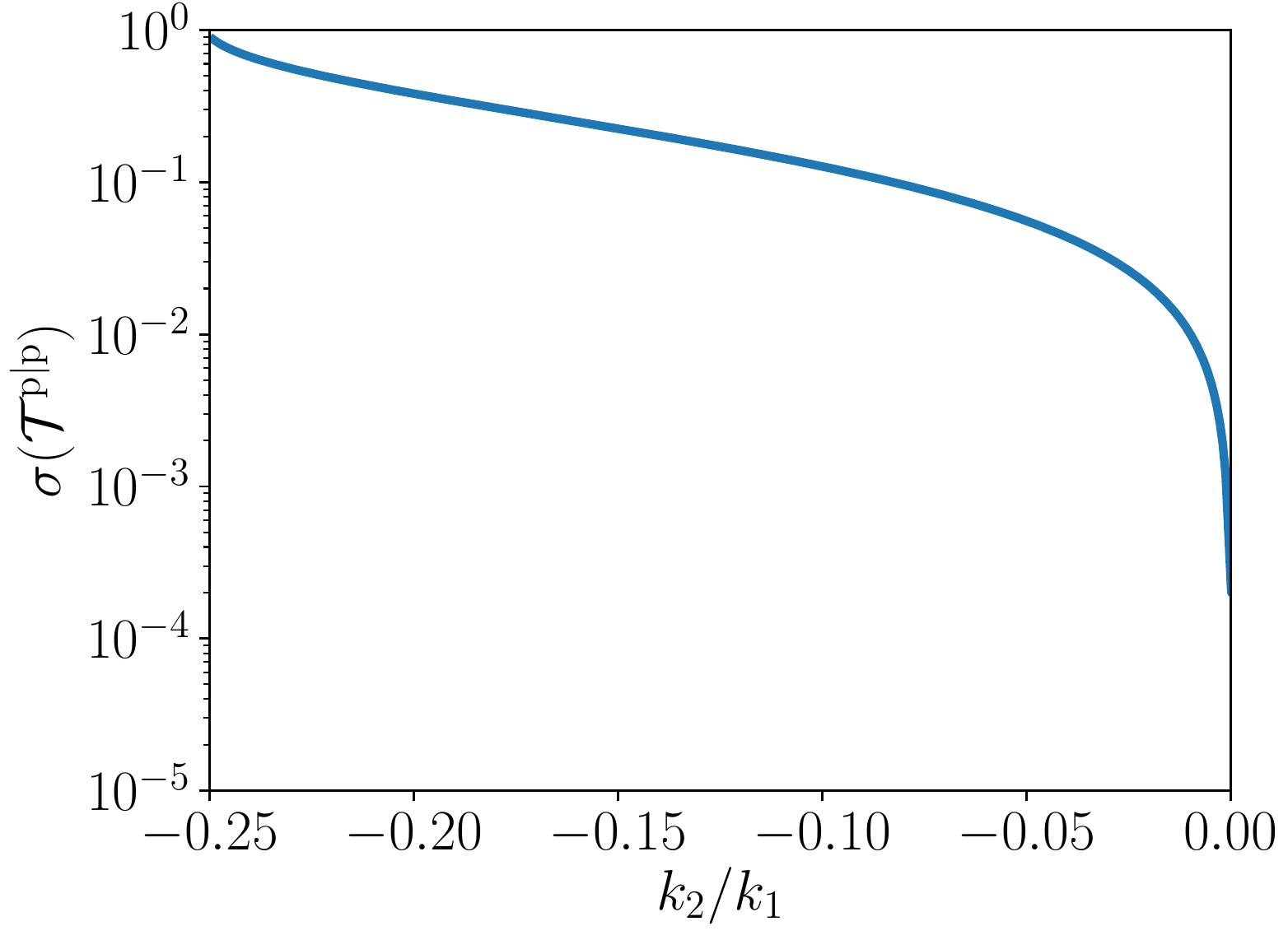}
 \end{minipage}\hfill
 \begin{minipage}{0.5\textwidth}
  \centering
  (b)\\[0.5em]
  \includegraphics[width=0.9\textwidth]{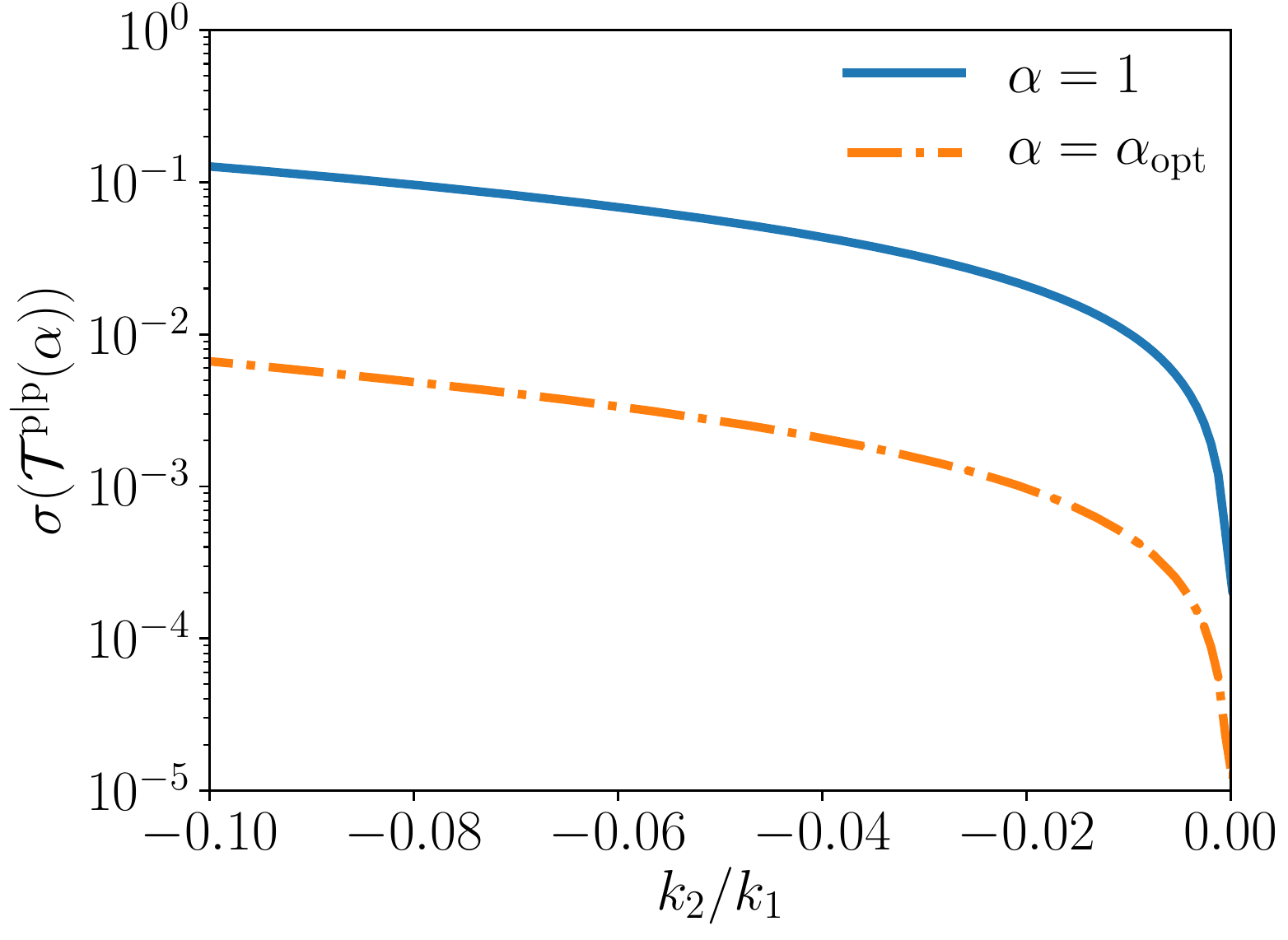}
 \end{minipage}
 \caption{(a) Spectral radius $\sigma$ of the iteration operator $\P$ corresponding to the one-dimensional problem \eqref{eq:conv.euler_lagrange} as a function of the nonlocality ratio $k_2/k_1$. (b) Spectral radius of $\P$ for the same problem with and without optimal relaxation (cf. Section \ref{sec:conv.relax_static})}
 \label{fig:sigma_vs_k2}
\end{figure}

In Figure \ref{fig:sigma_vs_k2} (a) the spectral radius $\sigma$ is shown as a function of $k_2/k_1 \in (-0.25,0.25]$ for an atomistic domain of size $M=10$. We observe the expected behavior: if $k_1 + 4k_2$ is close to $0$, $\sigma$ is close to but still strictly smaller than $1$, and $\sigma \rightarrow 0$ as $k_2 \rightarrow 0$.

\subsection{Relaxation}
\label{sec:conv.relaxation}

\subsubsection{Static relaxation}
\label{sec:conv.relax_static}

Even though the convergence rate of the Sinclair method is far superior than the convergence rate of the alternating Schwarz method, it is not optimal. This is in particular crucial whenever $\Lhnl$ and \Lh differ considerably. Although the solution in a region of interest in the atomistic domain may still be considered as good enough in such cases, intermediate solutions in the vicinity of the artificial interface will be non-smooth, slowing down the convergence.

One approach to accelerate the speed of convergence of domain decomposition solvers is relaxation (see, e.g., \citep{toselli_domain_2005}). The underlying idea is to control the transmission conditions between both problems in an optimal way by augmenting the artificial boundary conditions with a relaxation parameter. Here, we will relax the magnitude of the inhomogeneous force $\force_\mrm{inh}$. Therefore, we define the relaxation parameter $\alpha > 0$ and let
\begin{equation}
 \force_{\mrm{inh},k+1} = \force_{\mrm{inh},k+1}(\alpha) = \alpha\force_{\mrm{inh},k+1}
\end{equation}
The solution in $\lat^\c$ in the $k$+1-th iteration is now given by (cf. \eqref{eq:sinc.displ_upd2})
\begin{equation}
 \displ^\c_{k+1} = \displ^\c_{k+1}(\alpha) = \displ^\c_k - \inv{\clS^\cc}[\force_{\mrm{inh},k+1}(\alpha)] = \alpha\displ^\c_{k+1}(\alpha=0) + (1 - \alpha)\displ^\c_k.
\end{equation}
The anharmonic and the harmonic problem in each iteration are then defined as follows
\begin{empheq}[box=\fbox]{align}\label{eq:conv.sinc_relax}
 \textbf{(AH)}_{k+1} \; \left\{
 \begin{aligned}
  \; \L[\{\displ^\a_{k+1}, \displ^\p_{k+1}\}] &= \force_\mrm{ext}    &\;\,& \text{in} \; \lat^\a, \\
  \; \displ_{k+1}                             &= \displ^\c_k         &    & \text{in} \; \lat^\p,
 \end{aligned}
 \right.
 &&
 \textbf{(H)}_{k+1} \; \left\{
 \begin{aligned}
  \;     \Lh[\displ_{\indHarm,k+1}]
     &= \L^{|\a}[\displ_{\indHarm,0}] + \Lh^{|\cb}[\displ_0]              &\;\,& \text{in} \; \latA, \\
  \;     \Lh[\displ_{\indHarm,k+1}]
     &= -\Lh^{|\a}(\alpha)[\displ_{\indAHarm,k+1}] + \force_\mrm{ext} &\;\,& \text{in} \; \lat^\c, \\
  \; \displ &= \bar{\displ}                                             &    & \text{on} \; \lat^\I,
 \end{aligned}
 \right.
\end{empheq}
where
\begin{equation}
 \Lh^{|\a}(\alpha)[\displ_{\indAHarm,k+1}] = \alpha \sum_{i=1}^{k+1} \force_{\mrm{inh},i} \qquad \text{in} \; \latC.
\end{equation}

Similar to Section \ref{sec:conv.proj}, we now recast the iteration equation \eqref{eq:conv.sinc_relax} into the form $\displ - \displ_{k+1} = \projop(\alpha)[\displ - \displ_k]$, where we consider for brevity only the convergence to the solution $\displ$ of the coupled problem (eq. \eqref{eq:conv.lin_problem}) since we are only interested in modifying the convergence properties of \eqref{eq:conv.sinc_relax}. Here, $\projop(\alpha)$ is the iteration operator which now depends on the relaxation parameter $\alpha$. The following Lemma can therefore be considered as a generalization of Lemma \ref{lem:itop1} and \ref{lem:itop2} since $\P(\alpha)$ reduces to $\P$ as $\alpha \rightarrow 1$.

\begin{customlem}{2}[\textbf{Iteration operator for the relaxed Sinclair method}]\label{lem:projop_relax}
 Let $\displ$ be a unique solution to \eqref{eq:conv.lin_problem}. Then, using the iteration equation \eqref{eq:conv.sinc_relax}, the error in the $k$+1-th iteration can be written as
 \begin{equation}\label{eq:conv.proj_relax}
  \displ - \displ_{k+1} = \P(\alpha)[\displ - \displ_k] \qquad \text{in} \; \lat,
 \end{equation}
 with the iteration operator $\P(\alpha) : \clV(\lat) \rightarrow \clV(\lat)$ given by
 \begin{equation}\label{eq:conv.projop_relax}
  \P(\alpha)
  = \begin{pmatrix} \P^\aa & \P^\ac \\ \P^\ca & \P^\cc(\alpha) \end{pmatrix}
  = \begin{pmatrix} \nullOp^\aa & \inv{\Lhnl^\aa}\Lhnl^\ac \\ \nullOp^\ca & \Id^\cc - \inv{\clS^\cc} (\Lh^\cc(\alpha) - \Lh^\ca(\alpha)\inv{\Lhnl^\aa}\Lhnl^\ac) \end{pmatrix}.
 \end{equation}
 
 Moreover, $\P^\ac$ is given by \eqref{eq:conv.T^|c} and
 \begin{align}\label{eq:conv.projop_cc_relax}
  \P^\cc(\alpha) &=
  \begin{pmatrix} \P^\cp(\alpha) & \nullOp^{\c|\c\setminus\p} \end{pmatrix},
  \qquad \text{with} \quad
  \P^\cp(\alpha) = \P^\cp_1 + \alpha\P^\cp_2 + (1-\alpha)\P^\cp_3,
 \end{align}
 where $\P^\cp_1$ and $\P^\cp_2$ are given by \eqref{eq:conv.Tcp_splitting} and
 \begin{equation}
  \P^\cp_3 = (\G^\cip - \B^\cI\G^\Iip)\Lh^{\ip/\p}.
 \end{equation}
\end{customlem}

\begin{prf}[Sketch of the proof]
 We do not give a full proof for compactness as it would largely resemble the proof of Lemma \ref{lem:itop1} and \ref{lem:itop2} with the difference being that $\Lh^{\c|c}$ and $\Lh^{\c|a}$ now depend on $\alpha$. Having obtained \eqref{eq:conv.projop_relax}, the essential idea is to write the $\alpha$-dependent operators as
 \begin{align}
     \Lh^{\c|\c}(\alpha) =
     \begin{pmatrix} \alpha\Lh^{\ip|\c} \\ \Lh^{\c\setminus\ip|\c} \end{pmatrix},
  && \Lh^{\c|\a}(\alpha) =
     \begin{pmatrix} \alpha\Lh^{\ip|\a} \\ \Lh^{\c\setminus\ip|\a} \end{pmatrix}.
 \end{align}
 The next step is then to evaluate $\G^\cc\Lh^\cc(\alpha)$ which will be equivalent to $\G^\cc\Lh^\cc$ plus the remainder $\G^\cip\Lh^{\ip/\c}$, where $\G^\cip\Lh^{\ip/\p}$ is the finite boundary contribution occurring in $\P^\cp_3$. 
\end{prf}

From the structure of the iteration operator $\P(\alpha)$ it can be immediately deduced from Theorem \ref{thm:conv_rate} that the convergence rate of the relaxed Sinclair method depends on the spectral radius $\sigma$ of $\P^\pp(\alpha)$. Thus, the \emph{optimal relaxation parameter} is the one which minimizes $\sigma(\P^\pp(\alpha))$, that is,
\begin{equation}\label{eq:conv.alpha_opt}
 \alpha_\mrm{opt} := \arg \left\{ \underset{\alpha}{\min} \; \sigma(\P^\pp(\alpha)) \right\}.
\end{equation}
The behavior of $\sigma(\P^\pp(\alpha_\mrm{opt}))$ is exemplified in Figure \ref{fig:sigma_vs_k2}  (b) for the one-dimensional problem from Section \ref{sec:conv.stab.1d}, showing that $\sigma(\P^\pp(\alpha_\mrm{opt}))$ is more than an order of magnitude smaller than $\sigma(\P^\pp)$ in the selected interval $k_2/k_1 \in [-0.1,0]$.

\subsubsection{Dynamic relaxation}
\label{sec:dyn_relax}

For nonlinear problems, computing \emph{the} optimal relaxation parameter in advance may not be the optimal choice since the atomistic operator $\L$ potentially changes in every nonlinear iteration. It seems thus more practical to dynamically update $\alpha$ after \emph{every} global iteration.

To compute an approximation of the optimal $\alpha$, we choose to linearize the problem around $\displ^\l_{k+1/2} = \{\displ^\a_{k+1}, \displ^\c_k\}$. Let $\displ$ now be the solution to this linearized problem, it then follows from \eqref{eq:conv.proj_relax} that
\begin{equation}\label{eq:dyn_relax.linearization}
 \displ_{k+1} - \displ = \P(\displ_{k+1/2}; \alpha)[\displ_k - \displ].
\end{equation}
Equation \eqref{eq:dyn_relax.linearization} is still not practical since we do not want to solve the eigenvalue problem \eqref{eq:conv.alpha_opt} in every iteration. Therefore, we convert \eqref{eq:dyn_relax.linearization} into a problem which minimizes the difference between two \emph{iterates} by subtracting \eqref{eq:dyn_relax.linearization} from $\displ_{k+2} - \displ_k$ leading to
\begin{equation}\label{eq:dyn_relax.linearization2}
 \displ_{k+2} - \displ_{k+1} = \P(\displ_{k+1/2}; \alpha)[\displ_{k+1} - \displ_k].
\end{equation}
The \emph{optimal dynamic relaxation parameter} is then defined as the one which minimizes the maximum element of \eqref{eq:dyn_relax.linearization2}
\begin{equation}\label{eq:dyn_relax.alpha_opt_dyn}
     \alpha_\mrm{opt}^\mrm{dyn}
  := \arg{ \left\{ \underset{\alpha}{\min} \; \| \P^\pp(\displ_{k+1/2}; \alpha)[\displ_{k+1} - \displ_k] \|_{l^\infty} \right\} }.
\end{equation}

Evaluating $\alpha_\mrm{opt}^\mrm{dyn}$ yet requires the additional computation of a continuum and a linearized atomistic problem, which follows from the definition of $\P(\displ_{k+1/2}; \alpha)$ (cf. Lemma \ref{lem:projop_relax}). Clearly, to be efficient, this method thus necessitates that the assumptions
\begin{itemize}
 \item[(i)]
 solving \textbf{(H)} is significantly cheaper than solving \textbf{(AH)},
 \item[(ii)]
 solving the linearized atomistic problem is significantly cheaper than solving the fully nonlinear problem \textbf{(AH)},
\end{itemize}
hold, in which by ``significantly cheaper'' we roughly mean an order of magnitude. Assumption (i) has been verified in \citep{hodapp_lattice_2019}, where the elapsed time to solve \textbf{(H)} using an efficient $\scH$-matrix solver \citep{bebendorf_hierarchical_2008} was found to be of the same order than a \emph{single} evaluation of the atomistic force $\var{}{}{\Etot^\a}$. Assumption (ii) can be justified by considering the solution of a nonlinear problem as a sequence of many linear problems. Moreover, when using a solver which builds Hessians (or approximations thereof), we can \emph{reuse} $\L^\aa(\displ_{k+1/2})$ and $\L^\ap(\displ_{k+1/2})$. Since most of the time is usually spend on \emph{building} $\L^\aa$ and $\L^\ap$, the time for solving the linear system is well-compensated. Furthermore, we can employ the solution to the linearized atomistic problem as an initial guess to the subsequent nonlinear iteration. In this respect, we can view the proposed relaxation method as a predictor-correcter scheme in which the linearized (trial) step is used to correct the boundary condition on $\textbf{(AH)}_{k+2}$.

Algorithm \ref{algo:dyn_relax} shows the essential steps to compute $\alpha_\mrm{opt}^\mrm{dyn}$. This algorithm can be directly integrated into Algorithm \ref{algo:sinc} before line 13. In the results section, we refer to the Sinclair method with dynamic relaxation as \texttt{SincDynRelax}.

\begin{algorithm}[hbt]
 \SetAlgoSkip{bigskip}
 \LinesNumbered
 \SetKwInput{Input}{Input}
 \SetKwInput{Output}{Output}
 \SetKwBlock{Repeat}{repeat}{end}
 \setstretch{1.2}
 \setlength{\commentWidth}{0.35\textwidth}
 
 \newcommand{\atcp}[1]{\tcp*[r]{\makebox[\commentWidth]{#1\hfill}}}
 \newcommand{\nosemic}{\SetEndCharOfAlgoLine{\relax}}
 \newcommand{\dosemic}{\SetEndCharOfAlgoLine{\string;}}
 \newcommand{\pushline}{\Indp}
 \newcommand{\popline}{\Indm\dosemic}
 \let\oldnl\nl
 \newcommand{\nonl}{\renewcommand{\nl}{\let\nl\oldnl}}
 
 \caption{Dynamic relaxation (\texttt{DynRelax})}
 \label{algo:dyn_relax}
 \Input{optimal relaxation parameter $\alpha^\mrm{dyn}_{\mrm{opt},k}$ and solution $\displ_{k+1/2}$ from previous iteration,\newline
        inhomogeneous force $\force_{\mrm{inh},k+1}$, natural boundary condition $\displ^\I$}
  $       \force^\ip_{\mrm{inh},k+1} \,\leftarrow\, \alpha^\mrm{dyn}_{\mrm{opt},k}\force^\ip_{\mrm{inh},k+1}$
    \atcp{relax inhomogeneous force}
  $\displ^\p_{\mrm{trial},k+1} \,\leftarrow\, \displ^\p_k + \texttt{DBEM}(\bar{\displ}^\I - \displ^\I_k, -\force_\mrm{inh}^\ip)$
    \atcp{compute trial solution}
  \nosemic                      $w^\p_1 \,\leftarrow\, \P_1^\pp[\displ_{\mrm{trial},k+1} - \displ_k],
                          \qquad w^\p_3 \,\leftarrow\, \P_3^\pp[\displ_{\mrm{trial},k+1} - \displ_k],$ \\
  \pushline\dosemic\nonl $\qquad w^\p_2 \,\leftarrow\, \P_2^\pp(\displ_{k+1/2})[\displ_{\mrm{trial},k+1} - \displ_k]$ \\
  \popline $\alpha_{\mrm{opt},k+1}^\mrm{dyn} \,\leftarrow\,
            \arg{ \left\{ \underset{\alpha}{\min} \;
              \| (w^\p_1 + w^\p_3) + \alpha(w^\p_2 - w^\p_3) \|_{l^\infty}
            \right\} }$
    \atcp{update $\alpha$} 
 \Output{optimal dynamic relaxation parameter $\alpha_{\mrm{opt},k+1}^\mrm{dyn}$}
\end{algorithm}

\begin{rem}
 The dynamic relaxation method is optimal in the following sense
 \begin{equation}
  \begin{array}{@{}r@{{}\mathrel{}}c@{\mathrel{}{}}l@{}}
       \| \displ^\p_{k+2} - \displ^\p_{k+1} \|_{l^\infty}
   =   \| \P^\pp(\displ_{k+1/2}; \alpha_\mrm{opt}^\mrm{dyn})[\displ_{k+1} - \displ_k] \|_{l^\infty}
     & \le                                        & \| \P^\pp(\displ_{k+1/2}; \alpha_\mrm{opt}^\mrm{dyn})[\displ_{k+1} - \displ_k] \| \\
     & \overset{\eqref{eq:conv.alpha_opt}}{\le}   & \| \P^\pp(\displ_{k+1/2}; \alpha_\mrm{opt})[\displ_{k+1} - \displ_k] \| \\
     & \le                                        & \| \P^\pp(\displ_{k+1/2}; \alpha_\mrm{opt}) \| \| \displ^\p_{k+1} - \displ^\p_k \| \\
     & \lesssim                                   & \sigma(\P^\pp(\displ_{k+1/2}; \alpha_\mathrm{opt})) \| \displ^\p_{k+1} - \displ^\p_k \|.
  \end{array}
 \end{equation}
\end{rem}


\begin{new}%
\section{Numerical examples}
\label{sec:examples}

In this section, we present some selected numerical experiments for a linear and a nonlinear problem. Thereby, we focus solely on the convergence properties of the Sinclair method in order to validate Theorem \ref{thm:conv_rate}. A more detailed analysis of the modeling error in force-based A/C coupling methods can, e.g., be found in \citep{parks_connecting_2008,luskin_atomistic--continuum_2013}.

In the following we consider two-dimensional problems. For various examples in one dimension the reader is referred to the preprint \citep{hodapp_analysis_2019}.

Let $\latInf$ be a hexagonal lattice with basis vectors
\begin{align}
 v_1 = \begin{pmatrix} \latConst & 0 \end{pmatrix}^\sT,
 &&
 v_2 = 1/2 \begin{pmatrix} \latConst & \sqrt{3}\latConst \end{pmatrix}^\sT,
\end{align}
where $\latConst$ is the lattice constant. Our computational domain $\lat$ is the intersection of $\latInf$ with a disc with radius $r$ and the decomposition of $\lat$ into an atomistic and a continuum domain is given as follows
\begin{align}
 \latA := \{\, x \in \real^2 \,\vert\, \| x \| \le r_\rma \,\} \cap \latInf,
 &&
 \latC := \{\, x \in \real^2 \,\vert\, r_\rma < \| x \| \le r \,\} \cap \latInf.
\end{align}

Atomic interaction is described by a Morse potential \citep{morse_diatomic_1929} with corresponding site energy
\begin{equation}
 \Eato(\{\displ(\atoB) - \displ(\ato)\}) = \sum_{\atoB \in \intRgAto} D e^{-2a\big(d(\displ(\atoB) - \displ(\ato)) - r_0\big)} - 2D e^{-a\big(d(\displ(\atoB) - \displ(\ato)) - r_0\big)},
\end{equation}
where $d(\displ(\atoB) - \displ(\ato)) = \| (\displ(\atoB) - \displ(\ato)) + (\atoB - \ato) \|$, with $D$, $a$ and $r_0$ being free parameters. The parameters of the Morse potential are chosen as follows
\begin{align}
 D = 1, && a = 4.4, && r_0 = 1.
\end{align}
It should be noted that these parameters are chosen for testing purposes and not to model a specific material. The same set of parameters has been used for previous benchmark problems of the flexible boundary condition method for infinite problems \citep{hodapp_lattice_2019}. For our numerical tests we use a cut-off radius up to six nearest neighbors. The lattice constant is then given by $\latConst = 0.978$. This implies that the atomistic model is slightly nonlocal. More precisely, the ratio between the maximum absolute values of force constants between second and first nearest neighbors is $\approx$\,1/43, which is a realistic value for many metals.

For this potential we construct our ``continuum site energy'' \eqref{eq:con.Etot_h} by triangulating the hexagonal lattice and then use the standard $\mathbb{P}1$ interpolation to define $\displ(x)$.

In all of the following numerical experiments the energies and forces are computed using the molecular dynamics code LAMMPS (\href{https://lammps.sandia.gov}{lammps.sandia.gov}).

\subsection{Linear problem: Unit point force at the origin}
\label{sec:examples.pf}

Our first goal is to validate Theorem \ref{thm:conv_rate} for a fully harmonic problem when the atomistic site energy is linearized around the ground state. We consider a problem subject to a unit point force applied on the atom at $\begin{pmatrix} 0 & 0 \end{pmatrix}^\sT$. Thereby, the point force is rescaled such that the solution remains in the linear regime. For the test problem in this section we set $r_\rma = 5\latConst$ and $r = 50\latConst$. As a boundary condition on the outer boundary $\lat^\I$ the continuum Green function is used.

We compare two variants of the Sinclair method, that is, with the relaxation parameter $\alpha = 1$, denoted by \texttt{Sinc}, and with $\alpha = \alpha_\mrm{opt}$ (eq. \eqref{eq:conv.alpha_opt}), denoted by \texttt{SincRelax}. From Figure \ref{fig:pf_err_vs_it} it can be seen that the computed convergence rates $\sigma$ and $\sigma_\mrm{opt}$ are sharp upper bounds, where $\sigma_\mrm{opt}$ is the rate obtained when using the optimal relaxation parameter $\alpha_\mrm{opt}$. The small difference is due to the fact that atomistic problem is not very nonlocal.

\begin{figure}[hbt]
 \centering
 \includegraphics[width=0.4\textwidth]{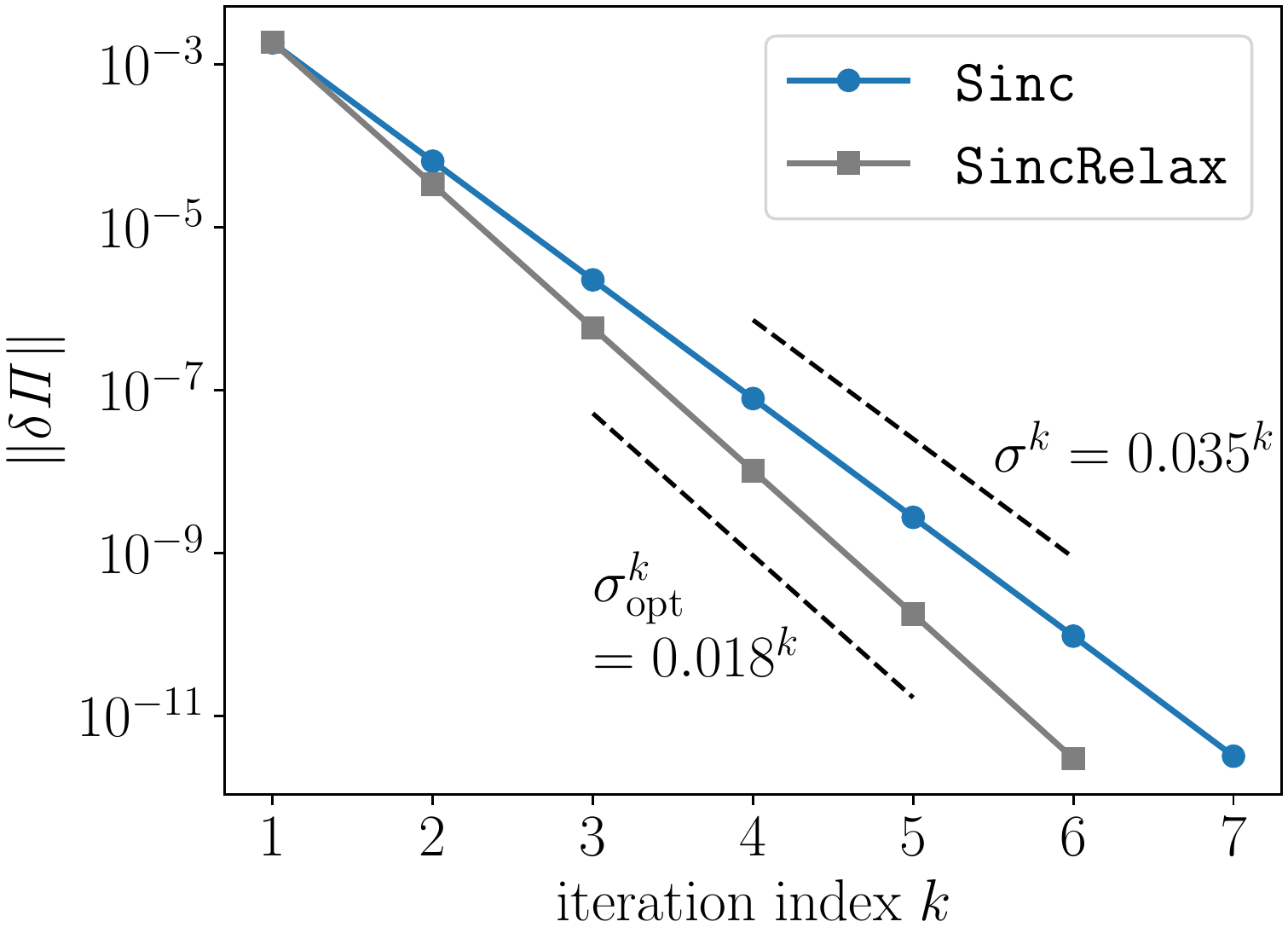}
 \caption{Convergence behavior of the Sinclair method for the linear problem}
 \label{fig:pf_err_vs_it}
\end{figure}

\subsection{Nonlinear problem: Microcrack}
\label{sec:examples.crack}

Our second goal is to validate Theorem \ref{thm:conv_rate} when the problem is linearized around some arbitrary state. To that end, we construct a nonlinear problem and check whether the observed asymptotic convergence rate agrees with the computed one, that is, the spectral radius $\sigma_\mrm{as}$ of $\P(\displ)$ when linearizing around the final solution $\displ$. 

To benchmark this behavior we now consider an effectively \emph{infinite problem} ($r \rightarrow \infty$) without external force but with some atoms removed as shown in Figure \ref{fig:domain_crack} (a), creating some kind of microcrack. The same test problem has been considered in \citep{shapeev_consistent_2011}. Upon relaxation the atoms will move to the positions indicated by the dashed circles.

\begin{figure}[hbt]
 \begin{minipage}{0.5\textwidth}\centering
  (a)
 \end{minipage}\hfill
 \begin{minipage}{0.5\textwidth}\centering
  (b)
 \end{minipage}\hfill\\[0.5em]
 \begin{minipage}{0.5\textwidth}
  \centering
  \includegraphics[width=0.54\textwidth]{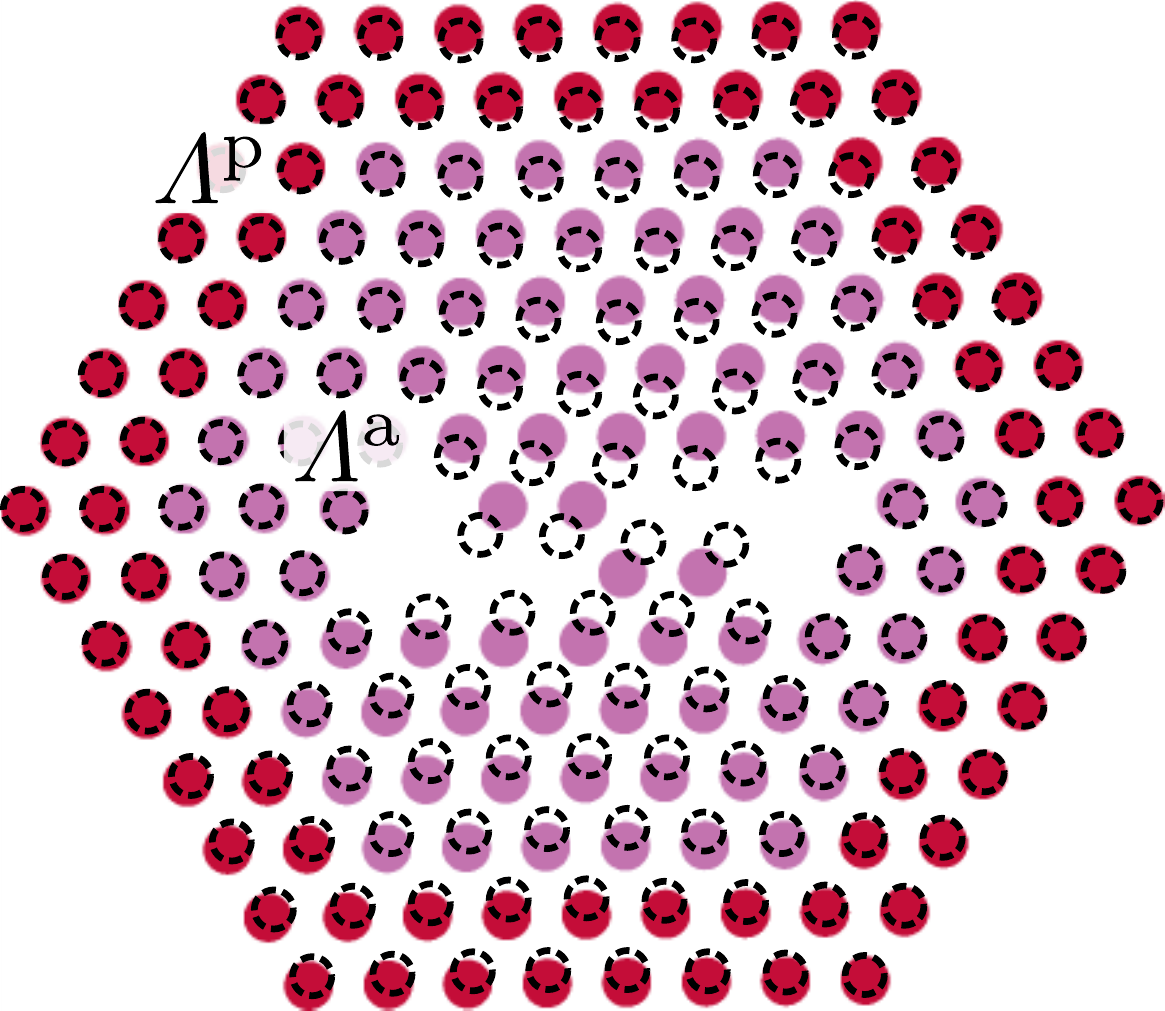}
 \end{minipage}\hfill
 \begin{minipage}{0.5\textwidth}
  \centering
  \includegraphics[width=0.8\textwidth]{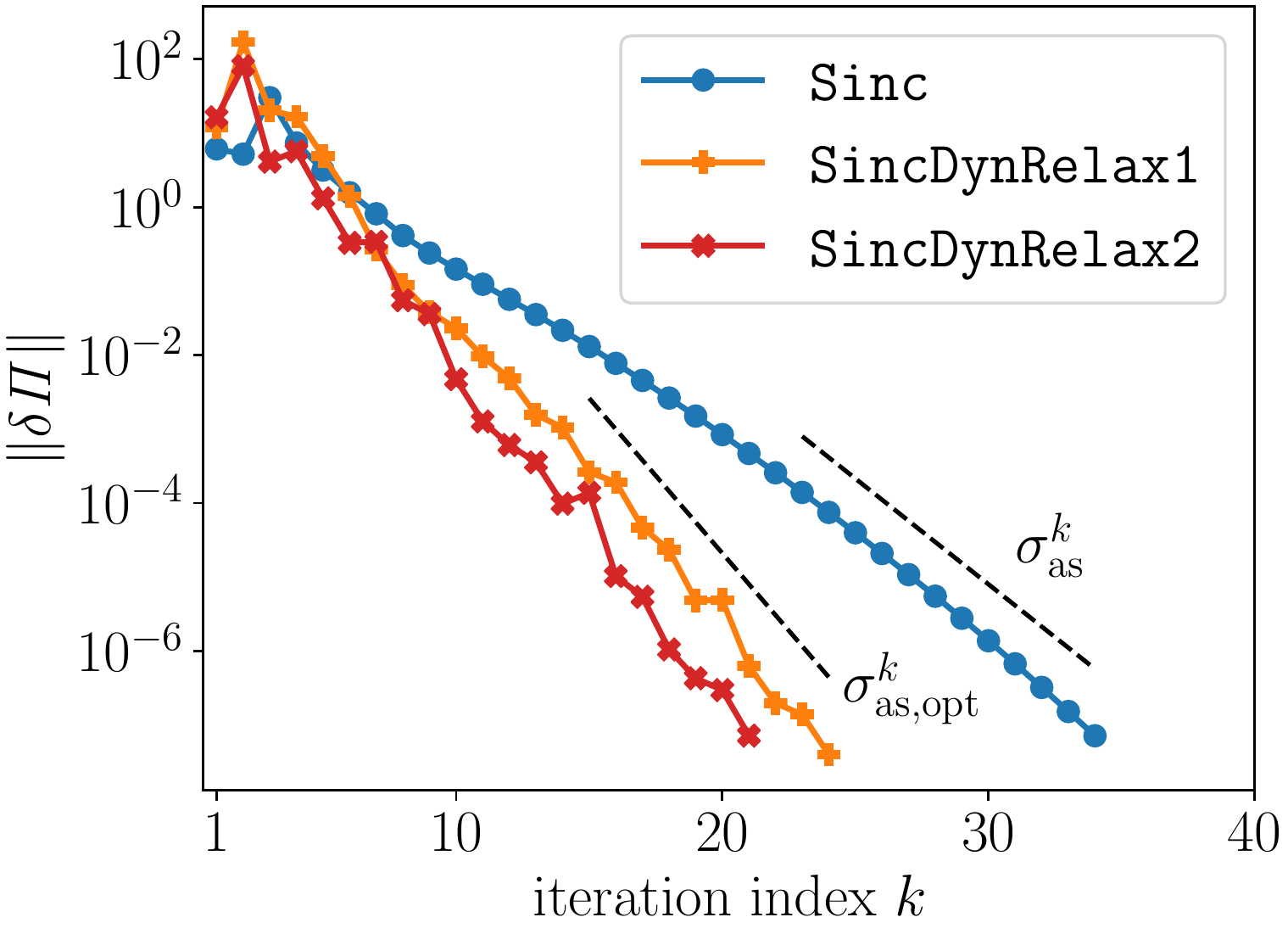}
 \end{minipage}
 \caption{(a) Atomistic domain domain $\latA$ (here shown for $r_\rma = 5\latConst$) and pad $\latP$ for the microcrack problem (the dashed circles refer to the relaxed atomic positions). (b) Convergence behavior of the different variants of the Sinclair method for the microcrack problem when $r_\rma = 4\latConst$}
 \label{fig:domain_crack}
\end{figure}

The nonlinear problem is then computed using Algorithm \ref{algo:sinc} with the atomistic problem being solved by means of the Hessian-free Newton-Raphson method implemented in LAMMPS. Further improvements of Algorithm \ref{algo:sinc} have been made by not relaxing the atomistic problem (line 8) to very low tolerances in each global iteration $k$. To be more precise, let $TOL$ be the global tolerance on the force norm $\| \var{}{}{\Etot}(\displ_{k+1}) \|$ and, when solving the atomistic problem, let $TOL_\rma$ be the tolerance on the atomistic force norm $\| \var{}{}{\EtotA}(\displ_{k+1/2}) \|$. We then start with $TOL = 10^0$ and $TOL_\rma = 10^{-2}$. If $TOL$ has been reached we set $TOL = 10^{-1}$ and $TOL_\rma = 10^{-3}$. We proceed in this way until $TOL = 10^{-7}$ and $TOL_\rma = 10^{-9}$. By not relaxing the atomistic problem each time to $TOL_\rma = 10^{-9}$ the number of force evaluations could be reduced by a factor of 2--3.

In addition to \texttt{Sinc}, we consider two variants of the Sinclair method with dynamic relaxation, one without updating $\Lhnl^\aa$ and $\Lhnl^\ap$ after the initial guess, denoted by \texttt{SincDynRelax1}, and one which updates $\Lhnl^\aa$ and $\Lhnl^\ap$ every global iteration $k$, denoted by \texttt{SincDynRelax2}. The number of iterations $N_\mrm{iter}$ and force evaluations $N_{\var{}{}{\Etot^\a}}$ for each variant are shown in Table \ref{tab:results_crack} for different sizes of the atomistic domain, alongside with the (optimal) asymptotic convergence rate and the energy error $\vert \Etot_\mrm{ref} - \Etot \vert$ with respect to the reference solution $\Etot_\mrm{ref}$.

\begin{table}[htb]
 \centering
 \begin{new}%
 \begin{tabular}{|c|c|c|c|c|c|c|c|c|c|}
  \hline
  &&&
  & \multicolumn{2}{c|}{\texttt{Sinc}}
  & \multicolumn{2}{c|}{\texttt{SincDynRelax1}}
  & \multicolumn{2}{c|}{\texttt{SincDynRelax2}} \\
  $\#\latA$ ($r_\rma$) & $\sigma_\mrm{as}$ & $\sigma_\mrm{as,opt}$ & $\vert \Etot_\mrm{ref} - \Etot \vert$
  & $N_\mrm{iter}$ & $N_{\var{}{}{\Etot^\a}}$
  & $N_\mrm{iter}$ & $N_{\var{}{}{\Etot^\a}}$
  & $N_\mrm{iter}$ & $N_{\var{}{}{\Etot^\a}}$ \\ \hline\hline
  53 ($4\latConst$) & 0.52 & 0.38 & 0.3
  & 34 & 3384 & 24 & 2591 & 21 & 2267 \\ \hline
  83 ($5\latConst$) & 0.42 & 0.3 & 0.2
  & 24 & 2684 & 16 & 1909 & 15 & 1666 \\ \hline
  179 ($7\latConst$) & 0.29 & 0.18 & 0.064
  & 17 & 2260 & 11 & 1589 & 11 & 1520 \\ \hline
 \end{tabular}
 \end{new}%
 \caption{converge rates, energy error and number of iterations and force evaluations for different variants of the Sinclair method}
 \label{tab:results_crack}
\end{table}

All variants converge faster with increasing domain size which is reasonable since the elastic solution is then expected to be more accurate on the interface $\latI$ reducing the spectral radius of $\P$ (cf. Section \ref{sec:conv.stab.general}). The performance of \texttt{SincDynRelax2} is slightly better than for \texttt{SincDynRelax1} in terms of the required number of iterations and force evaluations. Whether it is also more efficient depends, however, on the implementation (cf. Section \ref{sec:dyn_relax}). Moreover, in Figure \ref{fig:domain_crack} (b), the force norm is shown as a function of the iteration index $k$ for the coupled problem with $\#\latA$=53. We observe that the computed asymptotic convergence rates are sharp.

In order to demonstrate the efficiency of the Sinclair method, we compare the results with a fully atomistic simulation with the far-field boundary condition set to zero. To achieve approximately the same accuracy than the coupled problem with $\#\latA$=179, a domain with radius $r$=$37.5\latConst$ and a total number of 5117 atoms was necessary. Using the same Newton-Raphson solver as described above, solving this problem required 1302 force evaluations until $\| \var{}{}{\Etot} \| < 10^{-7}$. The ratio of the per-atom force evaluations (number of atoms $\times$ number of force evaluations) between the fully atomistic and the coupled problem with $\#\latA$=179 using \texttt{SincDynRelax1} is then $5117 \cdot 1302 / 179 \cdot 1589 \approx$\,23.4. Given that the ratio of the domain sizes is 5117/179\,$\approx$\,28.6, the Sinclair method can be considered as quite efficient when used as a standalone solver.
\end{new}%


\section{Conclusions}
\label{sec:summary}

We have developed and analyzed a new domain decomposition solver for force-based atomistic/continuum coupling. The proposed solver extends the method of Sinclair \citep{sinclair_flexible_1978}, developed in the 1970s for effectively infinite problems, to bounded domains. The novelty of the proposed method is the splitting of the global differential operator into a local anharmonic and an global harmonic part which stands in contrast to classical methods which partition the problem into separated (local) domains. We have analyzed the convergence properties of the method and shown that this splitting gives rise to significantly improved convergence rates over classical methods (e.g., alternating Schwarz). Moreover, a relaxation method has been proposed further reducing the number of necessary iterations.

The present work also incorporates several practical aspects. First, we have developed an implementation of the global harmonic problem using a discrete boundary element method which does not require an explicit discretization of the interior domain. In a previous publication \citep{hodapp_lattice_2019}, we have shown that this method can be efficiently combined with hierarchical approximations of the dense system matrices, reducing the computational complexity to $\# DOF\log{(\# DOF)}$ and allowing for large-scale simulations with possibly hundreds of thousands of real atoms. Second, we have shown that the harmonic problem does not require a priori knowledge of the exact behavior in the atomistic domain. Although an initial guess can improve the pre-asymptotic convergence behavior, it is not crucial. Third, we generally remark that domain decomposition solvers share the attractive advantage to be easily integrable into existing molecular dynamics codes which is a crucial requirement for practical application.

However, we have not discussed various other points which require additional attention, for example:
\begin{itemize}
 \item Stability analysis in two and three dimensions.
 \item Application of other acceleration techniques (e.g., overlapping subdomains).
 \item Adequacy of the method to be used as a preconditioner for monolithic Krylov subspace solvers.
\end{itemize}
Currently we are working on an efficient three-dimensional implementation including an approximation of the solution on the outer boundary (see Remark \ref{rem.outer_boundary}) in combination with an $\scH$-matrix solver. Results will reported in a future article.

Moreover, it is worth noting that the proposed method is not limited to atomistic/continuum coupling and can likewise be used to solve related multi-domain problems, e.g., hybrid quantum/classical mechanics problems (cf. \citep{woodward_flexible_2002}).


\section*{Acknowledgements}

Support of this project through partial funding by the Skoltech Next Generation Program (Skoltech-MIT NGP grant \#2016-7) and the Fond National Suisse (projects 200021-140506 and 191680) is highly acknowledged.

The author thanks both reviewers for their very detailed and constructive review which undoubtedly helped to improve the quality of this article. Moreover, the author would like to thank Prof. Alexander Shapeev for several valuable comments and suggestions.


\section*{Appendix}

\begin{appendices}

\begin{moved}%
\section{Proof of Theorem \ref{thm:rho_sinc}}
\label{sec:apdx.proof_thm_rho_sinc}

\begin{movedAndModified}%
\noindent%
In the appendix we will abbreviate the linearized atomistic operator $\Lhnl$ and write simply $\L$.

The strategy of the proof is (1) to derive an explicit representation of $\sigma(\P^\pp)$ in terms of the force constants $k_1$ and $k_2$, the size of the atomistic domain $M$, and elements of the inverse matrix corresponding to the atomistic Hessian $\L^{\a|\a}$. Then, in part (2), we will estimate $\sigma(\P^\pp)$ \emph{without} explicitly computing these inverse matrix elements. To that end, we require two preliminary results, Proposition \ref{prop:conv_rate_1d} and \ref{prop:inv_coeffs}, which are stated below. The section concludes with the proof of Theorem \ref{thm:rho_sinc}.

For part (1) of the proof we will use the fact that, for the one-dimensional system, we can show that the convergence rate only depends on the spectral properties of the inhomogeneous operator \Pinf (eq. \eqref{eq:conv.Tinf}). This considerably simplifies the evaluation of $\sigma(\P^\pp)$ since we do not have to be concerned with $\Pbnd$. We remark that, since we are solely interested in the convergence \emph{rate}, we consider only the convergence to solutions of $\displ$ of the coupled problem (eq. \ref{eq:conv.euler_lagrange}); and not the reference solution $\displ_\mrm{ref}$.
\end{movedAndModified}%

\begin{prop}\label{prop:conv_rate_1d}
 Let $\displ$ be the solution to the one dimensional problem \eqref{eq:conv.euler_lagrange}. For this problem, the convergence rate does not depend on the boundary operator \eqref{eq:conv.Schur_inv}, that is, we can bound the error in the $k$+1-th iteration as
 \begin{equation}\label{eq:conv.bound_1d}
  \| \displ - \displ_{k+1} \| \le C_1 \max{\{\tilde{\sigma}^{k+1}, C_2\tilde{\sigma}^k\}} \| \displ - \displ_0 \|,
 \end{equation}
 where $\tilde{\sigma}$ is the spectral radius of the inhomogeneous operator $\Pinf$ \eqref{eq:conv.Tinf} and $C_1,C_2 > 0$ are constants, independent of $k$.
\end{prop}

\begin{prf}
 See Appendix \ref{apdx:proof_prop_Tinf}.
\end{prf}

\begin{movedAndModified}%
For part (2) of the proof we require several results concerning estimates of the elements of the matrix corresponding to $\inv{\L^\aa}$.

For this purpose, we will introduce some additional notation. In what follows, we will denote a real matrix of size $N \times M$ by italic letters with two underscores, e.g, $\uuA \in \real^{N \times M}$. An individual element $(i,j)$ of this matrix $\uuA$ will be denoted by $(\uuA)_{i,j}$ or occasionally, for compactness, $A_{i,j}$ if it can be done unambiguously. The same holds for vectors. That is, a real vector of size $N$ will be denoted by italic letters with one underscore, e.g., $\uv \in \real^N$, and an individual element $(i)$ of this vector is denoted by $(\uv)_i$ or $v_i$.

Using this notation, we will denote $\uuL^\aa$ as the matrix corresponding to the operator $\L^\aa$. For the one-dimensional system in Section \ref{sec:conv.stab.1d} $\uuL^\aa$ is an $N^\a \times N^\a$ matrix, where $N^\a = 2M + 1$. We define $\uuL^\aa$ as follows
\begin{equation}\label{eq:apdx.Laa}
 \uuL^\aa =
 \begin{pmatrix}
  L^\aa(-M,-M)   & \cdots & L^\aa(-M,M) \\
  \vdots         & \ddots & \vdots \\
  L^\aa(M,-M)    & \cdots & L^\aa(M,M) \\
 \end{pmatrix},
 \qquad with \quad
 L^\aa_{\sst i,j} =
 \left\{
 \begin{aligned}
  \; 2k  &\quad&& \text{if} \; i = j, \\
    -k_1 &     && \text{if} \; \lvert i - j \rvert = 1, \\
    -k_2 &     && \text{if} \; \lvert i - j \rvert = 2, \\
    0    &     && \text{else},
 \end{aligned}
 \right.
\end{equation}
where $k = k_1 + k_2$. Consequently, the matrix corresponding to the inverse operator $\inv{\L^\aa}$ is $\inv{\uuL^\aa}$.
\end{movedAndModified}%

\begin{prop}\label{prop:inv_coeffs}
 Let $\L^\aa$ be given via \eqref{eq:conv.euler_lagrange} and $\uuL^\aa$ its matrix representation as defined in eq. \eqref{eq:apdx.Laa}. With $k_2 \in (-k_1/4,0)$ it then holds
 \begin{itemize}
  \item[\text{(a)}] Positivity of the inverse coefficients: $L^{\aa^{\st -1}}_{\sst 1,1}, \, L^{\aa^{\st -1}}_{\sst 1,N^\a} > 0$.
  \item[\text{(b)}] Ratio of the inverse coefficients:
  \begin{multicols}{2}
   \begin{itemize}[leftmargin=1cm]
    \item[\text{(b.1)}] \, $\dfrac{L^{\aa^{\st -1}}_{\sst 1,N^\a-1}}{L^{\aa^{\st -1}}_{\sst 1,N^\a}} < 3 - \dfrac{3}{2M+1}$,
    \item[\text{(b.2)}] \, $\dfrac{L^{\aa^{\st -1}}_{\sst 1,1}}{L^{\aa^{\st -1}}_{\sst 1,N^\a}} > 1$.
   \end{itemize}
  \end{multicols}
 \end{itemize}
\end{prop}

\begin{prf}
 See Appendix \ref{apdx:proof_prop_inv_coeffs}.
\end{prf}

With the previous results we can now prove Theorem \ref{thm:rho_sinc}:

\begin{prf}[of Theorem \ref{thm:rho_sinc}]
 First, let us recall that we only need to analyze the spectral radius of the inhomogeneous iteration operator $\Pinf^\pp$, defined in \eqref{eq:conv.Tinf}, which follows from Theorem \ref{thm:conv_rate} and Proposition \ref{prop:conv_rate_1d}.
 
 Therefore, we begin by defining the matrix representation $\tilde{\uuT}^\pp$ corresponding to $\Pinf^\pp$. Since the atomistic model comprises second-nearest neighbor interactions $\tilde{\uuT}^\pp$ is a 4$\times$4 matrix (since there are two pad atoms on both sides; cf. Figure \ref{fig:domain_decomp_1d}). We define $\tilde{\uuT}^\pp$ as follows
 \begin{equation}\label{eq:conv.stab_proof.projop}
  \begin{aligned}
   \tilde{\uuT}^\pp
   &=
   \begin{pmatrix}
    \tilde{T}(-M,-M)   & \tilde{T}(-M,-M+1)   & \tilde{T}(-M,M-1)   & \tilde{T}(-M,M) \\
    \tilde{T}(-M+1,-M) & \tilde{T}(-M+1,-M+1) & \tilde{T}(-M+1,M-1) & \tilde{T}(-M+1,M) \\
    \tilde{T}(M-1,-M)  & \tilde{T}(M-1,-M+1)  & \tilde{T}(M-1,M-1)  & \tilde{T}(M-1,M) \\
    \tilde{T}(M,-M)    & \tilde{T}(M,-M+1)    & \tilde{T}(M,M-1)    & \tilde{T}(M,M)
   \end{pmatrix}.
  \end{aligned}
 \end{equation}
 To evaluate the elements of $\tilde{\uuT}^\pp$, we require the lattice Green function $G(\ato - \atoB)$ (eq. \eqref{eq:con.GLh=I}). For a one-dimensional lattice it is given by $G(\ato - \atoB) = - \vert\, \ato -\atoB \,\vert / 2\bar{k}$. Using the interactions stencils \eqref{eq:conv.euler_lagrange} and the Green function, we then construct the matrices corresponding to the $\greenOp$- and $\diffOp$-operators in \eqref{eq:conv.Tinf} and compute $\tilde{\uuT}^\pp$. After performing some algebra, we can write the elements of $\tilde{\uuT}^\pp$ as
 \begin{equation}\label{eq:conv.stab_proof.projop_coeffs}
  \begin{aligned}
   \tilde{T}_{1,1} = \tilde{T}_{4,4} &= -k_2(L^{\aa^{\st -1}}_{\sst 1,1} + (2M + 3)L^{\aa^{\st -1}}_{\sst 1,N^\a})/2, \\
   \tilde{T}_{1,2} = \tilde{T}_{4,3} &=   -k_1(L^{\aa^{\st -1}}_{\sst 1,1} + (2M + 3)L^{\aa^{\st -1}}_{\sst 1,N^\a})/2
                        -  k_2(L^{\aa^{\st -1}}_{\sst 1,2} + (2M + 3)L^{\aa^{\st -1}}_{\sst 1,N^\a-1})/2 + 1, \\
   \tilde{T}_{1,3} = \tilde{T}_{4,2} &=   -k_1((2M + 3)L^{\aa^{\st -1}}_{\sst 1,1} + L^{\aa^{\st -1}}_{\sst 1,N^\a})/2
                        -  k_2((2M + 3)L^{\aa^{\st -1}}_{\sst 1,2} + L^{\aa^{\st -1}}_{\sst 1,N^\a-1})/2 + M + 1, \\
   \tilde{T}_{1,4} = \tilde{T}_{4,1} &= -k_2((2M + 3)L^{\aa^{\st -1}}_{\sst 1,N^\a} + L^{\aa^{\st -1}}_{\sst 1,N^\a})/2, \\
   \tilde{T}_{2,1} = \tilde{T}_{3,4} &= -(M+1)k_2L^{\aa^{\st -1}}_{\sst 1,N^\a}, \\
   \tilde{T}_{2,2} = \tilde{T}_{3,3} &= -(M+1)k_2L^{\aa^{\st -1}}_{\sst 1,N^\a-1} - (M+1)k_1L^{\aa^{\st -1}}_{\sst 1,N^\a} + 1/2, \\
   \tilde{T}_{2,3} = \tilde{T}_{3,2} &= -(M+1)k_1L^{\aa^{\st -1}}_{\sst 1,1} - (M+1)k_2L^{\aa^{\st -1}}_{\sst 1,2} + M + 1/2, \\
   \tilde{T}_{2,4} = \tilde{T}_{3,1} &= -(M+1)k_2L^{\aa^{\st -1}}_{\sst 1,1}.
  \end{aligned}
 \end{equation}
 Using \eqref{eq:conv.stab_proof.projop_coeffs}, it can be readily shown that the matrix has rank 2. Moreover, $\tilde{\uuT}^\pp$ is centrosymmetric and, therefore, its spectrum are the eigenvalues of the 2$\times$2 matrices (cf. \citep{cantoni_eigenvalues_1976})
 \begin{align}
  \uuC_1 =
  \begin{pmatrix}
   \tilde{T}_{1,1} - \tilde{T}_{1,4} & \tilde{T}_{1,2} - \tilde{T}_{1,3} \\
   \tilde{T}_{2,1} - \tilde{T}_{2,4} & \tilde{T}_{2,2} - \tilde{T}_{2,3}
  \end{pmatrix}, &&
  \uuC_2 =
  \begin{pmatrix}
   \tilde{T}_{1,1} + \tilde{T}_{1,4} & \tilde{T}_{1,2} + \tilde{T}_{1,3} \\
   \tilde{T}_{2,1} + \tilde{T}_{2,4} & \tilde{T}_{2,2} + \tilde{T}_{2,3}
  \end{pmatrix}.
 \end{align}
 From \eqref{eq:conv.stab_proof.projop_coeffs} it follows $\tilde{T}_{1,1} - \tilde{T}_{1,4} = \tilde{T}_{2,1} - \tilde{T}_{2,4}$ and $\tilde{T}_{1,2} - \tilde{T}_{1,3} = \tilde{T}_{2,2} - \tilde{T}_{2,3}$. The eigenvalues of $\uuC_1$ are then given by
 \begin{align}
  \lambda_{11} = (\tilde{T}_{1,1} - \tilde{T}_{1,4}) + (\tilde{T}_{1,2} - \tilde{T}_{1,3}), && \lambda_{12} = 0.
 \end{align}
 Hence, at least one eigenvalue of $\uuC_2$ must be equal to zero which implies that the possible nonzero eigenvalue is given by the trace of $\uuC_2$. Thus, the eigenvalues of $\uuC_2$ are
 \begin{align}
  \lambda_{21} = (\tilde{T}_{1,1} + \tilde{T}_{1,4}) + (\tilde{T}_{2,2} + \tilde{T}_{2,3}), && \lambda_{22} = 0.
 \end{align}
 Exploiting the identity
 \begin{equation}\label{eq:conv.stab_proof.id_hom_def}
  -(M+1)(k_1(L^{\aa^{\st -1}}_{\sst 1,1} - L^{\aa^{\st -1}}_{\sst 1,N^\a}) + k_2(L^{\aa^{\st -1}}_{\sst 1,2} - L^{\aa^{\st -1}}_{\sst 1,N^\a-1})
  -(M+2)k_2(k_1(L^{\aa^{\st -1}}_{\sst 1,1} - L^{\aa^{\st -1}}_{\sst 1,N^\a})
  = -M,
 \end{equation}
 which can be obtained by applying a homogeneous deformation to the crystal, the eigenvalues can be written as
 \begin{align}
  \lambda_{11} &= -k_2 (L^{\aa^{\st -1}}_{\sst 1,1} + L^{\aa^{\st -1}}_{\sst 1,N^\a}) + 2k_2 L^{\aa^{\st -1}}_{\sst 1,N^\a}, \\
  \label{eq:conv.stab_proof.lambda21}
  \lambda_{21} &= -k_2 (L^{\aa^{\st -1}}_{\sst 1,1} + L^{\aa^{\st -1}}_{\sst 1,N^\a}),
 \end{align}
 from which it follows that $\lambda_{11} < \lambda_{21}$. Using the statements (a) and (b.2) from Proposition \ref{prop:inv_coeffs}, we immediately find that $\lambda_{11} > 0$ and, therefore, it remains to check whether $\lambda_{21} = \sigma < 1$.
 
 For this purpose, using the identities \eqref{eq:conv.stab_proof.id_hom_def} and
 \begin{equation}
      (k_1 + k_2)(L^{\aa^{\st -1}}_{\sst 1,1} + L^{\aa^{\st -1}}_{\sst 1,N^\a})
    + k_2(L^{\aa^{\st -1}}_{\sst 1,2} + L^{\aa^{\st -1}}_{\sst 1,N^\a-1})
  = 1
 \end{equation}
 in \eqref{eq:conv.stab_proof.lambda21}, $\lambda_{21}$ becomes
 \begin{equation}
  \begin{aligned}
   \lambda_{21} &=   (M+1) \left(1 - 2 \left(   k_2 L^{\aa^{\st -1}}_{\sst 1,N^\a-1}
                                            + (k_1 + k_2) L^{\aa^{\st -1}}_{\sst 1,N^\a} \right) \right)
                   - 2k_2 L^{\aa^{\st -1}}_{\sst 1,N^\a} - M \\ 
                &=   1 - 2(M+1) \left( k_2 L^{\aa^{\st -1}}_{\sst 1,N^\a-1}
                                       + (k_1 + k_2) L^{\aa^{\st -1}}_{\sst 1,N^\a} \right)
                    - 2k_2 L^{\aa^{\st -1}}_{\sst 1,N^\a}.
  \end{aligned}
 \end{equation}
 Clearly, to satisfy stability we require that
 \begin{equation}
    2(M+1) \left(   k_2 L^{\aa^{\st -1}}_{\sst 1,N^\a-1} 
                + (k_1 + k_2) L^{\aa^{\st -1}}_{\sst 1,N^\a} \right)
  + 2k_2 L^{\aa^{\st -1}}_{\sst 1,N^\a} > 0.
 \end{equation}
 Setting $k_1 > -4k_2$ and solving for $L^{\aa^{\st -1}}_{\sst 1,N^\a-1} / L^{\aa^{\st -1}}_{\sst 1,N^\a}$, we obtain the condition
 \begin{equation}
    \frac{L^{\aa^{\st -1}}_{\sst 1,N^\a-1}}{L^{\aa^{\st -1}}_{\sst 1,N^\a}}
  < 3 - \frac{1}{M+1}.
 \end{equation}
 Invoking Proposition \ref{prop:inv_coeffs}, statement (b.1), completes the proof.
 \qed
\end{prf}

\section{Proof of Proposition \ref{prop:conv_rate_1d}}
\label{apdx:proof_prop_Tinf}

\noindent%
If we replace $\displ_\mrm{ref}$ with $\displ$, it follows immediately from Lemma \ref{lem:itop1} that $\displ - \displ_{k+1} = \P[\displ - \displ_k]$.

Using the additive split $\P = \Pinf + \Pbnd$ (eq. \eqref{eq:conv.Tinf+Tbnd}), we recall from Theorem \ref{thm:conv_rate} that the iteration operator in the $k$+1-th iteration can be written as follows (cf. eq. \eqref{eq:conv.u_k+1-u=P[u_0-u]})
\begin{equation}\label{eq:conv.projop_prod}
 \P^{k+1} = \prod_{i=1}^{k+1} \P = \prod_{i=1}^{k+1} (\Pinf + \Pbnd) = (\Pinf + \Pbnd)(\Pinf + \Pbnd) \, \cdots = (\Pinf^2  + \Pinf\Pbnd + \Pbnd\Pinf + \Pbnd^2) \, \cdots.
\end{equation}
We first analyze successive applications of $\Pbnd$ to itself. Therefore, note that an application of $\Pbnd$ to an arbitrary vector $v \in \clV$ gives (in one dimension only!) a homogeneous solution $w(\ato) = F\ato + C \; \text{in} \; \lat$, for some $F,C>0$, since the homogeneous problem does not contain any source terms. Since the coupling is consistent, applying this homogeneous solution $w$ as a boundary condition, i.e, setting $\displ^\p = w^\p$, to the atomistic problem gives the same solution in \latA, i.e., $(\inv{\L^\aa}\L^\ac)[w] = w^\a$, and, therefore, we have $(\L^\cc - \L^\ca(\inv{\L^\aa}\L^\ac))[w] = \L^\cc[w] - \L^\ca[w] = 0$. This implies that
\begin{equation}
 \forall\,v\in\clV \qquad (\Pbnd\Pbnd)[v] = 0 
\end{equation}
which can only hold iff $\Pbnd\Pbnd = \nullOp$. That is, the operator \Pbnd is nilpotent with index 2.

In addition, since \Pbnd generates homogeneous solutions, any inhomogeneous force vanishes and, therefore, it also holds $\Pinf\Pbnd = \nullOp$.

With $\Pbnd\Pbnd = \Pinf\Pbnd = \nullOp$, equation \eqref{eq:conv.projop_prod} reduces to
\begin{equation}
 \P^{k+1} = \Pinf^{k+1} + \Pbnd\Pinf^{k}.
\end{equation}
Assuming that $\Pinf^{k+1}$ has the eigendecomposition $\Pinf^{k+1} = \tilde{\clQ} \tilde{\clD}^{k+1} \tilde{\clQ}^{-1}$, we can then write the error in the $k$+1-th iteration for the one-dimensional system as
\begin{equation}
 \displ - \displ_{k+1}  = (\tilde{\clQ} \tilde{\clD}^{k+1} \tilde{\clQ}^{-1} + \Pbnd\tilde{\clQ} \tilde{\clD}^k \tilde{\clQ}^{-1}) [\displ - \displ_0],
\end{equation}
from which the upper bound follows as
\begin{equation}
 \| \displ - \displ_{k+1} \| \le (\tilde{\sigma}^{k+1} + \tilde{\sigma}^k \| \Pbnd \|) \| \tilde{\clQ} \| \| \tilde{\clQ}^{-1} \| \| \displ - \displ_0 \|.
\end{equation}
With $C_1 = \| \tilde{\clQ} \| \| \tilde{\clQ}^{-1} \|$ and $C_2 = \| \Pbnd \|$ we obtain \eqref{eq:conv.bound_1d}.
\qed

\section{Proof of Proposition \ref{prop:inv_coeffs}}
\label{apdx:proof_prop_inv_coeffs}

\subsection{Preparation}

Without loss of generality we set $k_1 = 1$ in the following. Moreover, for clarity let $K = 2M$.

\subsubsection{Representation of \texorpdfstring{$L^{\aa^{\st -1}}_{\sst i,j}$}{} in terms of cofactors of \texorpdfstring{$\uuL^\aa$}{}}

In what follows we make frequent use of the representation of the inverse coefficients
\begin{equation}\label{eq:appdx.proof_inv_coeffs.Laainv=C/detLaa}
 L^{\aa^{\st -1}}_{\sst i,j} = \frac{C_{i,j}}{\det{\uuL^\aa}}, \qquad \text{with} \; \uuC = (-1)^{i+j}\det{\uuM_{ij}},
\end{equation}
where $\uuC$ is denoted as the cofactor matrix of $\uuL^\aa$ and $\det{\uuM_{ij}}$ is the $(i,j)$-th minor of $\uuL^\aa$. The matrix $\uuM_{ij} \in \real^{K \times K}$ is obtained by removing the $i$-th row and the $j$-th column from $\uuL^\aa$.

\subsubsection{LU factorization of \texorpdfstring{$\uuM_{1,N^\a}$}{}}

By definition, the matrix $\uuM_{1,N^\a}$, obtained after removing the first row and the last column from $\uuL^\aa$, is the upper Hessenberg-Toeplitz matrix
\begin{equation}
 \uuM_{N^\a,1} =
 \begin{pmatrix}
  -1 & 2k & -1 & -k_2 &&& \\
  -k_2 & -1 & 2k & -1 & -k_2 && \\
  && \ddots &&&& \\
  && -k_2 & -1 & 2k & -1 & -k_2 \\
  &&& -k_2 & -1 & 2k & -1 & \\
  &&&& -k_2 & -1 & 2k \\
  &&&&& -k_2 & -1 \\
 \end{pmatrix}.
\end{equation}
Let $\uuM_{1,N^\a} = \uuL\uuU$, where
\begin{align}
 \uuL =
 \begin{pmatrix}
  1 &&& \\
  L_{2,1} & 1 && \\
  && \ddots & \\
  && L_{K,K-1} & 1
 \end{pmatrix},
 &&
 \uuU
  \begin{pmatrix}
  U_{1,1} & U_{1,2} & U_{1,3} & U_{1,4} &&& \\
  & \ddots &&&& \\
  && U_{K-3,K-3} & U_{K-3,K-2} & U_{K-3,K-1} & U_{K-3,K} \\
  &&& U_{K-2,K-2} & U_{K-2,K-1} & U_{K-2,K} & \\
  &&&& U_{K-1,K-1} & U_{K-1,K} \\
  &&&&& U_{K,K} \\
 \end{pmatrix}.
\end{align}
The elements of $\uuL$ and $\uuU$ can be computed using the following recursion formulas
\begin{align}\label{eq:appdx.proof_inv_coeffs.LU_rec_eqs}
    U_{1,j} = \left( \uuM_{1,N^\a} \right)_{1,j},
 && L_{i,i-1} = \frac{\left( \uuM_{1,N^\a} \right)_{i,i-1}}{U_{i-1,i-1}},
 && U_{i,j} = \left( \uuM_{1,N^\a} \right)_{i,j} - L_{i,i-1}U_{i-1,j}.
\end{align}

For our proof we further require an upper bound for the diagonal elements of $\uuU$. For this purpose we compute the first three components
\begin{equation}\label{eq:appdx.proof_inv_coeffs.U11_U22_U33}
 U_{1,1} = -1, \qquad U_{2,2} = -2kk_2 - 1, \qquad U_{3,3} = \frac{4kk_2 + k_2^2 + 1}{-2kk_2 - 1}.
\end{equation}
Subsequently, using \eqref{eq:appdx.proof_inv_coeffs.LU_rec_eqs}, we can write
\begin{equation}\label{eq:appdx.proof_inv_coeffs.UKK}
 \forall\,i = 4,...,K \qquad
 U_{i,i} = -1 + \frac{2kk_2}{U_{i-1,i-1}} - \frac{k_2^2}{U_{i-1,i-1}U_{i-2,i-2}} -1 \frac{k_2^4}{U_{i-1,i-1}U_{i-2,i-2}U_{i-3,i-3}}.
\end{equation}
From \eqref{eq:appdx.proof_inv_coeffs.U11_U22_U33} we deduce that for $k_2 \ge k_2^\ast$ and $i = 1,2,3$, $U_{i,i}(k_2) \le U_{i,i}(k_2^\ast)$. To generalize this result, we explicitly compute the diagonal components $\bar{U}_{i,i} = U_{i,i}(-1/4)$ for the limiting case when $k_2 = -1/4$, that is,
\begin{equation}
 \bar{U}_{1,1} = \frac{4}{4}, \quad \bar{U}_{2,2} = \frac{5}{8}, \quad \bar{U}_{3,3} = \frac{6}{12}, \quad \bar{U}_{4,4} = \frac{7}{16} \quad \bar{U}_{5,5} = \frac{8}{20}, \quad ..., \quad \bar{U}_{K,K} = -(1 + 3K^{-1})/4,
\end{equation}
where last expression can be readily checked by using $U_{K-1,K-1} = -(1 + 3(K-1)^{-1})/4$, $U_{K-2,K-2} = -(1 + 3(K-2)^{-1})/4$ and $U_{K-3,K-3} = -(1 + 3(K-3)^{-1})/4$ in \eqref{eq:appdx.proof_inv_coeffs.UKK}.

It follows that $\bar{U}_{K,K}$ is bounded by $\{-1,-1/4\}$. Using \eqref{eq:appdx.proof_inv_coeffs.U11_U22_U33}, it can then be shown inductively, e.g., by employing a straightforward contradiction argument, that
\begin{equation}\label{eq:appdx.proof_inv_coeffs.UKK_bound}
 \forall\, K>1 \qquad U_{K,K} < -(1 + 3K^{-1})/4.
\end{equation}

\subsection{Proof of statement \textnormal{(a)}}

Using the representation \eqref{eq:appdx.proof_inv_coeffs.Laainv=C/detLaa}, the inverse coefficients read
\begin{align}
 L^{\aa^{\st -1}}_{\sst 1,1}    = \frac{C_{1,1}}{\det{\uuL^\aa}}, &&
 L^{\aa^{\st -1}}_{\sst 1,N^\a} = \frac{C_{1,N^\a}}{\det{\uuL^\aa}}.
\end{align}
Since $\uuL^\aa$ is positive definite it follows that $\det{\uuL^\aa} > 0$. Therefore, it remains to check that both numerators are positive. For $C_{1,1}$ this is trivially true since positive definiteness of $\uuL^\aa$ holds independently of its size. For the second term we find that
\begin{equation}
 C_{1,N^\a} = (-1)^{1+N^\a} \det{\uuM_{N^\a,1}} = \det{\uuM_{N^\a,1}} = \det{\uuL\uuU} = \det{\uuU} > 0,
\end{equation}
where the latter inequality follows from the fact $\uuU \in \real^{K \times K}$ and that $\forall\,i = 1,...,K \; U_{i,i} < 0$.
\qed

\subsection{Proof of statement \textnormal{(b)}}

To prove the second statement, we introduce the shift matrix $\uuJ = \uuM^{-1}_{N^\a,1}\uuM_{N^\a-1,1}$ to obtain
\begin{equation}\label{eq:appdx.proof_inv_coeffs.-detJ*}
   \frac{L^{\aa^{\st -1}}_{\sst 1,N^\a-1}}{L^{\aa^{\st -1}}_{\sst 1,N^\a}}
 = -\frac{\det{\uuM_{N^\a-1,1}}}{\det{\uuM_{N^\a,1}}}
 = -\frac{\det{\uuM_{N^\a,1} \uuJ}}{\det{\uuM_{N^\a,1}}}
 = -\det{\uuJ}.
\end{equation}
The matrix $\uuJ$ is the upper triangular matrix
\begin{equation}
 \uuJ =
 \begin{pmatrix}
  \uuI & \uv \\
  \uNull^\sT & \det{\uuJ}
 \end{pmatrix}.
\end{equation}
Again, using the factorized representation of $\uuM_{N^\a,1}$ the determinant of $\uuJ$ follows as
\begin{equation}\label{eq:appdx.proof_inv_coeffs.detJ*(U)}
 \det{\uuJ} = -\frac{2k}{U_{K,K}} + \frac{k_2}{U_{K,K}U_{K-1,K-1}} + \frac{k_2^3}{U_{K,K}U_{K-1,K-1}U_{K-2,K-2}}.
\end{equation}
Combining \eqref{eq:appdx.proof_inv_coeffs.detJ*(U)} and \eqref{eq:appdx.proof_inv_coeffs.UKK_bound} in \eqref{eq:appdx.proof_inv_coeffs.-detJ*} we obtain the result stated in (b.1).

Statement (b.2) is a corollary of statement (b.1). To prove it, we use the fact that $\det{\uuM_{1,1}}$ can be represented as
\begin{equation}
 \det{\uuM_{1,1}} = \det{\uuM_{1,N^\a}} + \det{\uuA},
\end{equation}
with
\begin{equation}
 \uuA = \uuM_{1,1} + \uv\uw^\sT, \qquad \text{with} \;
 \uv = 
 \begin{pmatrix}
  0, ..., 0, 1 
 \end{pmatrix}^\sT, \;
 \uw = 
 \begin{pmatrix}
  1, k_2, 0, ..., 0
 \end{pmatrix}^\sT,
\end{equation}
which can be derived using the properties of determinants. Using the matrix determinant lemma, we can write
\begin{equation}
 \begin{aligned}
  \det{\uuA} &= (1 + \uw^\sT \uuM^{-1}_{1,1} \uv)\det{\uuM_{1,1}} \\
             &= \left( 1 + \left( \uuM^{-1}_{1,1} \right)_{1,K} + k_2 \left( \uuM^{-1}_{1,1} \right)_{1,K-1} \right) \det{\uuM_{1,1}} \\
             &> \left( 1 - 4k_2 \left( \uuM^{-1}_{1,1} \right)_{1,K} + k_2 \left( \uuM^{-1}_{1,1} \right)_{1,K-1} \right) \det{\uuM_{1,1}} \\
             &> 0,
 \end{aligned}
\end{equation}
\begin{samepage}
where the last inequality follows immediately from statement (b.1). Thus, we have
\begin{equation}
   \frac{L^{\aa^{\st -1}}_{\sst 1,1}}{L^{\aa^{\st -1}}_{\sst 1,N^\a}}
 = \frac{\det{\uuM_{1,1}}}{\det{\uuM_{N^\a,1}}}
 = \frac{\det{\uuM_{N^\a,1}} + \det{\uuA}}{\det{\uuM_{N^\a,1}}}
 > 1.
\end{equation}
\qed
\end{samepage}
\end{moved}%

\begin{new}%
\section{Domain indices}
\label{apdx:dom_idx}

\begin{table}[H]
 \centering
 \begin{new}%
 \begin{tabular}{|c|c|l|l|}
  \hline
  Domain & Index & Description & Introduced in \\ \hline\hline
  $\latInf$ & $\diagup$ & Infinite Bravais lattice & Section \ref{sec:atom} \\ \hline
  $\lat$ & $\bullet^\l$ & Computational domain & Section \ref{sec:atom} \\ \hline
  $\bar{\lat}$ & $\bullet^\lb$ & Computational domain + boundary & Section \ref{sec:atom} \\ \hline
  $\latA$ & $\bullet^\a$ & Atomistic domain & Section \ref{sec:acc} \\ \hline
  $\latC$ & $\bullet^\c$ & Continuum domain & Section \ref{sec:acc} \\ \hline
  $\latP$ & $\bullet^\p$ & Pad domain & Section \ref{sec:acc} \\ \hline
  $\latI$ & $\bullet^\i$ & Interface of $\latC$ to $\latA$ & Section \ref{sec:acc} \\ \hline
  $\lat^\I$ & $\bullet^\I$ & Outer interface of $\latC$ & Section \ref{sec:acc} \\ \hline
  $\lat^\cb$ & $\bullet^\cb$ & Continuum domain and outer interface ($\latC \cup \lat^\I$) & Section \ref{sec:sinc.description.iter_eq} \\ \hline
  $\lat^\ip$ & $\bullet^\ip$ & Layer of nodes where the inhomogeneous force is nonzero & Section \ref{sec:sinc.source_term} \\ \hline
  $\lat^\re$ & $\bullet^\re$ & Remainder domain $\latInf \setminus \bar{\lat}$ & Section \ref{sec:sinc.H_k+1.hom} \\ \hline
  $\lat^\Im$ & $\bullet^\Im$ & Layer of nodes in $\lat^\re$ which interact with $\lat^\I$ according to $\intRgAto^\mrm{h}$ & Section \ref{sec:sinc.H_k+1.hom} \\ \hline
  $\lat^\pr$ & $\bullet^\pr$ & Set of atoms in $\latA$ which interact with nodes from $\latP$ & Section \ref{sec:conv.proj}, Lemma \ref{lem:itop2} \\ \hline
 \end{tabular}
 \end{new}%
 \caption{Description of the discrete domains and their superscripted indices as they are used for the lattice functions and operators throughout the manuscript}
 \label{tab:dom_idx}
\end{table}
\end{new}%

\end{appendices}


\section*{References}
\bibliographystyle{elsarticle-harv}
\bibliography{references_sinclair-analysis}

\begin{thebibliography}{41}
\expandafter\ifx\csname natexlab\endcsname\relax\def\natexlab#1{#1}\fi
\expandafter\ifx\csname url\endcsname\relax
  \def\url#1{\texttt{#1}}\fi
\expandafter\ifx\csname urlprefix\endcsname\relax\def\urlprefix{URL }\fi

\bibitem[{Anciaux et~al.(2018)Anciaux, Junge, Hodapp, Cho, Molinari, and
  Curtin}]{anciaux_coupled_2018}
Anciaux, G., Junge, T., Hodapp, M., Cho, J., Molinari, J.-F., Curtin, W., Sep.
  2018. The {Coupled} {Atomistic}/{Discrete}-{Dislocation} method in 3d part
  {I}: {Concept} and algorithms. Journal of the Mechanics and Physics of Solids
  118, 152--171.
\newline\urlprefix\url{https://linkinghub.elsevier.com/retrieve/pii/S0022509617310098}

\bibitem[{Argon(2007)}]{argon_strengthening_2007}
Argon, A., Aug. 2007. Strengthening {Mechanisms} in {Crystal} {Plasticity}.
  Oxford University Press.
\newline\urlprefix\url{http://www.oxfordscholarship.com/view/10.1093/acprof:oso/9780198516002.001.0001/acprof-9780198516002}

\bibitem[{Bebendorf(2008)}]{bebendorf_hierarchical_2008}
Bebendorf, M., 2008. Hierarchical matrices: a means to efficiently solve
  elliptic boundary value problems. No.~63 in Lecture notes in computational
  science and engineering. Springer, Berlin, oCLC: ocn220011087.

\bibitem[{Benzi et~al.(2001)Benzi, Frommer, Nabben, and
  Szyld}]{benzi_algebraic_2001}
Benzi, M., Frommer, A., Nabben, R., Szyld, D.~B., Oct. 2001. Algebraic theory
  of multiplicative {Schwarz} methods:. Numerische Mathematik 89~(4), 605--639.
\newline\urlprefix\url{http://link.springer.com/10.1007/s002110100275}

\bibitem[{Cantoni and Butler(1976)}]{cantoni_eigenvalues_1976}
Cantoni, A., Butler, P., 1976. Eigenvalues and {Eigenvectors} of {Symmetric}
  {Centrosymmetrlc} {Matrlces}. Linear Algebra and its Applications 13,
  275--288.

\bibitem[{Curtin and Miller(2003)}]{curtin_atomistic/continuum_2003}
Curtin, W.~A., Miller, R.~E., 2003. Atomistic/continuum coupling in
  computational materials science. Modelling and Simulation in Materials
  Science and Engineering 11~(3), 33--68.
\newline\urlprefix\url{https://doi.org/10.1088%2F0965-0393%2F11%2F3%2F201}

\bibitem[{Dobson et~al.(2010)Dobson, Luskin, and
  Ortner}]{dobson_stability_2010}
Dobson, M., Luskin, M., Ortner, C., Jul. 2010. Stability, {Instability}, and
  {Error} of the {Force}-based {Quasicontinuum} {Approximation}. Archive for
  Rational Mechanics and Analysis 197~(1), 179--202, arXiv: 0903.0610.
\newline\urlprefix\url{http://arxiv.org/abs/0903.0610}

\bibitem[{Dobson et~al.(2011)Dobson, Luskin, and
  Ortner}]{dobson_iterative_2011}
Dobson, M., Luskin, M., Ortner, C., Sep. 2011. Iterative methods for the
  force-based quasicontinuum approximation: {Analysis} of a {1D} model problem.
  Computer Methods in Applied Mechanics and Engineering 200~(37-40),
  2697--2709.
\newline\urlprefix\url{https://linkinghub.elsevier.com/retrieve/pii/S0045782510002203}

\bibitem[{E and Ming(2007)}]{e_cauchyborn_2007}
E, W., Ming, P., Feb. 2007. Cauchy{\textendash}{Born} {Rule} and the
  {Stability} of {Crystalline} {Solids}: {Static} {Problems}. Archive for
  Rational Mechanics and Analysis 183~(2), 241--297.
\newline\urlprefix\url{http://link.springer.com/10.1007/s00205-006-0031-7}

\bibitem[{Ehrlacher et~al.(2016)Ehrlacher, Ortner, and
  Shapeev}]{ehrlacher_analysis_2016}
Ehrlacher, V., Ortner, C., Shapeev, A.~V., Dec. 2016. Analysis of {Boundary}
  {Conditions} for {Crystal} {Defect} {Atomistic} {Simulations}. Archive for
  Rational Mechanics and Analysis 222~(3), 1217--1268.
\newline\urlprefix\url{http://link.springer.com/10.1007/s00205-016-1019-6}

\bibitem[{Fellinger et~al.(2018)Fellinger, Tan, Hector, and
  Trinkle}]{fellinger_geometries_2018}
Fellinger, M.~R., Tan, A. M.~Z., Hector, L.~G., Trinkle, D.~R., Nov. 2018.
  Geometries of edge and mixed dislocations in bcc {Fe} from first-principles
  calculations. Physical Review Materials 2~(11).
\newline\urlprefix\url{https://link.aps.org/doi/10.1103/PhysRevMaterials.2.113605}

\bibitem[{Greengard and Rokhlin(1987)}]{greengard_fast_1987}
Greengard, L., Rokhlin, V., 1987. A {Fast} {Algorithm} for {Particle}
  {Simulations}. Journal of Computational Physics 73, 325--348.

\bibitem[{Hackbusch(1999)}]{hackbusch_sparse_1999}
Hackbusch, W., 1999. A {Sparse} {Matrix} {Arithmetic} based on {H}-{Matrices}.
  {Part} {I}: {Introduction} to {H}-{Matrices}. Computing 62, 89--108.

\bibitem[{Hodapp(2019)}]{hodapp_analysis_2019}
Hodapp, M., Dec. 2019. Analysis of a {Sinclair}-type domain decomposition
  solver for atomistic/continuum coupling. arXiv:1912.10530 [cond-mat]ArXiv:
  1912.10530.
\newline\urlprefix\url{http://arxiv.org/abs/1912.10530}

\bibitem[{Hodapp et~al.(2019)Hodapp, Anciaux, and Curtin}]{hodapp_lattice_2019}
Hodapp, M., Anciaux, G., Curtin, W., May 2019. Lattice {Green} function methods
  for atomistic/continuum coupling: {Theory} and data-sparse implementation.
  Computer Methods in Applied Mechanics and Engineering 348, 1039--1075.
\newline\urlprefix\url{https://linkinghub.elsevier.com/retrieve/pii/S0045782519300775}

\bibitem[{Kochmann and Venturini(2014)}]{kochmann_meshless_2014}
Kochmann, D.~M., Venturini, G.~N., Apr. 2014. A meshless quasicontinuum method
  based on local maximum-entropy interpolation. Modelling and Simulation in
  Materials Science and Engineering 22~(3), 034007.
\newline\urlprefix\url{http://stacks.iop.org/0965-0393/22/i=3/a=034007?key=crossref.1c71666a6a31fc5b99b546047cf19d95}

\bibitem[{Kohlhoff and Schmauder(1989)}]{kohlhoff_new_1989}
Kohlhoff, S., Schmauder, S., 1989. A {New} {Method} for {Coupled}
  {Elastic}-{Atomistic} {Modelling}. In: Atomistic {Simulation} of {Materials}.
  Springer US, Boston, MA, pp. 411--418.

\bibitem[{Li(2009)}]{li_efficient_2009}
Li, X., Sep. 2009. Efficient boundary conditions for molecular statics models
  of solids. Physical Review B 80~(10), 104112.
\newline\urlprefix\url{https://link.aps.org/doi/10.1103/PhysRevB.80.104112}

\bibitem[{Li(2012)}]{li_atomistic-based_2012}
Li, X., Jun. 2012. An atomistic-based boundary element method for the reduction
  of molecular statics models. Computer Methods in Applied Mechanics and
  Engineering 225-228, 1--13.
\newline\urlprefix\url{https://linkinghub.elsevier.com/retrieve/pii/S0045782512000825}

\bibitem[{Luskin and Ortner(2013)}]{luskin_atomistic--continuum_2013}
Luskin, M., Ortner, C., May 2013. Atomistic-to-continuum coupling. Acta
  Numerica 22, 397--508.
\newline\urlprefix\url{https://www.cambridge.org/core/product/identifier/S0962492913000068/type/journal_article}

\bibitem[{Martinsson(2002)}]{martinsson_fast_2002}
Martinsson, P.-G., 2002. Fast multiscale methods for lattice equations. {PhD}
  thesis, The University of Texas at Austin.

\bibitem[{Martinsson and Rodin(2009)}]{martinsson_boundary_2009}
Martinsson, P.-G., Rodin, G.~J., Aug. 2009. Boundary algebraic equations for
  lattice problems. Proceedings of the Royal Society A: Mathematical, Physical
  and Engineering Sciences 465~(2108), 2489--2503.
\newline\urlprefix\url{http://www.royalsocietypublishing.org/doi/10.1098/rspa.2008.0473}

\bibitem[{Morse(1929)}]{morse_diatomic_1929}
Morse, P.~M., Jul. 1929. Diatomic {Molecules} {According} to the {Wave}
  {Mechanics}. {II}. {Vibrational} {Levels}. Physical Review 34~(1), 57--64.
\newline\urlprefix\url{https://link.aps.org/doi/10.1103/PhysRev.34.57}

\bibitem[{Ortner and Zhang(2014)}]{ortner_energy-based_2014}
Ortner, C., Zhang, L., Sep. 2014. Energy-based atomistic-to-continuum coupling
  without ghost forces. Computer Methods in Applied Mechanics and Engineering
  279, 29--45.
\newline\urlprefix\url{https://linkinghub.elsevier.com/retrieve/pii/S0045782514002059}

\bibitem[{Parks et~al.(2008)Parks, Bochev, and Lehoucq}]{parks_connecting_2008}
Parks, M.~L., Bochev, P.~B., Lehoucq, R.~B., Jan. 2008. Connecting
  {Atomistic}-to-{Continuum} {Coupling} and {Domain} {Decomposition}.
  Multiscale Modeling \& Simulation 7~(1), 362--380.
\newline\urlprefix\url{http://epubs.siam.org/doi/10.1137/070682848}

\bibitem[{Pavia and Curtin(2015)}]{pavia_parallel_2015}
Pavia, F., Curtin, W.~A., Jul. 2015. Parallel algorithm for multiscale
  atomistic/continuum simulations using {LAMMPS}. Modelling and Simulation in
  Materials Science and Engineering 23~(5), 055002.
\newline\urlprefix\url{http://stacks.iop.org/0965-0393/23/i=5/a=055002?key=crossref.02ff7625c3b88ed47b97e830253f79d1}

\bibitem[{Quarteroni et~al.(2007)Quarteroni, Sacco, and
  Saleri}]{quarteroni_numerical_2007}
Quarteroni, A., Sacco, R., Saleri, F., 2007. Numerical mathematics, 2nd
  Edition. No.~37 in Texts in applied mathematics. Springer, Berlin ; New York.

\bibitem[{Rao et~al.(1998)Rao, Hernandez, Simmons, Parthasarathy, and
  Woodward}]{rao_greens_1998}
Rao, S., Hernandez, C., Simmons, J.~P., Parthasarathy, T.~A., Woodward, C.,
  Jan. 1998. Green's function boundary conditions in two-dimensional and
  three-dimensional atomistic simulations of dislocations. Philosophical
  Magazine A 77~(1), 231--256.
\newline\urlprefix\url{http://www.tandfonline.com/doi/abs/10.1080/01418619808214240}

\bibitem[{Rudd and Broughton(2000)}]{rudd_concurrent_2000}
Rudd, R.~E., Broughton, J.~Q., 2000. Concurrent {Coupling} of {Length} {Scales}
  in {Solid} {State} {Systems}. physica status solidi (b) 217, 41.

\bibitem[{Shapeev(2011)}]{shapeev_consistent_2011}
Shapeev, A.~V., Jul. 2011. Consistent {Energy}-{Based} {Atomistic}/{Continuum}
  {Coupling} for {Two}-{Body} {Potentials} in {One} and {Two} {Dimensions}.
  Multiscale Modeling \& Simulation 9~(3), 905--932.
\newline\urlprefix\url{http://epubs.siam.org/doi/10.1137/100792421}

\bibitem[{Shilkrot et~al.(2004)Shilkrot, Miller, and
  Curtin}]{shilkrot_multiscale_2004}
Shilkrot, L., Miller, R.~E., Curtin, W.~A., Apr. 2004. Multiscale plasticity
  modeling: coupled atomistics and discrete dislocation mechanics. Journal of
  the Mechanics and Physics of Solids 52~(4), 755--787.
\newline\urlprefix\url{https://linkinghub.elsevier.com/retrieve/pii/S0022509603001625}

\bibitem[{Shimokawa et~al.(2004)Shimokawa, Mortensen, Schi{\o}tz, and
  Jacobsen}]{shimokawa_matching_2004}
Shimokawa, T., Mortensen, J.~J., Schi{\o}tz, J., Jacobsen, K.~W., Jun. 2004.
  Matching conditions in the quasicontinuum method: {Removal} of the error
  introduced at the interface between the coarse-grained and fully atomistic
  region. Physical Review B 69~(21), 214104.
\newline\urlprefix\url{https://link.aps.org/doi/10.1103/PhysRevB.69.214104}

\bibitem[{Sinclair(1971)}]{sinclair_improved_1971}
Sinclair, J.~E., Dec. 1971. Improved {Atomistic} {Model} of a bcc {Dislocation}
  {Core}. Journal of Applied Physics 42~(13), 5321--5329.
\newline\urlprefix\url{http://aip.scitation.org/doi/10.1063/1.1659943}

\bibitem[{Sinclair(1975)}]{sinclair_influence_1975}
Sinclair, J.~E., Mar. 1975. The influence of the interatomic force law and of
  kinks on the propagation of brittle cracks. Philosophical Magazine 31~(3),
  647--671.
\newline\urlprefix\url{http://www.tandfonline.com/doi/abs/10.1080/14786437508226544}

\bibitem[{Sinclair et~al.(1978)Sinclair, Gehlen, Hoagland, and
  Hirth}]{sinclair_flexible_1978}
Sinclair, J.~E., Gehlen, P.~C., Hoagland, R.~G., Hirth, J.~P., Jul. 1978.
  Flexible boundary conditions and nonlinear geometric effects in atomic
  dislocation modeling. Journal of Applied Physics 49~(7), 3890--3897.
\newline\urlprefix\url{http://aip.scitation.org/doi/10.1063/1.325395}

\bibitem[{Tadmor et~al.(1996)Tadmor, Ortiz, and
  Phillips}]{tadmor_quasicontinuum_1996}
Tadmor, E.~B., Ortiz, M., Phillips, R., Jun. 1996. Quasicontinuum analysis of
  defects in solids. Philosophical Magazine A 73~(6), 1529--1563.
\newline\urlprefix\url{http://www.tandfonline.com/doi/abs/10.1080/01418619608243000}

\bibitem[{Thomson et~al.(1992)Thomson, Zhou, Carlsson, and
  Tewary}]{thomson_lattice_1992}
Thomson, R., Zhou, S.~J., Carlsson, A.~E., Tewary, V.~K., Nov. 1992. Lattice
  imperfections studied by use of lattice {Green}{\textquoteright}s functions.
  Physical Review B 46~(17), 10613--10622.
\newline\urlprefix\url{https://link.aps.org/doi/10.1103/PhysRevB.46.10613}

\bibitem[{Toselli and Widlund(2005)}]{toselli_domain_2005}
Toselli, A., Widlund, O.~B., 2005. Domain decomposition methods--algorithms and
  theory. No.~34 in Springer series in computational mathematics. Springer,
  Berlin, oCLC: ocm56879011.

\bibitem[{Tyrtyshnikov(1996)}]{tyrtyshnikov_mosaic-skeleton_1996}
Tyrtyshnikov, E., Jun. 1996. Mosaic-{Skeleton} approximations. Calcolo
  33~(1-2), 47--57.
\newline\urlprefix\url{http://link.springer.com/10.1007/BF02575706}

\bibitem[{Woodward and Rao(2002)}]{woodward_flexible_2002}
Woodward, C., Rao, S.~I., May 2002. Flexible \textit{{Ab} {Initio}} {Boundary}
  {Conditions}: {Simulating} {Isolated} {Dislocations} in bcc {Mo} and {Ta}.
  Physical Review Letters 88~(21).
\newline\urlprefix\url{https://link.aps.org/doi/10.1103/PhysRevLett.88.216402}

\bibitem[{Zhang(2005)}]{zhang_schur_2005}
Zhang, F. (Ed.), 2005. The {Schur} complement and its applications. No. v. 4 in
  Numerical methods and algorithms. Springer, New York.

\end{thebibliography}

\end{document}